\documentclass[10pt, reqno]{amsart}

\usepackage{hyperref} 
\usepackage{amssymb, mathrsfs}
\usepackage{amsmath, amsthm}
\usepackage{enumerate}
\usepackage{dsfont}
\usepackage{color}
\usepackage[a4paper]{geometry}
\usepackage[utf8]{inputenc}   
\usepackage{stmaryrd}
\usepackage{titlesec, wasysym}
\usepackage{appendix}
\usepackage{todonotes, fancyhdr}
\usepackage{chngcntr}
\usepackage{cancel}

\usepackage{tikz-cd}
\usetikzlibrary{calc}

\usepackage{setspace}
\renewcommand{\baselinestretch}{0.99}

\usepackage[all]{xy}
\numberwithin{subsection}{section}
\numberwithin{subsubsection}{subsection}
\numberwithin{equation}{section} 

\newenvironment{Dem}[1][\unskip]{%
    \begin{list}{\hspace{1.02cm}{\sf \textbf{ Proof #1 --}}}{%
        \setlength{\topsep}{0pt}%
        \setlength{\leftmargin}{0pt}%
        \setlength{\rightmargin}{0pt}%
        \setlength{\listparindent}{0pt}%
        \setlength{\itemindent}{0pt}%
        \setlength{\parsep}{0pt}%
        \addtolength{\leftmargin}{0pt}
        \addtolength{\rightmargin}{0pt}%
    } \item }{\hfill $\rhd$\end{list}\smallskip}

\newenvironment{Dem*}[1][\unskip]{%
    \begin{list}{\hspace{0cm}{\sf \textbf{{\small Proof #1 --}}}}{%
        \setlength{\topsep}{0pt}%
        \setlength{\leftmargin}{0pt}%
        \setlength{\rightmargin}{0pt}%
        \setlength{\listparindent}{0pt}%
        \setlength{\itemindent}{0pt}%
        \setlength{\parsep}{0pt}%
        \addtolength{\leftmargin}{20pt}%
        \addtolength{\rightmargin}{0pt}%
    } \item }{\hfill $\rhd$\end{list}\smallskip}

\renewcommand\thesection       {\arabic{section}}
\renewcommand\thesubsection    {\thesection{\boldmath $.$}\arabic{subsection}}
\renewcommand\thesubsubsection    {\thesection{\boldmath $.$}\arabic{subsection}{\boldmath $.$}\arabic{subsubsection}}

\titleformat{\section}[block]
{\filcenter\normalfont\sffamily\bfseries\Large}
{{\hspace{-0.7cm}}\thesection \hspace{0.2em} --\vspace{0.3cm}}{0.5em}{}

\titleformat{\subsection}[block]
{\filcenter\normalfont\sffamily\bfseries\large}  						  
{\hspace{-0.7cm}\thesubsection \hspace{0.5em} \vspace{0.3cm}}{.5em}{}  
\titlespacing{\subsection}{-0pc}{1.5ex plus .1ex minus .2ex}{0pc}

\titleformat{\subsubsection}[block]
{\normalfont\sffamily\bfseries}					  
{\thesubsubsection \vspace{0.3cm}}{.5em}{}  
\titlespacing{\subsection}{-0pc}{1.5ex plus .1ex minus .2ex}{0pc}

\newtheoremstyle{mystyle}
{3pt}               
{3pt}               
{\it }                      
{}                      
{\sffamily\bfseries}             
{}                      
{0.5em}                 
{#1 #2{\hspace{0.2cm}--\hspace{-0.2cm}}  }   

\theoremstyle{mystyle}

\newtheorem{thm}{Theorem}
\newtheorem*{thm*}{Theorem}

\newtheorem{cor}[thm]{\hspace{-0.15cm}  {Corollary} }
\newtheorem{lem}[thm]{\hspace{-0.14cm}  {Lemma} }
\newtheorem{prop}[thm]{\hspace{-0.13cm} {Proposition}}

\newtheoremstyle{mystyle2}
{3pt}               
{3pt}               
{\it }                      
{}                      
{\sffamily\bfseries}             
{}                      
{0.5em}                 
{\llap{#2 }#1{\hspace{0.2cm}--}}

\theoremstyle{mystyle2}

\newtheorem*{definition*}{Definition}
\newtheorem*{theorem*}{Theorem}
\newtheorem*{Remark*}{Remark}
\newtheorem*{lem*} {Lemma}
\newtheorem*{defn*} {Definition}
\newtheorem*{prop*} {Proposition}
\newtheorem*{cor*} {Corollary}

\newcommand{\ssk}{\smallskip}

\renewcommand{\epsilon}{\varepsilon}
\newcommand{\eps}{\epsilon}

\newcommand\bbE{\mathbb{E}}

\newcommand\bbN{\mathbb{N}}
\newcommand\bbP{\mathbb{P}}
\newcommand\bbR{\mathbb{R}}
\newcommand{\bbT}{\mathbb{T}}

\newcommand{\mcB}{\mathcal{B}}

\newcommand{\mcE}{\mathcal{E}}
\newcommand{\mcF}{\mathcal{F}}

\newcommand{\mcI}{\mathcal{I}}

\newcommand{\mcM}{\mathcal{M}}

\newcommand{\scrR}{\ensuremath{\mathscr{R}}}

\makeatletter
\newcommand*{\defeq}{\mathrel{\rlap{%
                     \raisebox{0.25ex}{$\m@th\cdot$}}%
                     \raisebox{-0.25ex}{$\m@th\cdot$}}%
                     =}
\makeatother

\makeatletter
\newcommand*{\eqdef}{=\mathrel{\rlap{%
                     \raisebox{0.25ex}{$\m@th\cdot$}}%
                     \raisebox{-0.25ex}{$\m@th\cdot$}}%
                     }
\makeatother

\usepackage{tikz-cd}
\usetikzlibrary{calc}
\usetikzlibrary{shapes.misc}
\usetikzlibrary{shapes.symbols}
\usetikzlibrary{snakes}
\usetikzlibrary{decorations}
\usetikzlibrary{decorations.markings}

\makeatletter
\pgfdeclareshape{crosscircle}
{
  \inheritsavedanchors[from=circle] 
  \inheritanchorborder[from=circle]
  \inheritanchor[from=circle]{north}
  \inheritanchor[from=circle]{north west}
  \inheritanchor[from=circle]{north east}
  \inheritanchor[from=circle]{center}
  \inheritanchor[from=circle]{west}
  \inheritanchor[from=circle]{east}
  \inheritanchor[from=circle]{mid}
  \inheritanchor[from=circle]{mid west}
  \inheritanchor[from=circle]{mid east}
  \inheritanchor[from=circle]{base}
  \inheritanchor[from=circle]{base west}
  \inheritanchor[from=circle]{base east}
  \inheritanchor[from=circle]{south}
  \inheritanchor[from=circle]{south west}
  \inheritanchor[from=circle]{south east}
  \inheritbackgroundpath[from=circle]
  \foregroundpath{
    \centerpoint%
    \pgf@xc=\pgf@x%
    \pgf@yc=\pgf@y%
    \pgfutil@tempdima=\radius%
    \pgfmathsetlength{\pgf@xb}{\pgfkeysvalueof{/pgf/outer xsep}}%
    \pgfmathsetlength{\pgf@yb}{\pgfkeysvalueof{/pgf/outer ysep}}%
    \ifdim\pgf@xb<\pgf@yb%
      \advance\pgfutil@tempdima by-\pgf@yb%
    \else%
      \advance\pgfutil@tempdima by-\pgf@xb%
    \fi%
    \pgfpathmoveto{\pgfpointadd{\pgfqpoint{\pgf@xc}{\pgf@yc}}{\pgfqpoint{-0.707107\pgfutil@tempdima}{-0.707107\pgfutil@tempdima}}}
    \pgfpathlineto{\pgfpointadd{\pgfqpoint{\pgf@xc}{\pgf@yc}}{\pgfqpoint{0.707107\pgfutil@tempdima}{0.707107\pgfutil@tempdima}}}
    \pgfpathmoveto{\pgfpointadd{\pgfqpoint{\pgf@xc}{\pgf@yc}}{\pgfqpoint{-0.707107\pgfutil@tempdima}{0.707107\pgfutil@tempdima}}}
    \pgfpathlineto{\pgfpointadd{\pgfqpoint{\pgf@xc}{\pgf@yc}}{\pgfqpoint{0.707107\pgfutil@tempdima}{-0.707107\pgfutil@tempdima}}}
  }
}
\makeatother

\colorlet{symbols}{black}    
\colorlet{testcolor}{green!60!black}
\colorlet{supcolor}{red!60!black}

\tikzset{
	root/.style={circle, fill=testcolor!70, draw=testcolor, inner sep=1pt, minimum size=0.5mm},
		kepsus/.style={semithick, ->},
	dot/.style={circle, draw=black, fill=black, inner sep=0pt, minimum size=0.3mm},
	noise/.style={circle, draw=black, fill=white, inner sep=0pt, minimum size=1mm},
	h/.style={circle, draw=black, fill=black, inner sep=0pt, minimum size=0.3mm},	
	var/.style={circle, fill=white, draw=purple, fill=purple, inner sep=0pt, minimum size=0.6mm},
	bdot/.style={circle, draw=black, fill=white, inner sep=0pt, minimum size=1mm},
	bluedot/.style={circle,draw=blue, fill=blue, inner sep=0pt, minimum size=2mm},
	noiseblue/.style={circle, fill=blue!20, draw=blue, inner sep=0pt, minimum size=1mm},
	dtestfcn/.style={ultra thick, densely dashed, testcolor,shorten >=1pt,shorten <=1pt,<-},
	testfcnx/.style={semithick,testcolor,shorten >=1pt,shorten <=1pt,<-, postaction={decorate,decoration={markings,mark=at position 0.6 with {\drawx}}}},
	testfcn/.style={semithick, testcolor, shorten >=1pt,shorten <=1pt,-},	
	K/.style= {semithick, shorten >=0pt,shorten <=0pt,-},
	kdashed/.style= {semithick,densely dashed,shorten >=1pt,shorten <=1pt,<->},
	DK/.style={thick, densely dotted, shorten >=0pt,shorten <=0pt},   
	ksup/.style={thick, supcolor, shorten >=1pt,shorten <=1pt},
	dots/.style={semithick,dotted,shorten >=1pt,shorten <=1pt},
	Deps/.style={semithick,draw=black!25,fill=black!25,shorten >=1pt,shorten <=1pt,->},
	kbase/.style={semithick,dotted,shorten >=1pt,shorten <=1pt,->},
	multx/.style={shorten >=1pt,shorten <=1pt,
		postaction={decorate,decoration={markings,mark=at position 0.5 with {\drawx}}}},
	kernelx/.style={semithick,shorten >=1pt,shorten <=1pt,->,
		postaction={decorate,decoration={markings,mark=at position 0.4 with {\drawx}}}},
	kernel1/.style={->,semithick,shorten >=1pt,shorten <=1pt,postaction={decorate,decoration={markings,mark=at position 0.45 with {\draw[--] (0,-0.1) -- (0,0.1);}}}},
	kernel2/.style={->,semithick,shorten >=1pt,shorten <=1pt,postaction={decorate,decoration={markings,mark=at position 0.45 with {\draw[--] (0.05,-0.1) -- (0.05,0.1);\draw[--] (-0.05,-0.1) -- (-0.05,0.1);}}}},
	kernelBig/.style={semithick,shorten >=1pt,shorten <=1pt,decorate, decoration={zigzag,amplitude=1.5pt,segment length = 3pt,pre length=2pt,post length=2pt}},
	rho/.style={dotted,semithick,shorten >=1pt,shorten <=1pt},
	renorm/.style={shape=circle,fill=white,inner sep=1pt},
	labl/.style={shape=rectangle,fill=white,inner sep=1pt},
	xi/.style={circle,fill=symbols!10,draw=symbols,inner sep=0pt,minimum size=1.2mm},
	xix/.style={crosscircle,fill=symbols!10,draw=symbols,inner sep=0pt,minimum size=1.2mm},
	xib/.style={circle,fill=symbols!10,draw=symbols,inner sep=0pt,minimum size=1.6mm},
	xibx/.style={crosscircle,fill=symbols!10,draw=symbols,inner sep=0pt,minimum size=1.6mm},
	not/.style={circle,fill=symbols,draw=symbols,inner sep=0pt,minimum size=0.5mm},
	>=stealth,
	graydot/.style={circle,fill=gray,inner sep=0pt, minimum size=1mm},
	zero/.style={circle,inner sep=0pt, minimum size=1mm, draw},
	kernelprimeeps/.style={densely dashed, semithick,shorten >=1pt,shorten <=1pt},
	smalldot/.style={circle,fill=black,draw=black, solid,inner sep=0pt,minimum size=0.5mm},
	}

\begin{document}

\begin{center}
{\LARGE\sffamily{ Random models for singular SPDEs  \vspace{0.5cm}}}
\end{center}

\begin{center}
{\sf I. BAILLEUL and Y. BRUNED}
\end{center}

\vspace{1cm}

\begin{center}
\begin{minipage}{0.8\textwidth}
\renewcommand\baselinestretch{0.7} \scriptsize \textbf{\textsf{\noindent Abstract.}} We give a proof of the convergence of the BHZ renormalized model associated with the generalized (KPZ) equation that does not require the full strength of the BPHZ renormalisation. Our approach is based on a convenient form of chaos decomposition. The other key ingredient is a generalisation of the Hairer-Quastel convergence theorem for Feynman diagrams with certain decorations encoding Taylor remainders.  With these ideas we are able to construct the model for the generalised KPZ equation.
\end{minipage}
\end{center}

\vspace{0.6cm}

\section{Introduction}
\label{SectionIntro}

\S 1. {\it Singular stochastic partial differential equations.} Let $\zeta$ stand for a spacetime distribution of negative regularity. Random or not we talk about it as a `noise'. We call {\sl singular}, an elliptic or parabolic semilinear partial differential equation (PDE) with real-valued unknown $u$, an equation of the form
\begin{equation} \label{EqModelEq}
\mathscr{L} u = f(u,\nabla u ; \zeta)
\end{equation}
where a formal application of Schauder estimates does not give sufficient regularity to $u$ for the term $f(u,\nabla u ; \zeta)$ to make sense. Recall indeed that loosely speaking the product of a distribution of regularity $a$ with a function of regularity $b$ is well-defined if and only if $a+b>0$. A typical example of a singular PDE is provided by the two dimensional parabolic Anderson model equation
\begin{equation*}
(\partial_t-\Delta) u = u \zeta,
\end{equation*}
where $\zeta$ stands for the random realization of a space white noise on the two dimensional torus. The distribution $\zeta$ has almost surely H\"older regularity $-1-\epsilon$ for all $\epsilon>0$; this gives $u$ an a priori H\"older regularity $1-\epsilon$ which is not sufficient to make sense of the product $u\zeta$. Coupled systems of singular PDEs are characterized by the same product problem. Such equations/systems typically appear as the probabilistic large scale limit of random microscopic equations/systems where the strength of the nonlinearity and the randomness balance each other. 

\ssk

\S 2. {\it Regularity structures.} The study of singular stochastic PDEs was launched by the groundbreaking works of M. Hairer on regularity structures \cite{Hai14} and Gubinelli, Imkeller \& Perkowski on paracontrolled calculus \cite{GIP}. While \cite{Hai14} and \cite{GIP} use different tools to build an equation-dependent setting where to make sense of it and provide a robust solution theory for a large family of equations, both works share the same guidelines.

\begin{enumerate}
	\item[1.] Work in a restricted space of potential solutions. The functions from this space are characterized by the fact that they `look like' some linear combination of reference objects that depend only on the noise $\zeta$. The choice of language used to give sense to the expression `look like' distinguishes regularity structures from paracontrolled calculus.   \vspace{0.1cm}

	\item[2.] Point 1 brings back the initial problem of making sense of the ill-defined products involving $u$ in \eqref{EqModelEq} to the problem of giving sense to some ill-defined products of quantities that only involve the noise $\zeta$, typically some polynomial expression $P(\zeta,\dots,\zeta)$ of $\zeta$. This is where working with {\sl random} distributions $\zeta$ saves the day. While the random variable $\omega\mapsto P\big(\zeta(\omega),\dots,\zeta(\omega)\big)$ does not make sense one can construct in a consistent way a {\sl random variable}, suggestively denoted here $\omega\mapsto P(\zeta,\dots,\zeta)(\omega)$, that plays the role of $P\big(\zeta(\omega),\dots,\zeta(\omega)\big)$. This is what {\sl renormalization} is about and the collection $\widehat\zeta(\omega)$ of all the random variables $P(\zeta,\dots,\zeta)(\omega)$ needed in the analysis of a given equation is called the {\sl enhanced noise}. On a technical level the above polynomial expressions of the noise come under the form of {\sl models over a regularity structure}.   \vspace{0.1cm}
	
	\item[3.] Replacing the formal quantities $P\big(\zeta(\omega),\dots,\zeta(\omega)\big)$ by the quantities $P(\zeta,\dots,\zeta)(\omega)$ in a naive formulation of the equation turns out to provide a well-defined equation that can be solved uniquely in an appropriate function space. We talk of a `renormalized equation' and its solution.   \vspace{0.1cm}

	\item[4.]  Of course a proper interpretation of the solution of the renormalized equation is needed. In the end one can prove that it is the limit, in a probabilistic sense, of solutions to equations of the form 
	$$
	\mathscr{L} u^\epsilon = f(u^\epsilon,\nabla u^\epsilon ; \zeta^\epsilon) + g^\epsilon(u^\epsilon)
	$$
	in which $\zeta$ has been regularized into $\zeta^\epsilon$ and some $\epsilon$-dependent, a priori diverging, counterterm $g^\epsilon(\cdot)$ has been added into the dynamics. The parameter $\epsilon>0$ goes to $0$ in the limit.
\end{enumerate}

\ssk

The analytic machinery introduced by Hairer in \cite{Hai14} was complemented in subsequent joint works with Bruned \& Zambotti \cite{BHZ}, Chandra \cite{CH} and Bruned, Chandra \& Chevyrev \cite{BCCH} to provide a robust solution theory for a large class of (systems of) equations. The work \cite{BHZ} provided a systematic study of the algebraic setting needed to analyse singular stochastic PDEs from a regularity structures point of view. The work \cite{CH} proved the convergence of a general algorithm used to construct the random variables $P(\zeta,\dots,\zeta)(\omega)$ as limits in some appropriate spaces of quantities of the form 
$$
P\big(\zeta^\epsilon(\omega),\dots,\zeta^\epsilon(\omega)\big) - \sum_Q c_Q^\epsilon \, Q\big(\zeta^\epsilon(\omega),\dots,\zeta^\epsilon(\omega)\big)
$$
for some finite sums of polynomial expressions $Q$ of $\zeta^\epsilon$ and $\epsilon$-dependent, a priori diverging, constants $c_Q^\epsilon$. We call {\sl BHZ renormalization} this construction, after the initials of the authors of the work \cite{BHZ}. (Amusingly these initials are very close to the initials BPHZ used to name the renormalization algorithm of Feynman graphs introduced by Bogoliubov \& Parasiuk in the 50s, and clarified later by some works of Hepp and Zimmermann. The BHZ renormalization algorithm is an elaboration of the BPHZ renormalization prescription on different types of objects. Indeed, the BPHZ renormalisation corresponds to the inductive extraction-contraction of subdiagrams of negative degree in a Feynman diagram whereas the BHZ renormalisation performs extraction-contraction of subtrees with negative degree in a given decorated tree. The renormalised decorated trees are then mapped to stochastic iterated integrals and with the computation of its variance one obtains a family of renormalised Feynman diagrams in the sense of BPHZ.) The work \cite{BCCH} proved point 4 above. In the end for any given equation of a well-identified class of subcritical singular elliptic/parabolic semilinear stochastic PDEs there exists, under a small parameter condition, a unique solution to that equation in an appropriate equation-dependent function space. This solution depends continuously on the realization $\widehat\zeta(\omega)$ of the random enhanced noise. To contrast things, note that $\widehat\zeta(\omega)$ itself is only a measurable function of $\zeta(\omega)$. We refer the reader to the overviews \cite{ChandraWeber, CorwinShen} of Chandra \& Weber and Corwin \& Shen for non-technical introductions to the domain of semilinear singular SPDEs, to the books \cite{FrizHairer, Berglund} of Friz \& Hairer and Berglund for a mildly technical introduction to regularity structures, and to Bailleul \& Hoshino's Tourist's Guide \cite{RSGuide} for a thorough tour of the analytic and algebraic sides of the theory. Our previous work \cite{BailleulBruned} gives a short proof of a generalization of point 4, and some aspects of the renormalization problem are described in Hairer's lecture notes \cite{HairerTakagi}.

\ssk

\S 3. {\it Renormalization.} The fantastic work \cite{CH} of Chandra \& Hairer on the convergence of the BHZ renormalized models is fairly difficult. It holds under a very general set of conditions on the law of the translation invariant (possibly multi-dimensional) noise that involve its cumulants. Recently Linares, Otto, Tempelmayr \& Tsatsoulis \cite{LOTT} devised a different approach to the renormalization problem within the setting devised by Otto, Sauer, Smith \& Weber \cite{OSSW2} for the study of a class of quasilinear singular stochastic PDEs. While the latter shares many features with the setting of regularity structures devised for the study of semilinear singular equations, the setting introduced in \cite{OSSW2} is different; it involves in particular no tree indexed objects even though its construction is inductive. More importantly for us they prove their convergence result for the renormalized regularized models assuming a spectral gap inequality for the law of the random noise rather than making assumptions on its cumulants. 

\ssk

Denote by ${\sf M}^\epsilon = ({\sf \Pi}^\epsilon, {\sf g}^\epsilon)$ the naive admissible model built from the regularized noise $\zeta^\epsilon$. For any symbol $\tau$ in the regularity structure the quantity ${\sf \Pi}^\epsilon(\tau)$ is a function that is a polynomial functional of the regularized noise, built in a canonical way from the combinatorial structure of $\tau$ and the equation under consideration. This analytical object typically diverges as $\epsilon>0$ is going to $0$. The map ${\sf g}^\epsilon$ associates to any symbol $\sigma$ in the regularity structure of positive `degree' a function on the sate space. All the objects will be properly introduced in Section \ref{SubsectionRS}. Denote by $\overline{\sf M}^\epsilon = (\overline{\sf \Pi}^\epsilon,\overline{\sf g}^\epsilon)$ the BHZ renormalized naive admissible model built from $\zeta^\epsilon$. On a technical side the convergence of $\overline{\sf \Pi}^\epsilon$ is encoded in inequalities of the form
\begin{equation} \label{EqMainEstimate}
\big\Vert \big\langle \overline{\sf \Pi}_{z}^\epsilon \tau - \overline{\sf \Pi}_{z}^{\epsilon'}\hspace{-0.07cm}\tau \,,\,\varphi_{z}^\lambda \big\rangle \big\Vert_{L^2(\Omega)} \lesssim o_{\epsilon\vee\epsilon'}(1)\,\lambda^{\vert\tau\vert}, 
\end{equation}
for all symbols $\tau$ of the regularity structure with negative degree $\vert\tau\vert$. These inequalities involve the recentered version $\overline{\sf \Pi}_{z}^{\epsilon}$ of the renormalized interpretation map $\overline{\sf \Pi}^{\epsilon}$. Think here of $\varphi_{z}^\lambda$ as a regular $\lambda$-approximation of a Dirac distribution at the point $z$ -- proper notations will be introduced below. Linares, Otto, Tempelmayr \& Tsatsoulis take profit from their spectral gap assumption to bound the $L^2(\Omega)$ norm of the quantity of interest
$$
\big\Vert \big\langle \overline{\sf \Pi}_{z}^\epsilon \tau - \overline{\sf \Pi}_{z}^{\epsilon'}\hspace{-0.07cm}\tau \,,\,\varphi_z^\lambda \big\rangle \big\Vert_{L^2(\Omega)} \leq \big\vert \bbE\big[\big\langle \overline{\sf \Pi}_{z}^\epsilon \tau - \overline{\sf \Pi}_{ z}^{\epsilon'}\hspace{-0.07cm}\tau \,,\,\varphi_{ z}^\lambda \big\rangle\big]\big\vert + \big\Vert \big\langle d\big(\overline{\sf \Pi}_{ z}^\epsilon \tau - \overline{\sf \Pi}_{z}^{\epsilon'}\hspace{-0.07cm}\tau\big)\,,\,\varphi_{z}^\lambda \big\rangle \big\Vert_{L^2(\Omega)}
$$
by its expectation and the operator norm of its Malliavin derivative. In the setting of \cite{LOTT} the control of the above expectation comes at low cost from their choice of renormalization procedure. The control of the $L^2$-norm of the derivative is way more involved and done inductively by seeing first $d\overline{\sf \Pi}^\epsilon$ as part of an extended model. The very recent work \cite{OST} of Otto, Seong and Tempelmayr gives a reader's guide digest of \cite{LOTT}; see also Broux, Otto \& Tempelmayr's recent lecture notes \cite{BOT24}. Hairer \& Steele implemented in \cite{HS23} Otto \& co's spectral gap strategy in the classical tree-based setting of regularity structures. A different, and somewhat more general, implementation was given in Bailleul \& Hoshino's work \cite{BH23}.

\ssk

\S4. {\it Our strategy for the (gKPZ) equation.} Let $\bbT$ stand for the one dimension torus and $\zeta$ stand for a spacetime white noise on $\bbR\times{\bbT}$. It has almost surely parabolic H\"older regularity $-3/2-\kappa$ for every $\kappa>0$. We fix $\kappa$ once and for all. The generalized (KPZ) equation has as unknown a function $u : [0,T]\times {\bbT}\rightarrow\bbR$ and reads
\begin{equation} \label{EqGKPZ}
\big(\partial_t - \partial_x^2\big)u = f(u)\zeta + g(u)(\partial_x u)^2,
\end{equation}
for some $C^4$ bounded real valued functions $f,g$ with bounded derivatives. The time horizon $T$ is not fixed and is part of the problem. As we expect from Schauder estimates that $u$ has parabolic regularity $1/2-\kappa$ the products $f(u)\zeta$ and $(\partial_x u)^2$ are not well-defined, let alone $g(u)(\partial_x u)^2$. The equation is singular.

We follow in the present work the elementary road of chaos decomposition in the form of {\it Stroock's formula}. For a symbol $\tau$ built from $\vert\tau\vert_\zeta$ noise symbols one has
\begin{equation} \label{EqStroock}
\big\Vert \big\langle \overline{\sf \Pi}_{z}^\epsilon\tau - \overline{\sf \Pi}_{z}^{\epsilon'}\hspace{-0.07cm}\tau \,,\,\varphi_{z}^\lambda \big\rangle \big\Vert_{L^2(\Omega)} \leq \sum_{n=0}^{\vert\tau\vert_\zeta-1} \big\Vert \bbE\big[\langle d^n(\overline{\sf \Pi}_{z}^\epsilon\tau - \overline{\sf \Pi}_{z}^{\epsilon'}\hspace{-0.07cm}\tau) , \varphi_{z}^\lambda \rangle\big] \big\Vert + \big\Vert \langle d^{\vert\tau\vert_\zeta}(\overline{\sf \Pi}_{z}^\epsilon\tau - \overline{\sf \Pi}_{z}^{\epsilon'}\hspace{-0.07cm}\tau) , \varphi_{z}^\lambda \rangle \big\Vert.
\end{equation}
The norms $\Vert\cdot\Vert$ are Hilbert-Schmidt type norms. We take profit here of the fact that the trees in the regularity structure of the (gKPZ) equation have at most four noises. Only trees with an even number of noises have a non-null expectation. This is where working with trees with at most four noises helps and reduces the work to considering a small number of situations despite the fact that the regularity structure of the (gKPZ) equation has $42$ trees of negative degree. Note that the full, $\vert\tau\vert_\zeta$-th order, derivative term in the right hand side is not random anymore and
$$
d^{\vert\tau\vert_\zeta}(\overline{\sf \Pi}_{ z}^\epsilon \tau - \overline{\sf \Pi}_{ z}^{\epsilon'}\hspace{-0.07cm}\tau) = d^{\vert\tau\vert_\zeta}({\sf \Pi}_{ z}^\epsilon\tau - {\sf \Pi}_{ z}^{\epsilon'}\hspace{-0.07cm}\tau).
$$
The quantities $\Vert \bbE[\langle d^k(\overline{\sf \Pi}_{z}^\epsilon\tau - \overline{\sf \Pi}_{z}^{\epsilon'}\hspace{-0.07cm}\tau) , \varphi_{z}^\lambda \rangle] \Vert^2$ are given by some (sums of) iterated integrals whose structure is encoded by a graph. Nodes of a graph represent integration variables. Edges of a graph represent some Taylor remainder of an edge-dependent function of its two node variables and another node of the graph around which one does the Taylor expansion. Such graphs are called Feynman graphs. Hairer \& Quastel considered in \cite{HairerQuastel} graphs where the Taylor expansions are done with respect to the same point for all edges. They found some conditions under which a graph has a bound of the form $o_{\varepsilon\vee\varepsilon'}(1) \lambda^{\vert\tau\vert}$. However some of the graphs involved in the study of the (gKPZ) equation fall outside of the scope of Hairer \& Quastel's criterion. This is for example the case of some graph associated with the tree $\mcI_1(\Xi\mcI(\Xi))$, with the notations of Section \ref{SubsectionRS}. Bruned \& Nadeem devised recently \cite{BN1} a variant of Hairer \& Quastel's criterion where each edge function can be expanded at an edge-dependent point. This seemingly minor difference offers a real flexibility. As in \cite{HairerQuastel}, Bruned \& Nadeem have two sets of conditions: The first set of conditions ensures the convergence of the iterated integrals as $\varepsilon>0$ goes to $0$, the second set of conditions ensures the $\lambda^{\vert\tau\vert}$ behaviour of the integral. For a graph with $n$ vertices each set of conditions involves around $2^n$ conditions. Now here is the point: Unlike in the setting of Hairer \& Quastel, changing the base point of a Taylor expansion to a nearby point, or changing the order of the expansion, changes only a very small number of expressions in the set of conditions. So one can check that some particular changes of that local type preserve the two sets of conditions. We talk of elementary changes. Rather than checking directly that a given graph like $\Vert \bbE[\langle d^k(\overline{\sf \Pi}_{z}^\epsilon\tau - \overline{\sf \Pi}_{z}^{\epsilon'}\hspace{-0.07cm}\tau) , \varphi_{z}^\lambda \rangle] \Vert$ satisfies the two sets of conditions we can then start from a graph that satisfies one of the two sets of conditions, by construction, and arrive at the targeted graph after finitely many elementary changes.

\ssk

Let us mention that our method applies to other equations such as two and three dimensional stochastic Yang-Mills dynamics as we are able to cover diagrams that contain only one subdivergence which is the case for these equations. If one has to face nasty or overlapping subdivergences then one has to extend our approach. It seems achievable but far from being trivial. But it could be a simplification over the work  \cite{CH}.

\medskip

\textbf{Organisation of the work.} We describe in Section \ref{SectionGKPZ} the generalized (KPZ) equation, its associated regularity structure and the BHZ renormalized model. We prove its convergence in Section \ref{SectionConvergence}. The useful notion of mirror graph is introduced in Section \ref{SubsectionMirror} and we talk in Section \ref{SubsectionMalliavin} of  Malliavin derivatives of trees and their expectation as this is what appears in Stroock's formula \eqref{EqStroock}. The proof of convergence of the BHZ model to a limit model is done in Section \ref{SubsectionProofConvergence}; its main ingredient is recalled in Section \ref{SubsectionWorkhorse}.

\medskip

{\textbf{Notations.}} {\it We denote by $z=(t,x)$ a generic spacetime point, with parabolic norm 
$$
\Vert z\Vert_{\frak{s}} \defeq \sqrt{\vert t\vert + \vert x\vert^2}. 
$$
For $k\in\bbN^2$ set
\begin{equation*}
\vert k\vert_{\frak{s}} = \vert (k_0,k_1)\vert_{\frak{s}} \defeq 2 k_0 + k_1.
\end{equation*} 
We will freely use all the notions and tools of regularity structures. We recall here Hairer's notation
$$
|\!|\!| K |\!|\!|_{a,m} \defeq \sup_{\vert k\vert_{\frak{s}}\leq m} \sup_z \Vert z\Vert_{\frak{s}}^{\vert k\vert_{\frak{s}}+a} \big\vert \partial^k K(z)\big\vert
$$
for the size of a kernel with a singularity at $0$. The heat kernel $P$ has a finite $|\!|\!| \cdot |\!|\!|_{1,m}$ size for any fixed $m\geq 0$. Pick a smooth non-increasing function $\chi : [0,\infty)\rightarrow [0,1]$ equal to $1$ on $[0,1/2]$ and $0$ on $[1,\infty)$. Rather than working with the heat kernel $P$ we work with 
$$
K\defeq(P-c)\chi,
$$ 
with the constant $c$ chosen for $K$ to have null mean -- a useful technical property that does not destroy the behaviour of $K$ as $z$ approaches the spacetime point $0$. The corresponding two convolutions operators coincide up to a smoothing operator and $|\!|\!| K |\!|\!|_{1,m}<\infty$.

Denote by $H$ the Hilbert space $L^2([0,1]\times \bbT)$. All random variables will be defined on a fixed, rich enough, probability space $(\Omega,\mcF,\bbP)$.}

\subsection*{Acknowledgements}

{\small
	 Y. B.  gratefully acknowledges funding support from the European Research Council (ERC) through the ERC Starting Grant Low Regularity Dynamics via Decorated Trees (LoRDeT), grant agreement No.\ 101075208. Views and opinions
	expressed are however those of the author(s) only and do not necessarily reflect those of
	the European Union or the European Research Council Executive Agency. Neither the
	European Union nor the granting authority can be held responsible for them.
}   

\subsection*{Author Contributions} All authors contributed equally to every step of the production.

\subsection*{Data Availibility} Data sharing not applicable to this article as no datasets were
generated or analysed during the current study.

\subsection*{Conflict of interest} The authors declare that they have no conflict of interest.

\bigskip

\section{The generalized (KPZ) equation}
\label{SectionGKPZ}

This section is dedicated to setting the stage for the proof of convergence of the BHZ renormalized model given in Section \ref{SectionConvergence}. We describe in Section \ref{SubsectionRS} the regularity structure associated with the generalized (KPZ) equation and the BHZ renormalized model. We introduce the useful notions of tree-like and mirror graphs in Section \ref{SubsectionMirror}. Malliavin derivatives of trees and their expectations are the object of Section \ref{SubsectionMalliavin}.

\medskip

\subsection{The (gKPZ) regularity structure and the BHZ renormalized model}
\label{SubsectionRS}

We describe in this section the regularity structure associated with the generalized (KPZ) equation and the BHZ renormalized model. A reader familiar with the theory can skip this section.

\medskip

\noindent \textit{{\S1. The (gKPZ) regularity structure.}} Let $\mcI$ and $\mcI_{(0,1)}$ stand for abstract `integration' symbols; they are formal placeholders for the convolution operator with $K$ and its space derivative. More generally, for $a\in\bbN^2$, we denote by $\mcI_a$ some abstract integration symbols encoding the convolution operator with $ \partial^a K $. We define also an abstract derivative $ \mathcal{D}^a $ satisfying the Leibniz rule and such that
\begin{equation} \label{derivation_symbol}
	\mathcal{D}^a \mathcal{I}_b(\tau) = \mathcal{I}_{a+b}(\tau), \quad \mathcal{D}^{a} X^k = X^{k-a}.
\end{equation}
This derivative corresponds to the differential operator $ \partial^a $. In the sequel we will denote $ \mathcal{D}^{(0,1)} $ by $ \mathcal{D}$. Let $\Xi$ stand for the symbol playing in the regularity structure the role of the noise $\zeta$. The regularity structure associated with the generalized (KPZ) equation \eqref{EqGKPZ} is generated by the following rule $\scrR$ that describes the local structure of its decorated trees. Write $[\mcI]_\ell$ for an $\ell$-tuple of operators $\mcI$. One has
$$
\scrR(\Xi) = \{\emptyset\},\quad \scrR(\mcI)=\scrR(\mcI_{(0,1)}) = \Big\{ [\mcI]_\ell, ([\mcI]_\ell,\mcI_{(0,1)}), \big([\mcI]_\ell, \mcI_{(0,1)}, \mcI_{(0,1)}\big) , ([\mcI]_\ell,\Xi)\Big\}_{\ell\in\bbN}
$$
Polynomial decorations $k = (k_0,k_1)\in\mathbb{N}^2$ are added on the nodes of these trees: $ k_0 $ is for time monomial and $ k_1 $ is for spatial monomial. We denote by $\mcB$ the family of all polynomially decorated trees satisfying the rule $\scrR$. By satisfying the rule, we mean that $\tau$ has all its nodes with their incoming edges described by the right hand side of $\scrR(\mcI)$ The empty tree is added to $\mcB$ and denoted by $\bf 1$. To define the degree $\vert\tau\vert$ of an element $\tau\in\mcB$ we set 
$$
\vert{\bf 1}\vert \defeq 0, \quad \vert\Xi\vert \defeq -3/2-\kappa,
$$
and define inductively for $k\in\bbN^2$ and $\sigma, \tau\in\mcB$
$$
\vert X^k\tau\vert \defeq \vert\tau\vert + | k |_{\frak{s}}, \quad \vert\mcI(\tau)\vert \defeq \vert\tau\vert + 2, \quad \vert\mcI_{(0,1)}(\tau)\vert \defeq \vert\tau\vert + 1, \quad \vert\sigma\tau\vert \defeq \vert\sigma\vert + \vert\tau\vert.
$$
One has $ X^k = X_0^{k_0} X_1^{k_1} $ where the $ X_i $ are identified with the canonical basis of $ \mathbb{N}^2 $. We call `{\it noises}' the elements of $\mcB$ of negative degree and denote their collection by $\mcB^-$. Recall that one can build inductively the set of noises using the operators 
\begin{equation} \label{EqNoiseStructure} \begin{split}
	 \mu & \mapsto \mcI_{(0,1)}(\mu),   \\
(\mu, \nu) &\mapsto \mu \mcI(\nu),   \\
(\mu, \nu) &\mapsto \mcI_{(0,1)}(\mu) \, \mcI_{(0,1)}(\nu),
\end{split} \end{equation}
starting from the initial collection of noises $\{\Xi, X_1\Xi\}$. Write $\circ$ for $\Xi$ and use the blue symbol $\begin{tikzpicture}[scale=1,baseline=-0.1cm] \node at (0,0)  [noiseblue] (1) {};\end{tikzpicture}$ for the product symbol $X_1\circ$. With a plain line for the operator $\mcI$ and a dotted line for the operator $\mcI_{(0,1)}$ we can list the elements of $\mcB^-$ according to their degree.
\begin{itemize}
	\item[--] Elements of $\mcB^-$ of degree $-1-2\kappa$:
	$$
	\begin{tikzpicture}[scale=0.3,baseline=0cm]
	\node at (0,0)  [noise] (1) {};
	\node at (0,1.1)  [noise] (2) {};	
	\draw[K] (1) to (2);
	\end{tikzpicture}
	\hspace{0.3cm}, \hspace{0.3cm}
	\begin{tikzpicture}[scale=0.3,baseline=0cm]	
	\node at (0,0)  [dot] (0) {};	
	\node at (0.8,0.8)  [noise] (noise1) {};
	\node at (-0.8,0.8)  [noise] (noise2) {};
	\draw[DK] (0) to (noise1);
	\draw[DK] (0) to (noise2);
	\end{tikzpicture}.
	$$
	
	\item[--] Elements of $\mcB^-$ of degree $-1/2-3\kappa$:
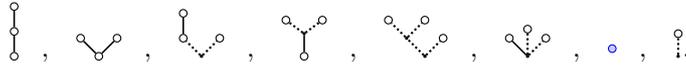
	$$
	\begin{tikzpicture}[scale=0.3,baseline=0cm]
	\node at (0,0)  [noise] (1) {};
	\node at (0,1.1)  [noise] (2) {};
	\node at (0,2.2)  [noise] (3) {};	
	\draw[K] (1) to (2);
	\draw[K] (2) to (3);	
	\end{tikzpicture}
	\hspace{0.3cm}, \hspace{0.3cm}
	\begin{tikzpicture}[scale=0.3,baseline=0cm]
	\node at (0,0)  [noise] (1) {};
	\node at (-0.8,0.8)  [noise] (2) {};
	\node at (0.8,0.8)  [noise] (3) {};
	\draw[K] (1) to (2);
	\draw[K] (1) to (3);	
	\end{tikzpicture}
	\hspace{0.3cm}, \hspace{0.3cm}
	\begin{tikzpicture}[scale=0.3,baseline=0cm]
	\node at (0,0)  [dot] (dot) {};
	\node at (-0.8,0.8)  [noise] (noise1) {};
	\node at (0.8,0.8)  [noise] (noise2) {};
	\node at (-0.8,1.9)  [noise] (noise3) {};	
	\draw[DK] (dot) to (noise1);
	\draw[DK] (dot) to (noise2);
	\draw[K] (noise1) to (noise3);
	\end{tikzpicture}
	\hspace{0.3cm}, \hspace{0.3cm}
	\begin{tikzpicture}[scale=0.3,baseline=0cm]
	\node at (0,0)  [noise] (noise1) {};
	\node at (0,1)  [dot] (dot) {};
	\node at (-0.8,1.6)  [noise] (noise2) {};
	\node at (0.8,1.6)  [noise] (noise3) {};	
	\draw[K] (noise1) to (dot);
	\draw[DK] (dot) to (noise2);
	\draw[DK] (dot) to (noise3);
	\end{tikzpicture}
	\hspace{0.3cm}, \hspace{0.3cm}
	\begin{tikzpicture}[scale=0.3,baseline=0cm]
	\node at (0,0)  [dot] (dot1) {};
	\node at (-0.8,0.8)  [dot] (dot2) {};
	\node at (-1.6,1.6)  [noise] (noise1) {};
	\node at (0,1.6)  [noise] (noise2) {};
	\node at (0.8,0.8)  [noise] (noise3) {};	
	\draw[DK] (dot1) to (dot2);
	\draw[DK] (dot2) to (noise1);
	\draw[DK] (dot2) to (noise2);
	\draw[DK] (dot1) to (noise3);
	\end{tikzpicture}
	\hspace{0.3cm}, \hspace{0.3cm}
	\begin{tikzpicture}[scale=0.3,baseline=0cm]
	\node at (0,0)  [dot] (dot) {};
	\node at (-0.8,0.8)  [noise] (noise1) {};
	\node at (0,1.2)  [noise] (noise2) {};
	\node at (0.8,0.8)  [noise] (noise3) {};	
	\draw[K] (dot) to (noise1);
	\draw[DK] (dot) to (noise2);
	\draw[DK] (dot) to (noise3);
	\end{tikzpicture}
	\hspace{0.3cm}, \hspace{0.3cm}
	\begin{tikzpicture}[scale=1,baseline=-0.1cm] \node at (0,0)  [noiseblue] (1) {};\end{tikzpicture}
	\hspace{0.3cm}, \hspace{0.3cm}
	\begin{tikzpicture}[scale=0.3,baseline=0cm]
	\node at (0,0)  [dot] (dot) {};
	\node at (0,1)  [noise] (noise) {};
	\draw[DK] (dot) to (noise);
	\end{tikzpicture}.
	$$

	\item[--] Elements of $\mcB^-$ of degree $-4\kappa$:
	\begin{equation*} \begin{split}
	&\begin{tikzpicture}[scale=0.3,baseline=0cm]	
	\node at (0,0) [noise] (dot) {};		
	\node at (-0.8,0.8) [noise] (noise1) {};
	\node at (0,1.6) [noise] (noise2) {};
	\node at (-0.8,2.4) [noise] (noise3) {};
	\draw[K] (dot) to (noise1);
	\draw[K] (noise1) to (noise2);
	\draw[K] (noise2) to (noise3);
	\end{tikzpicture}
	\hspace{0.2cm}, \hspace{0.2cm}
	\begin{tikzpicture}[scale=0.3,baseline=0cm]	
	\node at (0,0)  [dot] (0) {};	
	\node at (0.8,0.8)  [noise] (noise1) {};
	\node at (-0.8,0.8)  [noise] (noise2) {};
	\node at (0,1.6)  [noise] (noise3) {};
	\node at (-0.8,2.4)  [noise] (noise4) {};		
	\draw[DK] (0) to (noise1);
	\draw[DK] (0) to (noise2);
	\draw[K] (noise2) to (noise3);
	\draw[K] (noise3) to (noise4);
	\end{tikzpicture}
	\hspace{0.2cm}, \hspace{0.2cm}
	\begin{tikzpicture}[scale=0.3,baseline=0cm]	
	\node at (0,0)  [noise] (noise1) {};	
	\node at (0.8,0.8)  [noise] (noise2) {};
	\node at (0,1.6)  [dot] (dot) {};
	\node at (-0.8,2.4)  [noise] (noise3) {};
	\node at (0.8,2.4)  [noise] (noise4) {};
	\draw[K] (noise1) to (noise2);
	\draw[K] (noise2) to (dot);
	\draw[DK] (dot) to (noise3);
	\draw[DK] (dot) to (noise4);
	\end{tikzpicture}
	\hspace{0.2cm}, \hspace{0.2cm}
	\begin{tikzpicture}[scale=0.3,baseline=0cm]	
	\node at (0,0)  [noise] (noise1) {};	
	\node at (0,1)  [dot] (dot) {};
	\node at (-0.8,1.6)  [noise] (noise2) {};
	\node at (0.8,1.6)  [noise] (noise3) {};
	\node at (0,2.4)  [noise] (noise4) {};
	\draw[K] (noise1) to (dot);
	\draw[DK] (dot) to (noise2);
	\draw[DK] (dot) to (noise3);
	\draw[K] (noise3) to (noise4);
	\end{tikzpicture}
	\hspace{0.2cm}, \hspace{0.2cm}
	\begin{tikzpicture}[scale=0.3,baseline=0cm]	
	\node at (0,0)  [dot] (dot1) {};
	\node at (-0.8,0.8)  [noise] (noise1) {};	
	\node at (0.8,0.8)  [noise] (noise2) {};
	\node at (0,1.6)  [dot] (dot2) {};
	\node at (-0.8,2.4)  [noise] (noise3) {};
	\node at (0.8,2.4)  [noise] (noise4) {};
	\draw[DK] (dot1) to (noise1);
	\draw[DK] (dot1) to (noise2);
	\draw[K] (noise2) to (dot2);
	\draw[DK] (dot2) to (noise3);
	\draw[DK] (dot2) to (noise4);
	\end{tikzpicture}
	\hspace{0.2cm}, \hspace{0.2cm}
	\begin{tikzpicture}[scale=0.3,baseline=0cm]	
	\node at (0,0)  [noise] (noise1) {};
	\node at (0,1)  [dot] (dot1) {};	
	\node at (-0.8,1.6)  [noise] (noise2) {};
	\node at (0.8,1.6)  [dot] (dot2) {};		
	\node at (0,2.4)  [noise] (noise3) {};	
	\node at (1.6,2.4)  [noise] (noise4) {};
	\draw[K] (noise1) to (dot1);
	\draw[DK] (dot1) to (noise2);
	\draw[DK] (dot1) to (dot2);
	\draw[DK] (dot2) to (noise3);
	\draw[DK] (dot2) to (noise4);
	\end{tikzpicture}
	\hspace{0.2cm}, \hspace{0.2cm}
	\begin{tikzpicture}[scale=0.3,baseline=0cm]	
	\node at (0,0)  [dot] (dot1) {};	
	\node at (-0.8,0.8)  [dot] (dot2) {};		
	\node at (-1.6,1.6)  [noise] (noise1) {};
	\node at (0,1.6)  [noise] (noise2) {};		
	\node at (0,2.7)  [noise] (noise3) {};	
	\node at (0.8,0.8)  [noise] (noise4) {};
	\draw[DK] (dot1) to (dot2);
	\draw[DK] (dot1) to (noise4);
	\draw[DK] (dot2) to (noise1);
	\draw[DK] (dot2) to (noise2);
	\draw[K] (noise2) to (noise3);
	\end{tikzpicture}
	\hspace{0.2cm}, \hspace{0.2cm}
	\begin{tikzpicture}[scale=0.3,baseline=0cm]	
	\node at (0,0)  [dot] (dot1) {};	
	\node at (-0.8,0.8)  [dot] (dot2) {};		
	\node at (-1.6,1.6)  [dot] (dot3) {};	
	\node at (0.8,0.8)  [noise] (noise1) {};
	\node at (0,1.6)  [noise] (noise2) {};		
	\node at (-0.8,2.4)  [noise] (noise3) {};	
	\node at (-2.4,2.4)  [noise] (noise4) {};
	\draw[DK] (dot1) to (dot2);
	\draw[DK] (dot1) to (noise1);
	\draw[DK] (dot2) to (noise2);
	\draw[DK] (dot2) to (dot3);
	\draw[DK] (dot3) to (noise3);
	\draw[DK] (dot3) to (noise4);	
	\end{tikzpicture}
	\hspace{0.2cm}, \hspace{0.2cm}
	\begin{tikzpicture}[scale=0.3,baseline=0cm]	
	\node at (0,0)  [dot] (dot1) {};	
	\node at (-0.8,0.8)  [dot] (dot2) {};		
	\node at (0.8,0.8)  [dot] (dot3) {};	
	\node at (1.3,1.4)  [noise] (noise1) {};
	\node at (0.3,1.4)  [noise] (noise2) {};		
	\node at (-0.3,1.4)  [noise] (noise3) {};	
	\node at (-1.3,1.4)  [noise] (noise4) {};
	\draw[DK] (dot1) to (dot2);
	\draw[DK] (dot1) to (dot3);
	\draw[DK] (dot3) to (noise1);
	\draw[DK] (dot3) to (noise2);
	\draw[DK] (dot2) to (noise3);
	\draw[DK] (dot2) to (noise4);
	\end{tikzpicture}
	\hspace{0.2cm}, \hspace{0.2cm}
	\begin{tikzpicture}[scale=0.3,baseline=0cm]	
	\node at (0,0)  [dot] (dot) {};	
	\node at (-0.8,0.8)  [noise] (noise1) {};
	\node at (-0.8,1.8)  [noise] (noise2) {};		
	\node at (0.8,0.8)  [noise] (noise3) {};	
	\node at (0.8,1.8)  [noise] (noise4) {};
	\draw[DK] (dot) to (noise1);
	\draw[DK] (dot) to (noise3);
	\draw[K] (noise1) to (noise2);
	\draw[K] (noise3) to (noise4);
	\end{tikzpicture}
	\hspace{0.2cm}, \hspace{0.2cm}
	\begin{tikzpicture}[scale=0.3,baseline=0cm]	
	\node at (0,0)  [dot] (dot1) {};	
	\node at (0.8,0.8)  [dot] (dot2) {};		
	\node at (-0.8,0.8)  [noise] (noise1) {};
	\node at (-0.8,1.8)  [noise] (noise2) {};		
	\node at (0.3,1.8)  [noise] (noise3) {};	
	\node at (1.3,1.8)  [noise] (noise4) {};
	\draw[DK] (dot1) to (noise1);
	\draw[DK] (dot1) to (dot2);
	\draw[K] (noise1) to (noise2);
	\draw[DK] (dot2) to (noise3);
	\draw[DK] (dot2) to (noise4);
	\end{tikzpicture}
	\hspace{0.2cm}, \hspace{0.2cm}
	\\
	&\begin{tikzpicture}[scale=0.3,baseline=0cm]			
	\node at (0,0)  [noise] (noise1) {};
	\node at (-0.8,0.8)  [noise] (noise2) {};		
	\node at (0,1.1)  [noise] (noise3) {};	
	\node at (0.8,0.8)  [noise] (noise4) {};
	\draw[K] (noise1) to (noise2);
	\draw[K] (noise1) to (noise3);
	\draw[K] (noise1) to (noise4);	
	\end{tikzpicture}
	\hspace{0.2cm}, \hspace{0.2cm}
	\begin{tikzpicture}[scale=0.3,baseline=0cm]			
	\node at (0,0)  [dot] (dot) {};
	\node at (-0.9,0.5)  [noise] (noise1) {};
	\node at (-0.4,0.9)  [noise] (noise2) {};	
	\node at (0.4,0.9)  [noise] (noise3) {};	
	\node at (0.9,0.5)  [noise] (noise4) {};	
	\draw[K] (dot) to (noise1);
	\draw[K] (dot) to (noise2);
	\draw[DK] (dot) to (noise3);
	\draw[DK] (dot) to (noise4);	
	\end{tikzpicture}
	\hspace{0.2cm}, \hspace{0.2cm}
	\begin{tikzpicture}[scale=0.3,baseline=0cm]			
	\node at (0,0)  [noise] (noise1) {};
	\node at (0,1)  [noise] (noise2) {};	
	\node at (-0.8,1.6)  [noise] (noise3) {};	
	\node at (0.8,1.6)  [noise] (noise4) {};	
	\draw[K] (noise1) to (noise2);
	\draw[K] (noise2) to (noise3);
	\draw[K] (noise2) to (noise4);
	\end{tikzpicture}
	\hspace{0.2cm}, \hspace{0.2cm}
	\begin{tikzpicture}[scale=0.3,baseline=0cm]			
	\node at (0,0)  [dot] (dot) {};	
	\node at (0.8,0.8)  [noise] (noise1) {};
	\node at (-0.8,0.8)  [noise] (noise2) {};	
	\node at (-1.6,1.6)  [noise] (noise3) {};	
	\node at (0,1.6)  [noise] (noise4) {};		
	\draw[DK] (dot) to (noise1);
	\draw[DK] (dot) to (noise2);
	\draw[K] (noise2) to (noise3);
	\draw[K] (noise2) to (noise4);	
	\end{tikzpicture}
	\hspace{0.2cm}, \hspace{0.2cm}
	\begin{tikzpicture}[scale=0.3,baseline=0cm]			
	\node at (0,0)  [noise] (noise1) {};
	\node at (0,1)  [dot] (dot) {};	
	\node at (-0.8,1.6)  [noise] (noise2) {};	
	\node at (0,2)  [noise] (noise3) {};	
	\node at (0.8,1.6)  [noise] (noise4) {};		
	\draw[K] (noise1) to (dot);
	\draw[K] (dot) to (noise2);
	\draw[DK] (dot) to (noise3);
	\draw[DK] (dot) to (noise4);	
	\end{tikzpicture}
	\hspace{0.2cm}, \hspace{0.2cm}
	\begin{tikzpicture}[scale=0.3,baseline=0cm]			
	\node at (0,0)  [dot] (dot1) {};	
	\node at (0.8,0.8)  [noise] (noise1) {};
	\node at (-0.8,0.8)  [dot] (dot2) {};
	\node at (-1.6,1.6)  [noise] (noise2) {};
	\node at (-0.8,2)  [noise] (noise3) {};
	\node at (0,1.6)  [noise] (noise4) {};	
	\draw[DK] (dot1) to (noise1);
	\draw[DK] (dot1) to (dot2);
	\draw[K] (dot2) to (noise2);
	\draw[DK] (dot2) to (noise3);	
	\draw[DK] (dot2) to (noise4);	
	\end{tikzpicture}
	\hspace{0.2cm}, \hspace{0.2cm}
	\begin{tikzpicture}[scale=0.3,baseline=0cm]				
	\node at (0,0)  [noise] (noise1) {};
	\node at (-0.8,0.8)  [noise] (noise2) {};
	\node at (0.8,0.8)  [noise] (noise3) {};
	\node at (0,1.6)  [noise] (noise4) {};	
	\draw[K] (noise1) to (noise2);
	\draw[K] (noise1) to (noise3);	
	\draw[K] (noise3) to (noise4);	
	\end{tikzpicture}
	\hspace{0.2cm}, \hspace{0.2cm}
	\begin{tikzpicture}[scale=0.3,baseline=0cm]				
	\node at (0,0)  [dot] (dot) {};
	\node at (-0.8,0.8)  [noise] (noise1) {};
	\node at (-0.8,1.8)  [noise] (noise2) {};
	\node at (0,1)  [noise] (noise3) {};
	\node at (0.8,0.8)  [noise] (noise4) {};	
	\draw[K] (dot) to (noise1);
	\draw[K] (noise1) to (noise2);	
	\draw[DK] (dot) to (noise3);
	\draw[DK] (dot) to (noise4);	
	\end{tikzpicture}
	\hspace{0.2cm}, \hspace{0.2cm}
	\begin{tikzpicture}[scale=0.3,baseline=0cm]				
	\node at (0,0)  [dot] (dot) {};
	\node at (-0.8,0.8)  [noise] (noise1) {};
	\node at (-0.8,1.8)  [noise] (noise2) {};
	\node at (0,1)  [noise] (noise3) {};
	\node at (0.8,0.8)  [noise] (noise4) {};	
	\draw[DK] (dot) to (noise1);
	\draw[K] (noise1) to (noise2);	
	\draw[DK] (dot) to (noise3);
	\draw[K] (dot) to (noise4);	
	\end{tikzpicture}
	\hspace{0.2cm}, \hspace{0.2cm}
	\begin{tikzpicture}[scale=0.3,baseline=0cm]				
	\node at (0,0)  [noise] (noise1) {};
	\node at (0.8,0.8)  [noise] (noise2) {};
	\node at (-0.8,0.8) [dot] (dot) {};
	\node at (0,1.6)  [noise] (noise3) {};
	\node at (-1.6,1.6)  [noise] (noise4) {};	
	\draw[K] (noise1) to (noise2);
	\draw[K] (noise1) to (dot);	
	\draw[DK] (dot) to (noise3);
	\draw[DK] (dot) to (noise4);	
	\end{tikzpicture}
	\hspace{0.2cm}, \hspace{0.2cm}
	\begin{tikzpicture}[scale=0.3,baseline=0cm]				
	\node at (0,0) [dot] (dot1) {};
	\node at (-0.8,0.8) [dot] (dot2) {};	
	\node at (0.8,0.8)  [noise] (noise1) {};
	\node at (0,1)  [noise] (noise2) {};
	\node at (-1.6,1.6)  [noise] (noise3) {};
	\node at (-0.8,1.8)  [noise] (noise4) {};	
	\draw[DK] (dot1) to (noise1);
	\draw[DK] (dot1) to (noise2);
	\draw[K] (dot1) to (dot2);		
	\draw[DK] (dot2) to (noise3);
	\draw[DK] (dot2) to (noise4);	
	\end{tikzpicture}
	\hspace{0.2cm}, \hspace{0.2cm}
	\begin{tikzpicture}[scale=0.3,baseline=0cm]				
	\node at (0,0) [dot] (dot1) {};
	\node at (-0.8,0.8) [dot] (dot2) {};	
	\node at (0.8,0.8)  [noise] (noise1) {};
	\node at (0,1)  [noise] (noise2) {};
	\node at (-1.6,1.6)  [noise] (noise3) {};
	\node at (-0.8,1.8)  [noise] (noise4) {};	
	\draw[K] (dot1) to (noise1);
	\draw[DK] (dot1) to (noise2);
	\draw[DK] (dot1) to (dot2);		
	\draw[DK] (dot2) to (noise3);
	\draw[DK] (dot2) to (noise4);	
	\end{tikzpicture}.
	\end{split} \end{equation*}

	\item[--] Elements of $\mcB^-$ of degree $-2\kappa$:
	\begin{equation*} \begin{split}
	\begin{tikzpicture}[scale=0.3,baseline=0cm]				
	\node at (0,0) [noise] (noise1) {};
	\node at (0,1.2)  [noiseblue] (noise2) {};	
	\draw[K] (noise1) to (noise2);		
	\end{tikzpicture}
	\hspace{0.2cm}, \hspace{0.2cm}
	\begin{tikzpicture}[scale=0.3,baseline=0cm]				
	\node at (0,0) [noiseblue] (noise1) {};
	\node at (0,1.2)  [noise] (noise2) {};	
	\draw[K] (noise1) to (noise2);		
	\end{tikzpicture}
	\hspace{0.2cm}, \hspace{0.2cm}
	\begin{tikzpicture}[scale=0.3,baseline=0cm]				
	\node at (0,0) [dot] (dot) {};
	\node at (-0.8,0.8)  [noiseblue] (noise1) {};
	\node at (0.8,0.8)  [noise] (noise2) {};
	\draw[DK] (dot) to (noise1);
	\draw[DK] (dot) to (noise2);
	\end{tikzpicture}
	\hspace{0.2cm}, \hspace{0.2cm}
	X\,\begin{tikzpicture}[scale=0.3,baseline=0cm]				
	\node at (0,0) [dot] (dot) {};
	\node at (-0.8,0.8)  [noise] (noise1) {};
	\node at (0.8,0.8)  [noise] (noise2) {};
	\draw[DK] (dot) to (noise1);
	\draw[DK] (dot) to (noise2);
	\end{tikzpicture}
	\hspace{0.2cm}, \hspace{0.2cm}
	\begin{tikzpicture}[scale=0.3,baseline=0cm]				
	\node at (0,0) [dot] (dot) {};
	\node at (0,1)  [noise] (noise1) {};
	\node at (0,2)  [noise] (noise2) {};
	\draw[DK] (dot) to (noise1);
	\draw[K] (noise1) to (noise2);
	\end{tikzpicture}
	\hspace{0.2cm}, \hspace{0.2cm}
	\begin{tikzpicture}[scale=0.3,baseline=0cm]				
	\node at (0,0) [dot] (dot1) {};
	\node at (0,0.8) [dot] (dot2) {};
	\node at (-0.8,1.6)  [noise] (noise1) {};
	\node at (0.8,1.8)  [noise] (noise2) {};
	\draw[DK] (dot1) to (dot2);
	\draw[DK] (dot2) to (noise1);
	\draw[DK] (dot2) to (noise2);
	\end{tikzpicture}
	\hspace{0.2cm}, \hspace{0.2cm}
	\begin{tikzpicture}[scale=0.3,baseline=0cm]				
	\node at (0,0) [dot] (dot) {};
	\node at (-0.8,0.8)  [noise] (noise1) {};
	\node at (0.8,0.8)  [noise] (noise2) {};
	\draw[K] (dot) to (noise1);
	\draw[DK] (dot) to (noise2);
	\end{tikzpicture}
	\hspace{0.2cm}, \hspace{0.2cm}
	\begin{tikzpicture}[scale=0.3,baseline=0cm]				
	\node at (0,0) [noise] (noise1) {};
	\node at (0,1)  [dot] (dot) {};
	\node at (0,2)  [noise] (noise2) {};
	\draw[K] (noise1) to (dot);
	\draw[DK] (dot) to (noise2);	
	\end{tikzpicture}
	\hspace{0.2cm}, \hspace{0.2cm}
	\begin{tikzpicture}[scale=0.3,baseline=0cm]				
	\node at (0,0)  [dot] (dot1) {};
	\node at (-0.8,0.8)  [dot] (dot2) {};	
	\node at (0.8,0.8) [noise] (noise1) {};
	\node at (0,1.6) [noise] (noise2) {};
	\draw[DK] (dot1) to (dot2);
	\draw[DK] (dot1) to (noise1);
	\draw[DK] (dot2) to (noise2);
	\end{tikzpicture}.	
	\end{split} \end{equation*}
\end{itemize}
For $\tau\in T$ denote by $\vert\tau\vert_\zeta$ the number of noises in $\tau$. We have $\vert\tau\vert_\zeta\leq 4$ for the trees of $\mcB$ with $\vert\tau\vert<0$. We thus see from Stroock's formula \eqref{EqStroock} that for trees with four noises only the terms $k=0$ and $k=2$ will give a non-null contribution in the sum over $k\leq \vert\tau\vert_\zeta-1$. For trees with three nodes only the term $k=1$ will give a non-null contribution in the sum. This small number of noises in trees of negative degree is the reason why we can develop a direct approach to the convergence of the BHZ renormalized model that does not require the full strength of the BPHZ renormalisation.

We recall the definition of the space $ T^+ $ which is the linear span of 
\begin{equation*}
\bigg\lbrace X^k \prod_{i=1}^n \mathcal{I}^+_{a_i}(\tau_i) \,;\, \tau_i \in \mathcal{B}, \, |\mathcal{I}_{a_i}(\tau_i)| > 0 \bigg\rbrace.
\end{equation*}
We use a different symbol $ \mathcal{I}^+_a $ instead of $ \mathcal{I}_a $ to stress that $ T^+ $ is not a subspace of $ T $. Recall that the co-action $\Delta : T \rightarrow T \otimes T^+$ satisfies the induction relations 

\begin{equation*} \begin{split} 
\Delta(\bullet) &\defeq \bullet\otimes \textsf{\textbf{1}}, \quad \textrm{ for }\bullet\in\big\{\textsf{\textbf{1}},  X_i, \Xi \big\},   \\
\end{split}\end{equation*}
\begin{equation} \label{EqDefnDelta} \begin{split} 
\Delta(\mcI_a\tau) &\defeq (\mcI_a\otimes\textrm{Id})\Delta \tau + \sum_{\vert k +m\vert < | \mcI_a\tau |} \frac{X^k}{k!}\otimes \frac{X^m}{m!}\mcI^+_{a+k+m}(\tau),
\end{split}\end{equation}
Denote by $\big((T,\Delta),(T^+,\Delta^+)\big)$ the BHZ regularity structure associated to that equation. (We do not need here the details of the definition of the splitting map $\Delta^+$.) We denote by $T^-$ the algebra generated by the trees of $T$ of non-positive degree.   \vspace{0.3cm}

\noindent \textit{{\S2. The BHZ model.}} Let $\zeta^\epsilon$ stand for a regularization of $\zeta$ by convolution with a smooth compactly supported function $\rho^\epsilon$ of the form 
$$
\rho^\epsilon(z) = \epsilon^{-3}\rho(\epsilon^{-2}t,\epsilon^{-1}x),
$$ 
with $\rho\geq 0$ of unit integral. Denote by  ${\sf M}^\epsilon = \big({\sf \Pi}^\epsilon, {\sf g}^\epsilon\big)$ the naive admissible model on the preceding regularity structure associated with the smooth noise $\zeta^\epsilon$. We introduce the renormalized naive model $\overline{\sf M}^\varepsilon$ following Bruned's recursive point of view \cite{BrunedRecursive}. We define a  multiplicative map $\delta_r' : T\rightarrow T$ by requiring that
\begin{equation*} \label{coaction} \begin{split} 
\delta_r'(\bullet) &\defeq \bullet\otimes \textsf{\textbf{1}} + \textsf{\textbf{1}} \otimes \bullet, \quad \textrm{ for }\bullet\in\big\{\textsf{\textbf{1}},  X_i, \Xi \big\},   \\
\delta_r'(\mcI_a\tau) &\defeq (\mcI_a\otimes\textrm{Id})\delta_r' \tau + \sum_{k \in \mathbb{N}^2} \frac{X^k}{k!} \otimes \mcI_{a+k}(\tau).
\end{split}\end{equation*}
The infinite sum over $\bbN$ above makes perfect sense using the bigraduation introduced in Section 2.3 of \cite{BHZ}. Let $\pi_-$ be the projection map that keeps only decorated trees with negative degree. For a tree $\tau\in T$ set
\begin{equation*}
\delta_r(\tau) \defeq (\pi_-\otimes\textrm{Id})\delta_r'(\tau).
\end{equation*} 
As all the element of $\mcB^-$ have degree bigger than $-2$ one has $\delta_r\mcI(\tau) = \mcI(\tau)$ for all $\tau\in\mcB$. Using Sweedler's notations for any tree $\tau$ of  non-positive degree one writes
$$
\delta_r(\tau) = \tau \otimes \mathbf{1} + \sum \tau^{(1)} \otimes \tau^{(2)};
$$
the $ \tau^{(1)} $ have strictly less edges than $\tau$ due to the extraction of a subtree of negative degree. Let $\ell$ be a character of the algebra $T^-$ that is {\it null on planted trees and trees of the form} $X^k\sigma$. Set
\begin{equation} \label{def_R}
R_\ell \defeq (\ell\otimes \textrm{Id})\delta_r;
\end{equation}
this map satisfies $R_\ell\mcI_a = \mcI_a$. Define
$$
 {\sf \Pi}^{\!{ R_\ell  }}(\Xi) =   {\sf {\Pi}}^{\!{ R_\ell  }\times}(\Xi) = \zeta^{\eps},
$$
and the maps ${\sf \Pi}^{\!{ R_\ell  }}$ and ${\sf \Pi}^{\!{ R_\ell  }\times}$ defined inductively by the relations
$$
{\sf \Pi}^{\!{ R_\ell  }} = {\sf {\Pi}}^{\!{ R_\ell  }\times} { R_\ell  }, \quad {\sf {\Pi}}^{\!{ R_\ell  }\times}(\tau \sigma) = \big({\sf {\Pi}}^{\!{ R_\ell  }\times}(\tau)\big) \big({\sf {\Pi}}^{\!{ R_\ell  }\times}(\sigma)\big), \quad  {\sf {\Pi}}^{\!{ R_\ell  }\times}(\mathcal{I}_a\tau) = D^{a} K * \big({\sf \Pi}^{\!{ R_\ell  }}(\tau)\big)
$$
where $*$ stands for the spacetime convolution operator. It follows from this definition and the identity $R_\ell\mcI_a=\mcI_a$ that the map ${\sf \Pi}^{\!{ R_\ell  }}$ satisfies the admissibility condition 
$$
 {\sf \Pi}^{\!{ R_\ell  }}(\mcI_a\tau) = D^{a} K * ({\sf \Pi}^{\!{ R_\ell  }} \tau),
$$
The BHZ associated with $\ell$ model is then defined as follows. Set
\begin{equation*} \begin{aligned}
{\sf \Pi}_z^{\! R_\ell } \tau \defeq \left(  {\sf \Pi}^{\! R_\ell } \otimes  ({\sf g}_z^{\! R_\ell })^{-1} \right) \Delta \tau
\end{aligned} \end{equation*}
with
$$
({\sf g}_z^{\! R_\ell })^{-1}\big(\mathcal{I}^+_a\tau\big) \defeq - \big(D^{a} K * {\sf \Pi}^{\! R_\ell }_z\tau\big)(z).
$$
It was proved in Section 3.2 of \cite{BrunedRecursive} that the recentered model ${\sf \Pi}_z^{\! R_\ell }$ satisfies the following identity
\begin{equation*} \begin{aligned}
\big({\sf \Pi}_z^{\! R_\ell } \tau\big)(y) &= \Big({{\sf \Pi}}^{\! R_\ell \times}_z \big( R_\ell  \tau\big)\Big)(y),   \\
{{\sf \Pi}}^{\! R_\ell \times}_z (\tau \sigma) &= ({{\sf \Pi}}^{\! R_\ell \times}_z \tau) \, ({{\sf \Pi}}^{\! R_\ell \times}_z\sigma),  \\ 
\big( {\sf \Pi}_z^{\! R_\ell \times}(\mathcal{I}_a\tau)\big)(y) &= \big(D^{a} K * {\sf \Pi}_z^{\! R_\ell } \tau\big)(y) - \sum_{|k| \leq |\mathcal{I}_a\tau |}\frac{(y-z)^k}{k!}  \big(D^{a +k}  K * {\sf \Pi}_z^{\! R_\ell } \tau\big)(z).
\end{aligned} \end{equation*}
We see from the property of the $\tau^{(1)}$ in Sweedler's decomposition of $\delta_r(\tau)$ that there is a unique character $\ell$ on $T^-$ that is null on planted trees and trees of the form $X^k\sigma$ such that 
\begin{equation*} 
\mathbb{E}\big[ ({\sf \Pi}^{\!{ R_\ell  }} \tau)(0) \big] = 0
\end{equation*}
for all $\tau\in T^-$. This character is called the {\sl BHZ character}; it satisfies the identity
\begin{equation} \label{BPHZ_character}
\mathbb{E} \big[ ({\sf \Pi}^{\!{ R_\ell  }} \tau)(y) \big] = 0
\end{equation}
for all $\tau\in\mcB^-$ and $y$, from translation invariance of the law of $\zeta$ and of the operator $\partial_t-\partial_x^2$.

\subsection{Tree-like graphs and their mirror graphs}
\label{SubsectionMirror}

The functions ${\sf \Pi}^\epsilon\tau$ or $\overline{\sf \Pi}^\epsilon\tau$ will eventually be integrated against some smooth test functions $\varphi$. We introduce a graphical notation for that operation

$$
\langle f , \varphi \rangle \defeq \int f(y) \varphi(y)dy \eqdef \begin{tikzpicture}[scale=0.5,baseline=0cm]
	\node at (0,0)  [root] (root) {};
	\node at (0,1)  [noise] (1) {};
	
	\draw[testfcn] (root) to node[labl, pos=1.24] {\tiny {\color{black}$f$}} (1);
\end{tikzpicture}.
$$
We use a green line to denote the integration operation and a green bullet {\color{testcolor} $\bullet$} to represent the test function. The function $\varphi$ will not vary in our computations so it is harmless that our graphical notation does not refer to $\varphi$ explicitly. It turns out to be useful for the computations to be done below to represent the function ${\sf \Pi}^\epsilon\tau$ with the same graphical notation as $\tau$ itself. The $\epsilon$-dependence of the analytic object  disappears in this notation but the context makes it clear. So we will write $\begin{tikzpicture}[scale=0.4,baseline=0.1cm]
	\node at (0,0)  [root] (root) {};
	\node at (0,1)  [noise] (1) {};
	
	\draw[testfcn] (root) to node[labl, pos=1.24] {\small {\color{black}$\tau$}} (1);
\end{tikzpicture}$ for $\int ({\sf \Pi}^\epsilon\tau)(y) \, \varphi(y)dy$.

The definition of a {\sl mirror graph} is best illustrated by a self-explaining example rather than by a formal definition. Here is an example of a tree representation of a test operation together with its mirror graph.
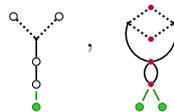
\begin{equation*}
\begin{tikzpicture}[scale=0.3,baseline=.8cm]
	\node at (0,0)  [root] (root) {};
	\node at (0,1)  [noise] (1) {};
	\node at (0,2)  [noise] (2) {};
	\node at (0,3)  [dot] (3) {};	
	\node at (-1,4)  [noise] (left) {};
	\node at (1,4)  [noise] (right) {};	
	
	\draw[testfcn] (1) to (root);    
	\draw[K] (1) to (2);
	\draw[K] (2) to (3);
	\draw[DK] (right) to (3);    
	\draw[DK] (left) to (3);   
\end{tikzpicture}
\hspace{0.3cm}, \hspace{0.3cm}
\begin{tikzpicture}[scale=0.3,baseline=.8cm]
	\node at (-0.5,0)  [root] (root1) {};
	\node at (0.5,0)  [root] (root2) {};	
	\node at (0,1)  [var] (1) {};
	\node at (0,2)  [var] (2) {};
	\node at (0,3)  [var] (right) {};	
	\node at (-1,3.7)  [dot] (3left) {};	
	\node at (1,3.7)  [dot] (3right) {};	
	\node at (0,4.4)  [var] (left) {};
	
	\draw[testfcn] (1) to (root1);    
	\draw[testfcn] (1) to (root2);
	\draw[K] (1) to[bend left=+50] (2);
	\draw[K] (1) to[bend left=-50] (2);
	\draw[K] (2) to[bend left=+50] (3left);
	\draw[K] (2) to[bend left=-50] (3right);
	\draw[DK] (3left) to (right);
	\draw[DK] (3right) to (right);
	\draw[DK] (3left) to (left);
	\draw[DK] (3right) to (left);	
\end{tikzpicture}
\end{equation*}
We merge here the pairs of noises that correspond via the mirror symmetry. In the mirror graph a {\color{purple} purple $\bullet$} represents the `merging' of two regularized noises, in our case the spacetime convolution operator with the kernel $\rho^\epsilon*\rho^\epsilon$. This kernel converges to $\delta_0(z'-z)$ as the regularization parameter $\epsilon>0$ goes to $0$. A plain edge in a mirror graph represents a spacetime convolution with the function $K$ and a dotted line a spacetime convolution with the kernel $\partial_x K$. In the above example the diagram on the right hand side represents a double integral of the form 
\begin{equation} \label{EqDoubleIntegral}
\iint \varphi(y) \varphi(y')\,\mcM^\epsilon(\tau)(y,y')\,dydy',
\end{equation}
with

\begin{equation*} \begin{split}
\mcM^\epsilon(\tau)(y,y') = \int (&\rho^\epsilon*\rho^\epsilon)(y-y')K(y-z_1)K(y'-z'_1)(\rho^\epsilon*\rho^\epsilon)(z_1-z'_1)K(z_1-z_2)K(z'_1-z'_2)   \\
&\times \partial^{(0,1)} K(z_2-z_{31}) \partial^{(0,1)} K(z_2-z_{32}) \partial^{(0,1)} K(z'_2-z'_{31}) \partial^{(0,1)} K(z'_2-z'_{32})   \\
&\times (\rho^\epsilon*\rho^\epsilon)(z_{31}-z'_{31})(\rho^\epsilon*\rho^\epsilon)(z_{32}-z'_{32}) \, dz_1 dz_2 d z_{31} dz_{32} dz'_1  dz'_2 dz'_{31} dz'_{32}.
\end{split} \end{equation*}
Its formal limit as $\epsilon$ goes to $0$ would be a simple integral of the form 
$$
\int \varphi(y)^2\,\mcM^0(\tau)(y,y)\,dy,
$$
with $\mcM^0$ defined as $\mcM^\varepsilon$ with $\delta_0$ in place of $\rho^\epsilon*\rho^\epsilon$. We use here the generic notation $\mcM$ for `{\it mirror}'. Here are two other examples of tree representations of test operations together with their mirror graphs.

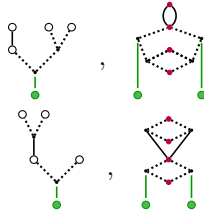
\begin{equation*}
\begin{tikzpicture}[scale=0.3,baseline=.4cm]
	\node at (0,0)  [root] (root) {};
	\node at (0,1)  [dot] (1) {};
	\node at (-1,2)  [noise] (2) {};
	\node at (1,2)  [dot] (3) {};
	\node at (-1,3)  [noise] (4) {};
	\node at (0.6,3)  [noise] (5) {};
	\node at (1.6,3)  [noise] (6) {};	
		
	\draw[testfcn] (1) to (root);
	\draw[DK] (1) to (2);
	\draw[DK] (1) to (3);
	\draw[K] (2) to (4); 
	\draw[DK] (3) to (5);
	\draw[DK] (3) to (6);
\end{tikzpicture}
\hspace{0.3cm}, \hspace{0.3cm}
\begin{tikzpicture}[scale=0.3,baseline=.4cm]
	\node at (-1.4,0)  [root] (root1) {};
	\node at (1.4,0)  [root] (root2) {};	
	\node at (0,1)  [var] (var1) {};
	\node at (0,2)  [var] (var2) {};
	\node at (0,3)  [var] (var3) {};
	\node at (0,4)  [var] (var4) {};		
	\node at (-1,1.5)  [dot] (dotl1) {};
	\node at (-1.4,2.5)  [dot] (dotl2) {};	
	\node at (1,1.5)  [dot] (dotr1) {};
	\node at (1.4,2.5)  [dot] (dotr2) {};
	
	\draw[testfcn] (dotl2) to (root1);
	\draw[testfcn] (dotr2) to (root2);
	\draw[DK] (var1) to (dotl1);
	\draw[DK] (var1) to (dotr1);
	\draw[DK] (var2) to (dotl1);
	\draw[DK] (var2) to (dotr1);
	\draw[DK] (dotl1) to (dotl2);
	\draw[DK] (dotl1) to (var2);	
	\draw[DK] (dotr1) to (dotr2);
	\draw[DK] (dotr1) to (var2);	
	\draw[DK] (dotl2) to (var3);
	\draw[DK] (dotr2) to (var3);
	\draw[K] (var3) to[bend left=-50] (var4);
	\draw[K] (var3) to[bend left=+50] (var4);
\end{tikzpicture}
\end{equation*}

\begin{equation*}
\begin{tikzpicture}[scale=0.3,baseline=.4cm]
	\node at (0,0)  [root] (root) {};
	\node at (0,1)  [dot] (dot1) {};
	\node at (-1,2)  [noise] (n1) {};
	\node at (1,2)  [noise] (n2) {};
	\node at (-1,3)  [dot] (dot2) {};
	\node at (-1.5,4)  [noise] (n3) {};
	\node at (-0.5,4)  [noise] (n4) {};
		
	\draw[testfcn] (dot1) to (root);
	\draw[DK] (dot1) to (n1);
	\draw[DK] (dot1) to (n2);
	\draw[K] (n1) to (dot2); 
	\draw[DK] (dot2) to (n3);
	\draw[DK] (dot2) to (n4);
\end{tikzpicture}
\hspace{0.3cm}, \hspace{0.3cm}
\begin{tikzpicture}[scale=0.3,baseline=.4cm]
	\node at (-1,0)  [root] (root1) {};
	\node at (1,0)  [root] (root2) {};	
	\node at (0,1)  [var] (var1) {};
	\node at (0,2)  [var] (var2) {};
	\node at (0,2.8)  [var] (var3) {};
	\node at (0,3.8)  [var] (var4) {};		
	\node at (-1,1.5)  [dot] (dotl1) {};
	\node at (-1,3.3)  [dot] (dotl2) {};	
	\node at (1,1.5)  [dot] (dotr1) {};
	\node at (1,3.3)  [dot] (dotr2) {};
	
	\draw[testfcn] (dotl1) to (root1);
	\draw[testfcn] (dotr1) to (root2);
	\draw[DK] (var1) to (dotl1);
	\draw[DK] (var2) to (dotl1);
	\draw[DK] (var1) to (dotr1);
	\draw[DK] (var2) to (dotr1);
	\draw[K] (dotl2) to (var2);
	\draw[K] (dotr2) to (var2);
	\draw[DK] (var3) to (dotl2);
	\draw[DK] (var4) to (dotl2);	%
	\draw[DK] (var3) to (dotr2);
	\draw[DK] (var4) to (dotr2);
\end{tikzpicture}
\end{equation*}
Here again the two diagrams on the right hand side represent a double integral of the form \eqref{EqDoubleIntegral}. More generally mirror diagrams are defined in the obvious way from any connected finite graphs built with the bricks $\circ, \begin{tikzpicture}[scale=1,baseline=-0.1cm] \node at (0,0)  [noiseblue] (1) {};\end{tikzpicture}, {\color{testcolor} \bullet}$ and the solid and dotted edges. We talk of a {\sl tree-like graph} if it corresponds to the quantity obtained by testing a polynomial functional of the regularized noise against a test function $\varphi$. We call it `tree-like' rather than `tree' as its associated analytic object does not necessarily have a tree structure. As an example, the following quantities are tree-like graphs
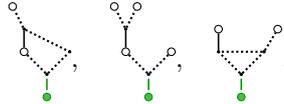
\begin{equation*}
\begin{tikzpicture}[scale=0.3,baseline=.4cm]
	\node at (0,0)  [root] (root) {};
	\node at (0,1)  [dot] (dot1) {};
	\node at (-1,2)  [noise] (n1) {};
	\node at (1,2)  [dot] (dot2) {};
	\node at (-1,3)  [dot] (dot3) {};
	\node at (-1.5,4)  [noise] (n3) {};
		
	\draw[testfcn] (dot1) to (root);
	\draw[DK] (dot1) to (n1);
	\draw[DK] (dot1) to (dot2);
	\draw[K] (n1) to (dot3); 
	\draw[DK] (dot3) to (n3);
	\draw[DK] (dot3) to (dot2);
\end{tikzpicture},\quad 
\begin{tikzpicture}[scale=0.3,baseline=.4cm]
	\node at (0,0)  [root] (root) {};
	\node at (0,1)  [dot] (dot1) {};
	\node at (-1,2)  [noise] (n1) {};
	\node at (1,2)  [noise] (n2) {};
	\node at (-1,3)  [dot] (dot2) {};
	\node at (-1.5,4)  [noise] (n3) {};
	\node at (-0.5,4)  [noise] (n4) {};
		
	\draw[testfcn] (dot1) to (root);
	\draw[DK] (dot1) to (n1);
	\draw[DK] (dot1) to (n2);
	\draw[K] (n1) to (dot2); 
	\draw[DK] (dot2) to (n3);
	\draw[DK] (dot2) to (n4);
\end{tikzpicture},
\quad
\begin{tikzpicture}[scale=0.3,baseline=.4cm]
	\node at (0,0)  [root] (root) {};
	\node at (0,1)  [dot] (1) {};
	\node at (-1,2)  [dot] (2) {};
	\node at (1,2)  [dot] (3) {};
	\node at (-1,3)  [noise] (4) {};
	\node at (1.6,3)  [noise] (6) {};	
		
	\draw[testfcn] (1) to (root);
	\draw[DK] (1) to (2);
	\draw[DK] (1) to (3);
	\draw[K] (2) to (4); 
	\draw[DK] (3) to (2);
	\draw[DK] (3) to (6);
\end{tikzpicture}.
\end{equation*}
On can define more generally tree-like graphs where each edge represents a kernel that may depend not only on the two integration variables associated with the vertices of the edge but also on the other integration variables. This is what happens when we replace a kernel by a Taylor remainder of that kernel based at some point different from the two vertices variables corresponding to the edge associated with the kernel. This more general setting will be studied in Section \ref{SubsectionWorkhorse}.

\ssk

 Mirror graphs will be involved below when estimating the norm of the expectation of the $n$th order Malliavin derivative of a polynomial of the noise for $n\neq 0$.


\subsection{Malliavin derivatives of trees and their expectation} 
\label{SubsectionMalliavin}

For $0\leq n\leq \vert\tau\vert_\zeta-1$ denote by $d^n\overline{\sf \Pi}^\epsilon_x\tau$ the $n$th order Malliavin derivative of $\overline{\sf \Pi}^\epsilon_x\tau$; this is an element of $L\big(H^{\otimes n},L^2(\Omega,\bbP)\big)$. On the example of the tree 
$$
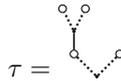
\tau=\begin{tikzpicture}[scale=0.3,baseline=0cm]
	\node at (0,0)  [dot] (dot1) {};
	\node at (-1,1)  [noise] (n1) {};
	\node at (1,1)  [noise] (n2) {};
	\node at (-1,2)  [dot] (dot2) {};
	\node at (-1.5,3)  [noise] (n3) {};
	\node at (-0.5,3)  [noise] (n4) {};
		
	\draw[DK] (dot1) to (n1);
	\draw[DK] (dot1) to (n2);
	\draw[K] (n1) to (dot2); 
	\draw[DK] (dot2) to (n3);
	\draw[DK] (dot2) to (n4);
\end{tikzpicture}
$$ 
here is  a piece of $\big\langle\big(d{\sf \Pi}^\epsilon\tau\big)(h) , \varphi\big\rangle$

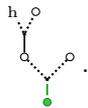
\begin{equation*}
\begin{tikzpicture}[scale=0.3,baseline=.4cm]
	\node at (0,0)  [root] (root) {};
	\node at (0,1)  [dot] (dot1) {};
	\node at (-1,2)  [noise] (n1) {};
	\node at (1,2)  [noise] (n2) {};
	\node at (-1,3)  [dot] (dot2) {};
	\node at (-1.5,4)  [noise] (n3) {};
	\node at (-0.5,4)  [noise] (n4) {};
		
	\draw[testfcn] (dot1) to (root);
	\draw[DK] (dot1) to (n1);
	\draw[DK] (dot1) to (n2);
	\draw[K] (n1) to (dot2); 
	\draw[DK] (dot2) to node[labl, pos=1.24] {\tiny h} (n3);
	\draw[DK] (dot2) to (n4);
\end{tikzpicture}\;.
\end{equation*}
As $\zeta$ has the same law as $-\zeta$ only quantities with an even number of noises have possibly non-null expectation. The first Malliavin derivative of a $4$-linear function of the noise being $3$-linear it has null expectation. Here are some pieces of $\big\langle\big(d^2{\sf \Pi}^\epsilon\tau\big)(h_1,h_2) , \varphi\big\rangle$ for the preceding tree

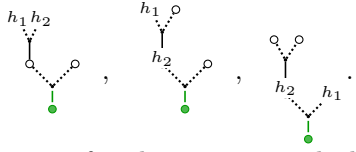
\begin{equation*}
\begin{tikzpicture}[scale=0.3,baseline=.4cm]
	\node at (0,0)  [root] (root) {};
	\node at (0,1)  [dot] (dot1) {};
	\node at (-1,2)  [noise] (n1) {};
	\node at (1,2)  [noise] (n2) {};
	\node at (-1,3)  [dot] (dot2) {};
	\node at (-1.5,4)  [noise] (n3) {};
	\node at (-0.5,4)  [noise] (n4) {};
		
	\draw[testfcn] (dot1) to (root);
	\draw[DK] (dot1) to (n1);
	\draw[DK] (dot1) to (n2);
	\draw[K] (n1) to (dot2); 
	\draw[DK] (dot2) to node[labl, pos=1.24] {\tiny $h_1$} (n3);
	\draw[DK] (dot2) to node[labl, pos=1.24] {\tiny $h_2$} (n4); 
\end{tikzpicture} 
\hspace{0.3cm}, \hspace{0.3cm}
\begin{tikzpicture}[scale=0.3,baseline=.4cm]
	\node at (0,0)  [root] (root) {};
	\node at (0,1)  [dot] (dot1) {};
	\node at (-1,2)  [noise] (n1) {};
	\node at (1,2)  [noise] (n2) {};
	\node at (-1,3.4)  [dot] (dot2) {};
	\node at (-1.5,4.4)  [noise] (n3) {};
	\node at (-0.4,4.4)  [noise] (n4) {};
		
	\draw[testfcn] (dot1) to (root);
	\draw[DK] (dot1) to (n1);
	\draw[DK] (dot1) to (n2);
	\draw[K] (n1) to node[labl, pos=0.05] {\tiny $h_2$} (dot2); 
	\draw[DK] (dot2) to node[labl, pos=1.24] {\tiny $h_1$} (n3);
	\draw[DK] (dot2) to (n4);
\end{tikzpicture}
\hspace{0.3cm}, \hspace{0.3cm}
\begin{tikzpicture}[scale=0.3,baseline=.8cm]
	\node at (0,0)  [root] (root) {};
	\node at (0,1)  [dot] (dot1) {};
	\node at (-1,2)  [noise] (n1) {};
	\node at (1,2)  [noise] (n2) {};
	\node at (-1,3.5)  [dot] (dot2) {};
	\node at (-1.5,4.4)  [noise] (n3) {};
	\node at (-0.4,4.4)  [noise] (n4) {};
		
	\draw[testfcn] (dot1) to (root);
	\draw[DK] (dot1) to (n1);
	\draw[DK] (dot1) to node[labl, pos=1.24] {\tiny $h_1$}(n2);
	\draw[K] (n1) to node[labl, pos=0.05] {\tiny $h_2$} (dot2); 
	\draw[DK] (dot2) to (n3);
	\draw[DK] (dot2) to (n4);  
\end{tikzpicture}.
\end{equation*}
We are interested in the expectation of such quantities, which produce elements of the space of Hilbert-Schmidt operators on $H^{\otimes n}$ that are tree-like graph maps where the former noise arguments $\circ, \begin{tikzpicture}[scale=1,baseline=-0.1cm] \node at (0,0)  [noiseblue] (1) {};\end{tikzpicture}$ from Section \ref{SubsectionMirror} are now some elements of $H$. We will use a common notation $\Vert\cdot\Vert$ for the norms on these spaces. The square norm of a tree-like graph map is bounded above by its corresponding mirror diagram, where the pairing of an element of $H$ with its mirror element produces a vertex in the mirror diagram. As an example here are some pieces of $\bbE\big[ \big\langle\big(d^2{\sf \Pi}^\epsilon\tau\big)(h_1,h_2) , \varphi \big\rangle\big]$, for the same tree $\tau$ as above
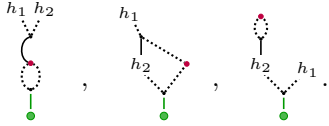
\begin{equation} \label{graph_examples}
\begin{tikzpicture}[scale=0.35,baseline=.4cm]
	\node at (0,0)  [root] (root) {};
	\node at (0,1)  [dot] (dot1) {};
	\node at (0,2)  [var] (var) {};
	\node at (0,3)  [dot] (dot2) {};
	\node at (-0.5,4)  [noise] (n3) {};
	\node at (0.5,4)  [noise] (n4) {};
		
	\draw[testfcn] (dot1) to  (root);
	\draw[DK] (dot1) to [bend left=+70] (var);
	\draw[DK] (dot1) to [bend left=-70] (var);
	\draw[K] (var) to [bend left=0] (dot2); 
	\draw[DK] (dot2) to node[labl, pos=1.24] {\tiny $h_1$} (n3);
	\draw[DK] (dot2) to node[labl, pos=1.24] {\tiny $h_2$} (n4);   
\end{tikzpicture} 
\hspace{0.3cm}, \hspace{0.3cm}
\begin{tikzpicture}[scale=0.3,baseline=.4cm]
	\node at (0,0)  [root] (root) {};
	\node at (0,1)  [dot] (dot1) {};
	\node at (-1,2)  [noise] (n1) {};
	\node at (1,2.3)  [var] (var) {};
	\node at (-1,3.5)  [dot] (dot2) {};
	\node at (-1.5,4.5)  [noise] (n3) {};
		
	\draw[testfcn] (dot1) to (root);
	\draw[DK] (dot1) to (n1);
	\draw[DK] (dot1) to (var);
	\draw[K] (n1) to node[labl, pos=0.05] {\tiny $h_2$} (dot2); 
	\draw[DK] (dot2) to node[labl, pos=1.24] {\tiny $h_1$} (n3);
	\draw[DK] (dot2) to (var);   
\end{tikzpicture}
\hspace{0.3cm}, \hspace{0.3cm}
\begin{tikzpicture}[scale=0.3,baseline=.4cm]
	\node at (0,0)  [root] (root) {};
	\node at (0,1)  [dot] (dot1) {};
	\node at (-1,2)  [noise] (n1) {};
	\node at (1,2)  [noise] (n2) {};
	\node at (-1,3.4)  [dot] (dot2) {};
	\node at (-1,4.4)  [var] (var) {};
		
	\draw[testfcn] (dot1) to (root);
	\draw[DK] (dot1) to (n1);
	\draw[DK] (dot1) to node[labl, pos=1.24] {\tiny $h_1$}(n2);
	\draw[K] (n1) to node[labl, pos=0.05] {\tiny $h_2$} (dot2); 
	\draw[DK] (dot2) to [bend left=+50] (var);
	\draw[DK] (dot2) to [bend left=-50] (var);   
\end{tikzpicture}.
\end{equation}
The loop in the third diagram is a constant annihilated by the kernel $K$, so the diagram represents the null function. We see on that example that looking at $\bbE\big[ \big\langle\big(d^2\overline{\sf \Pi}^\epsilon\tau\big)(h_1,h_2) , \varphi\big\rangle\big]$ for the BHZ renormalized model $\overline{\sf \Pi}^\epsilon$ removes exactly the diverging contribution of the first graph in the sum giving $\bbE\big[ \big\langle\big(d^2{\sf \Pi}^\epsilon\tau\big)(h_1,h_2) , \varphi\big\rangle\big]$;  one can see from Section \ref{SectionConvergence} that the other terms are converging as the regularization parameter $\epsilon$ is sent to $0$. 

\ssk

For a given tree $\tau$ the $L^2(\Omega)$ random variable $\big\langle \overline{\sf \Pi}_z^\epsilon\tau , \varphi^\lambda_z\big\rangle$ is Malliavin smooth and its $k$th order Malliavin derivative coincides with its Fr\'echet $k$th order derivative. Its chaos decomposition is given by Stroock's formula
\begin{equation} \label{EqStroockFormula}
\big\langle \overline{\sf \Pi}_z^\epsilon\tau , \varphi^\lambda_z\big\rangle = \big\langle \bbE\big[\hspace{0.02cm}\overline{\sf \Pi}_z^\epsilon\tau\big] , \varphi^\lambda_z\big\rangle + \sum_{n=1}^{\vert\tau\vert_\zeta} \frac{1}{n!}\, I_n\Big(\big\langle \bbE\big[d^n\overline{\sf \Pi}_z^\epsilon\tau \big] , \varphi_z^\lambda\big\rangle\Big)
\end{equation}
where $I_p$ is the canonical isometry from the symmetric tensor space $\{L^2([0,T]\times \bbT)\}^{\odot n}$ onto the $n$th Wiener chaos. See for instance Section 2.7.2 in Nourdin \& Peccati's book \cite{NourdinPeccati}. If $\Phi^{\epsilon,\lambda}_{\tau,\varphi}\in L^2\big(([0,T]\times\bbT)^n\big)$ stands for the kernel of the real-valued map on $\{L^2([0,T]\times\bbT)\}^n$
$$
\big(h_1,\dots,h_n\big) \mapsto \bbE\big[\big\langle \big( d^n \overline{\sf \Pi}_z^\epsilon\tau\big)(h_1,\dots,h_n) , \varphi_z^\lambda\big\rangle\big]
$$
it is a consequence of Jensen's inequality that 
$$
\bigg\Vert I_n\bigg(\Big\langle \bbE\big[d^n\overline{\sf \Pi}_z^\epsilon\tau \big] , \varphi_z^\lambda\Big\rangle\bigg)\bigg\Vert_{L^2(\Omega)}^2 \leq \int_{([0,T]\times\bbT)^n} \Phi^{\epsilon,\lambda}_{\tau,\varphi}(z_1,\dots,z_n)^2\,dz_1\dots dz_n.
$$
See e.g. Section 3 of Mourrat, Weber \& Xu's work \cite{MourratWeberXu} for this nice remark that was already used in Hairer's seminal work \cite{Hai14}. Mirror graphs are graphical notations for the above right hand sides. As particular cases, for a tree-like graph $\tau$ with two noises one has
$$
\textsc{Var}\big(\begin{tikzpicture}[scale=0.4,baseline=0.1cm]
	\node at (0,0)  [root] (root) {};
	\node at (0,1)  [noise] (1) {};	
	\draw[testfcn] (root) to node[labl, pos=1.24] {\small {\color{black}$\tau$}} (1); 
\end{tikzpicture}\big)^2 \leq \iint \varphi(y) \varphi(y')\,\mcM^\epsilon(\tau)(y,y')\,dydy',
$$
and for an arbitrary tree $\tau$ one has
$$
\bigg\Vert I_{\vert\tau\vert_\zeta}\bigg(\Big\langle d^{\vert\tau\vert_\zeta}\overline{\sf \Pi}_z^\epsilon\tau , \varphi_z^\lambda\Big\rangle\bigg)\bigg\Vert \leq \iint \varphi_z^\lambda(y) \varphi_z^\lambda(y')\,\mcM^\epsilon(\tau)(y,y')\,dydy'.
$$
The term $I_{\vert\tau\vert_\zeta}(\cdot)$ is actually deterministic and the notation $\Vert\cdot\Vert$ stands here for the $L^2$ norm on $(L^2([0,T]\times\bbT))^{\vert\tau\vert_\zeta}$.

\ssk

\noindent \textbf{{Noise derivative operators --}} We formalize the noise differentiation at the level of the regularity structure by introducing some derivative operators on trees that replace a noise with new noise symbols $\Xi_j$ with degree equal to 
$$ 
| \Xi_j | = -\kappa.
$$ 
In the context of the generalized (KPZ) equation we are interested in some stochastic iterated integrals having at most four noises, so the index $ j $ will belongs to the finite set $ \lbrace 1,2,3,4 \rbrace $. We define the derivative operators $ D_{\Xi_j} $ by

\begin{equation*} \begin{split}  
D_{\Xi_j}\Xi \defeq& \; \Xi_j,   \\
\end{split} \end{equation*}
\begin{equation*} \begin{split}  
D_{\Xi_j} \bullet  \defeq& \; 0, \quad \textrm{ for }\bullet\in\Big\{\textsf{\textbf{1}} ; X_i ; \Xi_{\ell}, 1\leq\ell\leq 4 \Big\},   \\
D_{\Xi_j} \left( \tau \tau' \right) =& \; \left( D_{\Xi_j} \tau \right) \tau' + \tau \left( D_{\Xi_j} \tau' \right),    \\ 
D_{\Xi_j} \mathcal{I}_a(\tau) = & \; \mathcal{I}_a( D_{\Xi_j}\tau).
\end{split} \end{equation*}
We define $\widehat{\mcB}^-$ as the decorated trees obtained by applying the iterated derivative operator $ D_{\Xi_{j_1}} \cdots D_{\Xi_{j_k}} $ where $ j_1,..., j_k $ are distinct elements of $ \lbrace 1,2,3,4 \rbrace $, to an element of $ \mcB^- $. We assume that $ \mcB^- \subset \widehat{\mcB}^- $. We can see the elements of $\mcB^- $ as some elements of $ \widehat{\mcB}^- $ where one is able to replace $k$ noises $\Xi$ by $ \Xi_{j_1},...,\Xi_{j_k} $. Each $ \Xi_{j_k} $ appears at most one time in an element of $\widehat{\mcB}^-$.  Denote by $\widehat{\mcB}^-_j$ the set of elements of $\widehat{\mcB}^-$ with a $\Xi_j$ symbol. No element of $ \widehat{\mcB}^-_{j} $ having at least one noise of type $  \Xi_j$ has a subtree of negative degree containing this special noise. (Indeed, this decorated trees is coming from a decorated tree belonging to $ \mcB^- $ where $\Xi_j$ is replaced $ \Xi $. Then a negative subtree different from a single noise that contains $ \Xi_j $ will be a negative subtree that contain $ \Xi $. The difference in degree between the two subtrees is  $ |\Xi| - |\Xi_j| =  - \frac{3}{2}$. This means that one has a subtree different from $ \Xi $ which is of degree below $ - \frac{3}{2} $. This is absurd.)
 Then one has
\begin{equation*}
\delta_r D_{\Xi_j} = \big( \textrm{Id} \otimes D_{\Xi_j} \big)\delta_r,
\end{equation*}
which implies that
\begin{equation*}
R_{\ell} D_{\Xi_j} = D_{\Xi_j} R_{\ell}.
\end{equation*}
We extend the definition of the pre-model and the model to $ \widehat{\mcB}^- $ by associating to each noise $ \Xi_j $ an element $ h_j$ of the Hilbert space $H$. For the pre-model we set
$$
\big({\sf \Pi}^{\!R_{\ell}}\Xi_j\big)(h_1,\dots, h_4) =   \big({\sf {\Pi}}^{\!R_{\ell}\times} \Xi_j\big)(h_1,\dots, h_4) = h_j,
$$
For the model we use a different degree where
$$
|\Xi_j|_+ \defeq -\frac{3}{2} - \kappa,
$$ 
and $\vert\cdot\vert_+$ is defined as $\vert\cdot\vert$ on the rest of the decorated trees. Given that $R_\ell\mcI_a=\mcI_a$, one has for any $\tau\in\widehat{\mcB}^-$
\begin{equation} \label{recentering_Xi}
\begin{aligned}
\big( {\sf \Pi}_z^{\!R_{\ell}}(\mathcal{I}_a\tau)\big)(y,h_1,\dots,h_4) &= \big( {\sf \Pi}_z^{\!R_{\ell\times}}(\mathcal{I}_a\tau)\big)(y,h_1,\dots,h_4)   \\
&= \big(D^{a} K * {\sf \Pi}_z^{\!R_{\ell}} \tau\big)(y,h_1,\dots,h_4)   \\ 
&\qquad- \sum_{|k| < |\mathcal{I}_a\tau |_+}\frac{(y-z)^k}{k!}  \big(D^{a+k}  K * {\sf \Pi}_z^{\!R_{\ell}} \tau\big)(z,h_1,\dots,h_4).
\end{aligned}
\end{equation}
Proposition 4.1 in \cite{BrunedNadeem} says that $d\overline{\sf \Pi}^\epsilon_z\tau = \overline{\sf \Pi}^\epsilon_z\big( D_{\Xi_1}\tau\big)$. One can proceed by induction on $k$ to prove the following fact.

\medskip

\begin{lem} \label{lemma_commutation}
For $\tau\in\mcB^-, 0\leq k\leq \vert\tau\vert_\zeta$ and all $z$ one has
$$
d^k\overline{\sf \Pi}^\epsilon_z\tau = \overline{\sf \Pi}^\epsilon_z\big(D_{\Xi_k}\dots D_{\Xi_1}\tau\big).
$$
\end{lem}

\medskip

We extend the character $\ell$ to the regularity structure built for the generalized (KPZ) with the five noises $\Xi, \Xi_1, \dots, \Xi_4$ with the same expression on $ \mcB^- $ and being zero on $ \widehat{\mcB}^- \setminus \mcB^- $. Indeed, an element of $ \widehat{\mcB}^- \setminus \mcB^- $ is of positive degree if it is not the single noise $  \Xi_j$. Moreover the subtree $ \Xi_j $ does not correspond to a subdivergence and it is interpreted as $h_j$.

\section{Convergence of the BHZ model}
\label{SectionConvergence}

We prove in this section the convergence of the BHZ model from \S2 of Section \ref{SubsectionRS} after recalling a useful criterion to get scaling bounds on mirror graphs. For a test function $\varphi$ write $\varphi_z$ for $\varphi(\cdot-z)$ and $\varphi_z^\lambda$ for $\lambda^{-3}\varphi\big(\frak{s}(\lambda)(\cdot-z)\big)$, where $\frak{s}(\lambda)(t,x) \defeq (\lambda^{-2}t,\lambda^{-1} x)$.


\subsection{A workhorse}
\label{SubsectionWorkhorse}

We will use as our workhorse a refinement due to Bruned \& Nadeem \cite{BN1} of a result of Hairer \& Quastel \cite{HairerQuastel} giving some scaling bounds on Feynman graphs -- Theorem A.3 of \cite{HairerQuastel}. We need such a refinement as renormalisations of the type \eqref{Hairer_Quastel_reno} (two incoming edges on a subdivergent diagram) are not covered by the convergence theorem of Hairer \& Quastel. We will use this refinement in our setting to deal with the mirror graphs that come from estimating the terms $\Vert\bbE[ \langle d^k\overline{\sf \Pi}^\epsilon_z\tau , \varphi_z^\lambda\rangle]\Vert$ with $1<\vert k\vert_{\frak{s}} < \vert\tau\vert_\zeta$. 

We restrict its statement to our setting, in which our graphs ${\bf G} = ({\bf V},{\bf E})$ always have two green vertices, with each of them attached to the remainder of the graph by one green edge;  the latter represent integration against a scaled test function $\varphi^\lambda$. Except from the green edges, each edge $ e  = (e_+,e_-) $ in our graphs represents a kernel of the form
\begin{equation} \label{Taylor_expansion}
L_e(z_{e_-}, z_{e_+} ;  z_{v_e}) = K_e(z_{e_+} - z_{e_-}) - \sum_{\vert j\vert_{\frak{s}}< r_e} \frac{(z_{e_+} - z_{v_e})}{j!} \partial^jK_e(z_{v_e} - z_{e_-}),
\end{equation}
the Taylor remainder of an edge-dependent function $K_e(\cdot - z_{e_-})$ around an edge-dependent point $z_{v_e}$. Note the strict inequality $<r_e$ in the right hand side of \eqref{Taylor_expansion}, so there is no Taylor expansion at all if $r_e=0$. We will have
$$
r_e\in\{0,1,2\},
$$
so we have at most some first order Taylor remainders. Our edges will be either $K, \partial K:= \partial^{(0,1)} K$ or $\rho^\epsilon*\rho^\epsilon$. We note that we always have $r_e=0$ for edges for which $K_e=\rho^\epsilon*\rho^\epsilon$. We assign the green integration edge a number $a_e=0$ and set 
\begin{align*}
a_e \defeq \left\{
\begin{aligned}  1\; & \textrm{if } K_e=K,    \\   2\; &\textrm{if } K_e=\partial K,   \\  3\; &\textrm{if } K_e = \rho^\epsilon*\rho^\epsilon.
\end{aligned}   
\right.
\end{align*}
This is the scaling degree of divergence of the kernel $K_e$ near $0$. The Taylor remainder \eqref{Taylor_expansion} can be coded via a decoration $\gamma_e=(a_e,r_e,v_e)$ on the associated edge.

\ssk

 Let us explain how one casts iteratively a decorated tree into one of these graphs. Let us consider $ \tau = X^k \Xi_j \prod_{i=1}^n \mathcal{I}_{a_i}(\tau_i) $ with $a_i \in \left\lbrace  (0,0), (0,1) \right\rbrace  $. We suppose that we have constructed the labelled graphs $ G_i $ associated to $ \tau_i $. We denote by $ \rho_i $  the node of $ G_i $ which is connected to the green vertex. Then, we identify all the green vertices coming from the graphs $G_i$ into a single node. We introduce a new node $\rho $ which we connect to each node $ \rho_i $ via an edge labelled by $(b_i, r_i,v_0)$ where $ v_0 $ is the green vertex, $ b_i $  correspond to some $a_e$ depending on the value of $a_i$, $r_e$ is computed as the degree of the graph $ G $. This degree corresponds to the degree of the decorated tree $\tau$ associated to $G$ -- see \cite[Prop. 4.3]{BN1}. In fact, the Taylor remainder given in \eqref{Taylor_expansion} matches exactly the recentering procedure described in \eqref{recentering_Xi}. We also connect the node $\rho$ to a new node via an edge decorated by $ (3 + \kappa,0,v_0) $ corresponding to a noise type edge. Finally, we add an edge between $\rho$ and $v_0$ decorated by $(-|k|,0,v_0)$ corresponding to the monomial $X^k$. In general the graphs associated in this way with some decorated trees are not trees due to the monomial decorations that create some cycles. Below, we illustrate this construction:
\begin{equation*}
	\begin{tikzpicture}[baseline=0.1cm]
		\node at (0,-1)  [root, label=below:\tiny{$ v_0 $}] (root) {};
		\node at (-1,1) [dot, label=left:$\tilde{G}_1$] (a) {};
		\node at (1,1) [dot, label=right:$\tilde{G}_n$] (b) {};
		\node at (-1.5,0) [dot, label=right:$$] (e) {};
		\node at (0,0.9) [label=center:$\cdots$] (d) {};
		\node at (0,0) [dot, label=right:\tiny{$\rho$}] (c) {};
		\draw[testfcn] (c) to (root);
		\draw[semithick, - >] (a) to  (c);
		\draw[semithick, - >] (e) to  (c);
		\draw[semithick, - >] (b) to  (c);
		\node at (0.4,0.5) [ label=right:\tiny{$(b_n,r_n,v_0)$}] () {};
         \node at (-1,-0.3) [ label=center:\tiny{$(3 + \kappa,0,v_0)$}] () {};
          \node at (-0.5,0.5) [ label=left:\tiny{$(b_1,r_1,v_0)$}] () {};
          	\draw[semithick, < - ] (c) to[bend left=60]  (root);
          	 \node at (0,-0.5) [ label=right:\tiny{$(-|k|,0,v_0)$}] () {};
	\end{tikzpicture}
	\end{equation*}
	In the graph above, $ \tilde{G}_i $ is obtained from $G_i$ by removing the green vertex and the edges connected to it. Here, in the graph, we have assumed that the $\tau_i$ have no monomials $X^k$. Otherwise, one has to add extra edges between the $\tilde{G}_i$ and the vertex $v_0$.

\ssk

For a subset $V\subset {\bf V}$ of vertices in our graph we define ${\bf E}_{\textrm{int}}(V)$ as the set of edges $e$ in the graph with $e_+$ and $e_-$ in $V$. We define ${\bf E}^{\uparrow}(V)$, resp. ${\bf E}^{\downarrow}(V)$, as the set of edges of the graph with $e_-\in V$, resp. $e_+\in V$. Last we define ${\bf E}^{\uparrow}_{r >0}(V)$ as the set of edges $e$ of ${\bf E}^{\uparrow}(V)$ such that $r_e > 0$. The following two conditions jointly define \textbf{\textsf{Condition (C)}}.


\begin{itemize}
	\item \textsf{\textbf{Integrability/power counting.}} {\it For every subset $V$ of vertices of our graph with at least two elements we have}
	\begin{equation} \label{condition_integrability}
	\sum_{e\in {\bf E}_{\textrm{int}}(V)} a_e + \sum_{e\in {\bf E}^{\uparrow}_{r > 0}(V)} {\bf 1}_{v_e\in V} (a_e+r_e-1) - \sum_{e\in {\bf E}^{\downarrow}(V)} {\bf 1}_{v_e\in V} r_e < 3(\vert V\vert -1).
	\end{equation}
	
	\item \textsf{\textbf{Recentering.}} {\it For every subset $V$ of vertices of our graph \emph{not containing the green vertices} one has} 
	\begin{equation} \label{condition_recentering}
	\sum_{e \in {\bf E}_{\textrm{int}}(V)} a_e + \sum_{e\in {\bf E}^{\downarrow}(V)} \Big({\bf 1}_{\{v_e\in V\}\cup\{r_e=0\}}(a_e+r_e-1) - (r_e-1)\Big) + \sum_{e\in {\bf E}^{\uparrow}(V)} \big(a_e+ r_e{\bf 1}_{v_e\notin V}\big) > 3\vert V\vert.
	\end{equation}
\end{itemize}

\ssk

The number $3$ in the inequalities \eqref{condition_integrability} and \eqref{condition_recentering} is the scaling dimension of the covariance of the one dimensional spacetime white noise. We denote by $\mcE_\star\subset {\bf E}$ the set of non-green edges of our graph ${\bf G}$ and write 
$$
\mathscr{I}^{\bf G}(\lambda) \defeq \int\prod_{e\in\mcE_\star} L_e(z_{e_-}, z_{e_+} ;  z_{v_e}) \, \varphi^\lambda(z_1)\varphi^\lambda(z_2) \, dz
$$
with $z_1,z_2$ the integration variables associated with the two green vertices and $dz$ a shorthand notation for all the integration variables, with $z_1$ and $z_2$ included. We denote by $V_\star$ the set of non-green vertices of $\bf G$, and set
$$
\alpha \defeq 3\vert V_\star\vert - \sum_{e\in\mcE_\star} a_e.
$$
The following statement is Theorem 3.1 in \cite{BN1}.

\ssk

\begin{thm} \label{ThmWorkhorse}
If our graph $\bf G$ satisfies Condition \textsf{\textbf{(C)}} then there exists a constant $c\in(0,\infty)$ depending only on $\bf G$ such that 
$$
\vert \mathscr{I}^{\bf G}(\lambda)\vert \leq c \lambda^\alpha \prod_{e\in\mcE_\star} \Vert K_e\Vert_{a_e,1}.
$$
\end{thm}

\ssk

Condition \eqref{condition_integrability} checks the integrability of the Feynman diagram. It can be understood as giving the integrability on all subsets of variables by selecting a subset $V\subset{\bf  V}$ of vertices. It is also an equivalent of Weinberg's power counting given by \cite{Wein} on Feynman diagrams. Condition \eqref{condition_recentering} allows to get the correct scaling of $ \lambda $ in the end, as expressed in Theorem \ref{ThmWorkhorse}. This is the reason why the green nodes, which are already scaled in $\lambda$, are excluded from the subsets $V$.

\ssk

For a spacetime point $z$ denote by $\frak{T}_z$ the translation operator 
$$
h\in H \mapsto h(\cdot+z) \in H
$$
on $H$; this is an isometry of $H$. As one has the equality in law
$$
\big\langle d^k\overline{\sf \Pi}^\epsilon_z, \varphi^\lambda_{z}\big\rangle(h_1,\dots,h_k) {\overset{\textrm{law}}{=}} \big\langle d^k\overline{\sf \Pi}^\epsilon_0, \varphi^\lambda_0\big\rangle(\frak{T}_z h_1,\dots,\frak{T}_z h_k)
$$
the mirror graph of $\bbE\big[\big\langle d^k\overline{\sf \Pi}^\epsilon_0, \varphi^\lambda_0\big\rangle\big]$ provides a $z$-independent upper bound for the norm of 
$$
\bbE\big[\big\langle d^k\overline{\sf \Pi}^\epsilon_z , \varphi^\lambda_z\big\rangle\big].
$$
We now check that the mirror graphs associated to the fully differentiated quantity $\big\langle d^{\vert\tau\vert_\zeta} \overline{\sf \Pi}^\epsilon_z\tau , \varphi_z^\lambda\big\rangle$ satisfy condition \textbf{\textsf{(C)}} for every $\tau\in \mcB^-$.  One first notices that one has
$$
d^{\vert\tau\vert_\zeta} \overline{\sf \Pi}^{\epsilon}_z\tau = d^{\vert\tau\vert_\zeta}{\sf \Pi}^\epsilon_z\tau
$$
for every $\tau\in\mcB^-$ as a consequence of the fact that $\overline{\sf \Pi}^{\epsilon}_z\tau - {\sf \Pi}^\epsilon_z\tau$ has a number of noises strictly smaller than $ \vert\tau\vert_\zeta $. 

\ssk

\begin{cor} \label{PropMirrorDiagrams}
One has for all $0<\lambda\leq 1$ and every $\tau\in \mcB^-$ the estimates
$$
\big\Vert \big\langle d^{\vert\tau\vert_\zeta} \overline {\sf \Pi}^\epsilon_z\tau , \varphi_z^\lambda\big\rangle \big\Vert \lesssim \lambda^{\vert\tau\vert},
$$
for an implicit multiplicative constant independent of $\epsilon>0$ and $z$.
\end{cor}

\ssk

\begin{Dem}
The \textbf{\textsf{integrability bound}} \eqref{condition_integrability} is satisfied due to the fact that one does not see any subdivergence in a full mirror graph without contraction. 

For the \textbf{\textsf{recentering bound}} \eqref{condition_recentering} one can proceed by {\it induction} on the construction of the diagram as in Section 4 in \cite{BN1}. 

One first notices that $ \big\Vert \big\langle d^{\vert\tau\vert_\zeta} \overline {\sf \Pi}^\epsilon_z\tau , \varphi_z^\lambda\big\rangle \big\Vert $ is built from two labelled graphs $G_1$ and $G_2$, each associated to $d^{\vert\tau\vert_\zeta} \overline {\sf \Pi}^\epsilon_z\tau$ as above, as in the example
\begin{equation*}
	\begin{tikzpicture}[scale=0.3,baseline=.4cm]
		\node at (-1,0)  [root] (root1) {};
		\node at (1,0)  [root] (root2) {};	
		\node at (0,1)  [var] (var1) {};
		\node at (0,2)  [var] (var2) {};
		\node at (0,2.8)  [var] (var3) {};
		\node at (0,3.8)  [var] (var4) {};		
		\node at (-1,1.5)  [dot] (dotl1) {};
		\node at (-1,3.3)  [dot] (dotl2) {};	
		\node at (1,1.5)  [dot] (dotr1) {};
		\node at (1,3.3)  [dot] (dotr2) {};
		
		\draw[testfcn] (dotl1) to (root1);
		\draw[testfcn] (dotr1) to (root2);
		\draw[DK] (var1) to (dotl1);
		\draw[DK] (var2) to (dotl1);
		\draw[DK] (var1) to (dotr1);
		\draw[DK] (var2) to (dotr1);
		\draw[K] (dotl2) to (var2);
		\draw[K] (dotr2) to (var2);
		\draw[DK] (var3) to (dotl2);
		\draw[DK] (var4) to (dotl2);	%
		\draw[DK] (var3) to (dotr2);
		\draw[DK] (var4) to (dotr2);
	\end{tikzpicture} \rightarrow 	G_1 = G_2 =  \begin{tikzpicture}[scale=0.3,baseline=.4cm]
	\node at (0,0)  [root] (root) {};
	\node at (0,1)  [dot] (dot1) {};
	\node at (-1,2)  [noise] (n1) {};
	\node at (1,2)  [noise] (n2) {};
	\node at (-1,3)  [dot] (dot2) {};
	\node at (-1.5,4)  [noise] (n3) {};
	\node at (-0.5,4)  [noise] (n4) {};
		\draw[K] (n1) to (dot2); 
		\draw[DK] (dot2) to node[labl, pos=1.24] {\tiny $h_4$} (n3);
	\draw[DK] (dot2) to node[labl, pos=1.24] {\tiny $h_3$} (n4); 
	\draw[testfcn] (dot1) to (root);
		\draw[DK] (dot1) to node[labl, pos=1.24] {\tiny $h_1$}(n2);
	\draw[K] (n1) to node[labl, pos=0.05] {\tiny $h_2$} (dot2); 
	\draw[DK] (dot1) to (n1);
	\draw[DK] (dot1) to (n2);
\end{tikzpicture}. \vspace{0.15cm}
\end{equation*} 
Here $ G_1 $ and $G_2$ are some labelled trees where we put the same $ h_i $ to  make clear that these nodes are identified  in $G$.

 For an arbitrary decorated graph $G=(V,E)$ we define a variant of \eqref{condition_recentering}
\begin{equation} \label{new_bound} \begin{split}
	 \sum_{e \in {\bf E}_{\textrm{int}}(V)} a_e &+ \sum_{e\in {\bf E}^{\downarrow}(V)} \big({\bf 1}_{\{v_e\in V\}\cup\{r_e=0\}}(a_e+r_e-1) - (r_e-1)\big) + \sum_{e\in {\bf E}^{\uparrow}(V)} \big(a_e+ r_e{\bf 1}_{v_e\notin V}\big)  \\ 
	 &> 3\vert V_i\vert + \frac{3}{2}\vert V_l\vert
\end{split} \end{equation}
where $V_i \subset V$  is the subsets of the inner nodes of $V$ and $V_l \subset V$  is the subsets of leaves of $V$.  (Unlike \eqref{condition_recentering}, this condition only deals with the full graph $G$ rather than with a number of its subgraphs.)  It turns out that this bound is satisfied  by the graphs $G_i$ associated with some decorated trees. The key point for proving that fact inductively is the choice of the $r_e$ that guarantees that we get the correct bound in the end.  Here, we assume that $\tau = \mathcal{I}_{a_1}(\tau_1)$ and we denote by $G$ (resp. $G_1$) the graph associated to $\tau$ (resp. $\tau_1$) with root $\rho$ (resp. $\rho_1$). If we suppose that $v_0 \notin V$ and that all the nodes of $G$ are contained in $V$ one can rewrite the previous inequality into
\begin{equation*}
		\sum_{e \in {\bf E}(V)} a_e  + r_{e_{\rho}}> 3\vert V_i\vert + \frac{3}{2}\vert V_l\vert
\end{equation*}
where $ e_{\rho} $ is the edge connecting $ \rho $ to $\rho_1$  which gives 
\begin{equation*}
	 r_{e_{\rho}}> 3\vert V_i\vert + \frac{3}{2}\vert V_l\vert - 	\sum_{e \in E(V)} a_e = |\tau |.
	\end{equation*}
The last equality is immediate to check as for example a noise type edge gives a contribution $ a_e = 3 + \kappa $ and the leaf associated to it gives a contribution $ \frac{3}{2} $. One has
\begin{equation*}
	\frac{3}{2} -  (3 + \kappa) = - \frac{3}{2}- \kappa
\end{equation*}
that corresponds to the degree of the noise. 

Now the two graphs $G_1$ and $G_2$ used for constructing a mirror diagram both satisfy the bound \eqref{new_bound}; it follows that the mirror diagram satisfies the recentering condition \eqref{condition_recentering}. 
 Indeed, when summing the two bounds for $G_1$ and $G_2$, the leaves are identified so they add up with the correct factor $\frac{1}{2}$.
\end{Dem}

\medskip

We now recall for the reader's convenience two stability results from \cite{BN1} that we will use -- Proposition 5.4 and 5.7 therein.  Let $ {\bf G} = ({\bf V},{\bf E}) $ be a labelled graph  that satisfies the recentering condition \eqref{condition_recentering} but which has some subgraphs $G = (V,E)$ that  do not satisfy the integrability condition \eqref{condition_integrability}.  For $e\in {\bf E}^{\downarrow}(V)$ with decoration $\gamma_e=(a_e,r_e,v_e)$, we look at the effect on the integrability condition of a change of base point $(v_e\rightarrow v)$ in the corresponding Taylor expansion. Set 
$$
L_e(z_e^+,z_e^- ; z_{v_e}) = K_e(z_e^+-z_e^-) - \sum_{\vert j\vert<r_e} \frac{(z_e^+-z_{v_e})^j}{j!} \, \partial^j K_e(z_{v_e} - z_e^-).
$$
This quantity is equal to 
\begin{equation*} \begin{split}
&K_e(z_e^+ - z_e^-) - \sum_{\vert k\vert<r'_e} \frac{(z_e^+ - z_v)^k}{k!} \, \partial^k K_e(z_v - z_e^-)   \\
&\qquad+ \sum_{\vert k\vert<r'_e} \frac{(z_e^+ - z_v)^k}{k!} \, \bigg(\partial^kK_e(z_v-z_e^-) - \sum_{\vert\ell\vert<r_e-\vert k\vert}  \frac{(z_v-z_{v_e})^\ell}{\ell!} \, \partial^{k+\ell}K_e (z_{v_e}-z_e^-)\bigg).
\end{split} \end{equation*}
We rewrite this change of base point formula in the form of the following diagramatic representation
\begin{equation} \label{eqdecomp}
\begin{tikzpicture}[baseline=0.5mm]
	\node at (0,0) [dot, label=right:$v$] (a) {}; 
	\node at (0,1) [dot, label=right:$v_e^+$] (b) {};
	\node at (0,2) [dot, label=right:$v_e^-$] (c) {};
	\draw[semithick, - >] (c) to node[labl,pos=0.45] {\tiny $(e,\gamma_e)$} (b);
	\draw[snake=zigzag, segment amplitude=0.5pt,segment length = 1mm, line after snake = 1mm, - >] (b) to (a);
\end{tikzpicture} = \begin{tikzpicture}[baseline=0.5mm]
	\node at (0,0) [dot, label=right:$v$] (a) {}; 
	\node at (0,1) [dot, label=right:$v_e^+$] (b) {};
	\node at (0,2) [dot, label=right:$v_e^-$] (c) {};
	\draw[semithick, - >] (c) to node[labl,pos=0.45] {\tiny $(e,\gamma'_e)$} (b);
	\draw[snake=zigzag, segment amplitude=0.5pt,segment length = 1mm, line after snake = 1mm, - >] (b) to (a);
\end{tikzpicture}+ \sum_{|k|_{\frak{s}}<r_e'}
\begin{tikzpicture}[baseline=0.1cm]
	\node at (-1,1) [dot, label=left:$v_e^-$] (a) {};
	\node at (1,1) [dot, label=right:$v_e^+$] (b) {};
	\node at (0,0) [dot, label=right:$v$] (c) {};
	\draw[semithick, - >] (a) to node[labl,pos=0.4] {\tiny $(\overline{e},\overline\gamma_k)$} (c);
	\draw[snake=zigzag, segment amplitude=0.5pt,segment length = 1mm, line after snake = 1mm, - >] (b) to (c);
	\draw[semithick, - >] (c) to[bend left=60] node[labl,pos=0.48] {\tiny $(\underline{e},\underline\gamma_k)$} (b);
\end{tikzpicture}
\end{equation}
where the symbol $\begin{tikzpicture}[baseline=-1mm]
	\node at (0,0) (a) {}; 
	\node at (1,0) (b) {};
	\draw[snake=zigzag, segment amplitude=0.5pt,segment length = 1mm, line after snake = 1mm, - >] (b) to (a);
\end{tikzpicture}$ means that there exists a path between $v_2$ and $v$ {\it within} $G$. We represent here only a portion of the graph $\bf G$. The vertices $v_e^-,v_e^+, v$ may thus also be attached to some other vertices in $\bf G$ which have not been pictured here. These vertices and edges of $\bf G$ are not touched by the re-expansion procedure pictured in \eqref{eqdecomp}. Here we set
\begin{equation*} \begin{split}
\gamma_e &= (a_e,r_e,v_e)   \\
\gamma'_e &= (a_e, r'_e, v)   \\
\overline\gamma_k &= \Big(a_e + |k|_{\frak{s}}, \max\big(r_e-|k|_{\frak{s}},0\big) , v_e\Big)   \\
\underline\gamma_k &= (-|k|_{\frak{s}},0,v_e)
\end{split} \end{equation*}
where 
$$
r_e'=\max\big(\lceil-|G|_{\mathfrak{s}} \rceil, r_e\big) \qquad \textrm{and}\qquad  |G|_{\mathfrak{s}} = - 3(\vert V\vert -1) + \sum_{e\in E} a_e
$$ 
is the degree of divergence of the graph $G$.

\ssk

 The change of base point improves the integrability condition for the first term in the right hand side. We need to control what happened of the recentering condition. The next two statement is about the stability of the recentering condition \eqref{condition_recentering} with respect to the previous re-expansion procedure.

\ssk

\ssk	

\begin{prop} \label{PropPreservedRecenteringBounds}
If the recentering condition \eqref{condition_recentering} is satisfied in \eqref{eqdecomp} for the terms with the edge  $(v_e^-,v_e^+)$ labelled by $\gamma_e$ and $\gamma'_e$ then  each term in the sum over $k$ also satisfies \eqref{condition_recentering}.
\end{prop}

\ssk

\begin{Dem}
The proof depends on the subdivergent subgraph. It has to be performed on three subgraphs listed in \cite[Eq. 5.18]{BN1}: they are the only subdivergences that appear in the study of the (gKPZ) equation. One can look at the proof of \cite[Prop. 5.7]{BN1}  for the treatment of the first two subgraphs and at  \cite[Sec. 5.3]{BN1} for the treatment of the third one.
\end{Dem}

\ssk

 We denote by ${\bf G}(\gamma_e')$ the labelled graph corresponding to the first term in the right hand side of \eqref{eqdecomp}. Each term in the sum over $k$ in \eqref{eqdecomp} corresponds to a graph ${\bf G}'$ different from $\bf G$, with $k$-dependent labels but fixed vertex set $\bf E$. We check which subsets $V'$ of $\bf E$ in the labelled graphs ${\bf G}'_k$ satisfy assumption \eqref{condition_recentering}.

\ssk

\begin{prop} \label{prop:invcond3}
 Assume ${\bf G}$ and ${\bf G}(\gamma_e')$ both satisfy the recentering condition \eqref{condition_recentering} for all $V\subset {\bf E}$. Then every subsets $V'\subset {\bf E}$ of vertices of ${\bf G}'_k$ such that $V' \cap\{v,v_e^-,v_e^+\} \neq \{v\}$ satisfies the condition \eqref{condition_recentering}.
\end{prop}

\ssk

\begin{Dem}
Given that the condition \eqref{condition_recentering} is satisfied on terms with $\gamma_e$ and $\gamma'_e$ it suffices to prove that the contribution of $\overline\gamma_k$ and $\underline\gamma_k$ to the left hand side of \eqref{condition_recentering} is greater than the minimum of the contributions of $\gamma_e$ and $\gamma'_e$. To this end we fix for a given subset $V\subset {\bf V}$ of vertices the set $\widetilde{V} \defeq V\cap\{v,v_e^-,v_e^+\}$ and consider the contribution of each of the four edges for every possible $\widetilde{V}$. The situation is summarized in the following tabular.
	\\
	\vspace{3mm}
	\begin{center}
		\begin{tabular}{|c|c|c|c|c|}
			\hline 
			\vphantom{\huge{A}} $\widetilde{V}$ & $\overline\gamma_k$ & $\underline\gamma_k$ & $\gamma'_e$ & $\gamma_e$ \\
			\hline
			$\substack{{\color{white} 0} \\ {\color{white} 0} \\ \{v\} \\ {\color{white} 0} \\ {\color{white} 0}}$ & $\substack{- {\bf 1}_{\{r_e-|k|_{\mathfrak{s}}>0\}}(r_e-|k|_{\mathfrak{s}}-1) \\ + {\bf 1}_{\{r_e-|k|_{\mathfrak{s}}= 0\}}(a_e+|k|_{\mathfrak{s}})}$ & $-|k|_{\mathfrak{s}} $ & $ 0 $ & $ 0 $\\
			\hline
			$\{v_e^-\} $ & $ a_e+|k|_{\mathfrak{s}}+\text{max}\big(r_e-|k|_{\mathfrak{s}},0\big) $ & $ 0 $ & $ a_e+r'_e $ & $ a_e+r_e $   \\
			\hline
			$\{v_e^+\} $ & $ 0 $ & $ -|k|_{\mathfrak{s}} $ & $ -(r'_e-1) $ & $ -(r_e-1) $   \\
			\hline
			$\{v,v_e^-\}$ & $ a_e+|k|_{\mathfrak{s}} $ & $ -|k|_{\mathfrak{s}} $ & $ a_e $ & $ a_e + r_e $   \\
			\hline
			$\substack{{\color{white} 0} \\ {\color{white} 0} \\ \{v,v_2\} \\ {\color{white} 0} \\ {\color{white} 0}}$ & $\substack{- {\bf 1}_{\{r_e-|k|_{\mathfrak{s}}>0\}}(r_e-|k|_{\mathfrak{s}}-1) \\ + {\bf 1}_{\{r_e-|k|_{\mathfrak{s}}=0\}}(a_e+|k|_{\mathfrak{s}})} $ & $ -|k|_{\mathfrak{s}} $ & $ -(r'_e-1) $ & $ -(r_e-1) $   \\
			\hline
			$\{v_e^-,v_e^+\} $ & $ a_e+|k|_{\mathfrak{s}}+\text{max}\big(r_e-|k|_{\mathfrak{s}},0\big) $ & $ -|k|_{\mathfrak{s}} $ & $ a_e $ & $ a_e $   \\
			\hline
			$\{v,v_e^-,v_e^+\}$ & $ a_e+|k|_{\mathfrak{s}} $ & $ -|k|_{\mathfrak{s}} $ & $ a_e$ & $ a_e $   \\
			\hline
		\end{tabular}
	\end{center}
	\vspace{3mm}
	Hence the sum of the contributions of $\overline\gamma_k$ and $\underline\gamma_k$ is indeed greater than the minimum between $\gamma_e$ and $\gamma_e'$ except for $\widetilde{V}=\{v\}$, when $r_e-|k|_{\mathfrak{s}}>0$. In that case one has $-(r_e - 1) \leq 0$, which in turn implies the required result.
\end{Dem}
	
\ssk

 Proposition \ref{PropPreservedRecenteringBounds} and Proposition \ref{prop:invcond3} are used at the very end of the proof of Proposition \ref{PropBoundDerivatives}. We  also need to look at the following decomposition in the sequel
	\begin{equation}
\label{eq:decomp3}
	\begin{tikzpicture}[baseline=0cm,scale=0.8]
		\node at (0,-1) [dot,label=left:$v_e^+$] (a) {};
		\node at (0,1) [dot,label=left:$v_e^-$] (b) {}; 
		\draw[kepsus] (b) to node[labl,pos=0.45] {\tiny $(e,\gamma_e)$} (a);
	\end{tikzpicture} = 
	\begin{tikzpicture}[baseline=0cm,scale=0.8]
		\node at (0,-1) [dot,label=left:$v_e^+$] (a) {};
		\node at (0,1) [dot,label=left:$v_e^-$] (b) {};
		\draw[kepsus] (b) to node[labl,pos=0.45] {\tiny $(e,\gamma_e')$} (a);
	\end{tikzpicture} + \sum_{\vert k\vert < r_e} \begin{tikzpicture}[baseline=0cm,scale=0.8]
		\node at (1.5,-1) [dot,label=below:$v_e$] (a) {};
		\node at (1.5,1) [dot,label=right:$v_e^-$] (b) {};
		\node at (0,0.5) [dot,label=left:$v_e^+$] (c) {};
		\draw[kepsus] (b) to node[labl,pos=0.4] {\tiny $(a_e + \vert k\vert,0)$} (a);
		\draw[kepsus] (c) to node[labl,pos=0.45] {\tiny $(-\vert k\vert,0)$} (a);
	\end{tikzpicture}
	\end{equation}
	where we have the following labels for $\gamma_e = (a_e,r_e,v_e), \, \gamma_e' = (a_e,0,0)$ and $v_e^+, v_e^-$ belong to a labelled graph $\bf G$ that satisfies the bound \eqref{condition_recentering}  for all $V$ not containing the green vertices.

\ssk	

\begin{prop} \label{prop_c}
 Under this assumption, the $V$-dependent condition \eqref{condition_recentering} is satisfied for all the terms of the previous decomposition for all subsets  $V$ of $\bf E$ that either do not contain $v_e^-$ or contain the pair $\{v_e^-,v_e^+\}$.
\end{prop}

\ssk

\begin{Dem}
Let $\overline{V} \subset V$ where  $V$ is the set of nodes of the graph described in \eqref{eq:decomp3}. Setting $\widetilde V=\overline V\cap\{v_e^+,v_e^-\}$, we have the following cases.
		$$
		\begin{array}{|c|c|c|c|c|}
			\hline
			\widetilde V & \gamma_e & \gamma_e' &  (v_e^+,0) &  (v_e^-,0)    \\
			\hline
			\{v_e^+,v_e^-\} & a_e & a_e & -\vert k\vert & a_e + \vert k\vert   \\
			\hline
			\{v_e^+\} & -(r_e - 1) & a_e & -\vert k\vert & 0   \\
			\hline
			\{v_e^-\} & a_e + r_e & a_e & 0 & a_e + \vert k\vert   \\
			\hline
		\end{array}
		$$
For the first two rows of the previous table, the contribution of $\gamma_e'$  and the sum of  $(0,v_e^+)$ and  $(0,v_e^-)$ are greater than the contribution of $\gamma_e$. It is for $\{v_e^-\}$ and $r_e > 0$ that we need to be more careful and  should notice on a case by cases basis that $a_e$ is sufficient for the required bound instead of $a_e + r_e$.
\end{Dem}
	
\ssk	
	
The previous proposition leaves many cases to check by hand.

\medskip

\subsection{Proof of convergence}
\label{SubsectionProofConvergence}

\begin{prop} \label{PropBoundDerivatives}
There exists a positive finite constant $C$ such that one has the bound
\begin{equation} \label{EqScalingBoundExpectation}
\big\Vert\bbE\big[ \big\langle d^n\overline{\sf \Pi}^\epsilon_z\tau , \varphi_z^\lambda\big\rangle\big]\big\Vert \leq C\lambda^{\vert\tau\vert}
\end{equation}
for all $\tau\in\mcB^-$ and $0\leq n\leq \vert\tau\vert_\zeta-1$, and all test functions $\varphi$ of a given regularity and support in the parabolic ball of radius $1/2$,  uniformly in $0<\epsilon\leq 1$.
\end{prop}

\medskip

\begin{Dem} 
 We check that we can use Theorem \ref{ThmWorkhorse} to evaluate each term $\big\Vert\bbE\big[ \big\langle d^n\overline{\sf \Pi}^\epsilon_z\tau , \varphi_z^\lambda\big\rangle\big]\big\Vert$ by a sum of Feynman graphs each of which satisfies the integrability conditions \eqref{condition_integrability} and the recentering bounds \eqref{condition_recentering}. One has $\vert\tau\vert_\zeta\leq 4$ for all $\tau\in\mcB^-$, so $n\leq 3$ in any case. The case $n=3$ is trivial as it concerns only the trees with four noises, in which case the third Malliavin derivative of $\overline{\sf \Pi}^\epsilon_z\tau$ has only one noise and null expectation.      The same thing happens more generally when $\vert\tau\vert_\zeta-n$ is odd. We are thus left with the cases where $\tau\in\mcB^-$ and $0\leq n\leq 2$ for which $\vert\tau\vert_\zeta-n$ is even, that is $\vert\tau\vert_\zeta\in \{2,4\}$ and $n=0$, or $\vert\tau\vert_\zeta=3$ and $n=1$. The trees with $\vert\tau\vert_\zeta=2$ have degree $-1-2\kappa$ or $-2\kappa$. The trees with $\vert\tau\vert_\zeta=3$ have degree $-1/2-3\kappa$.

\medskip

 \textbf{(a)} We analyse first the case $n=0$ where there is \textbf{no Malliavin derivative}. Recall that
\begin{equation} \label{EqSum}
	\begin{aligned}
 (\overline{\sf \Pi}^\epsilon_z\tau)(y)  =  \left(  (\overline{\sf \Pi}^{\! \epsilon}\cdot)(y) \otimes  (\overline{\sf g}_z^{\!\epsilon})^{-1}(\cdot) \right) \Delta \tau = \sum_{\sigma\leq\tau} (\overline{\sf g}^\epsilon_z)^{-1}(\tau/\sigma)\,(\overline{\sf \Pi}^\epsilon\sigma)(y).
	\end{aligned}
\end{equation}
where we have used the notation
	\begin{equation*}
		\Delta \tau = \sum_{\sigma \leq \tau} \sigma \otimes \tau/\sigma 
	\end{equation*}
where $ \sigma $ runs over the set of all subtrees of $\tau$ with the same root as $\tau$  such that $\tau/\sigma$ has non-negative degree, and there is an implicit sum over the polynomial decorations. The subtree $\sigma$ could have more monomial decorations than $\tau$. By the notation $  \tau/\sigma $, we denote the decorated tree where $ \sigma  $ is contracted to a single node in $ \tau $. Also edge decorations can be removed in this process.   \vspace{0.1cm}

The case $ \sigma = {\bf 1}$ does not happen in the above decomposition of $\tau$ as $\tau/ {\bf 1} = \tau$ must be of non-negative degree  by definition while $\vert\tau\vert<0$. In the sum \eqref{EqSum}, the term $\sigma=\tau$ has null mean from the  defining property \eqref{BPHZ_character} of the BPHZ renormalization map.  This particular case also deals with the trees $\tau$ that have no $\sigma \leq \tau$ other than $\tau$, like $\begin{tikzpicture}[scale=0.3,baseline=0cm]	
	\node at (0,0)  [dot] (0) {};	
	\node at (0.8,0.8)  [noise] (noise1) {};
	\node at (-0.8,0.8)  [noise] (noise2) {};
	\draw[DK] (0) to (noise1);
	\draw[DK] (0) to (noise2);
	\end{tikzpicture}$, $\begin{tikzpicture}[scale=0.3,baseline=0cm]
	\node at (0,0)  [dot] (dot1) {};
	\node at (-0.8,0.8)  [dot] (dot2) {};
	\node at (-1.6,1.6)  [noise] (noise1) {};
	\node at (0,1.6)  [noise] (noise2) {};
	\node at (0.8,0.8)  [noise] (noise3) {};	
	\draw[DK] (dot1) to (dot2);
	\draw[DK] (dot2) to (noise1);
	\draw[DK] (dot2) to (noise2);
	\draw[DK] (dot1) to (noise3);
	\end{tikzpicture}$, $\begin{tikzpicture}[scale=0.3,baseline=0cm]				
	\node at (0,0)  [dot] (dot1) {};
	\node at (-0.8,0.8)  [dot] (dot2) {};	
	\node at (0.8,0.8) [noise] (noise1) {};
	\node at (0,1.6) [noise] (noise2) {};
	\draw[DK] (dot1) to (dot2);
	\draw[DK] (dot1) to (noise1);
	\draw[DK] (dot2) to (noise2);
	\end{tikzpicture}$ or a number of other trees.
 We now concentrate on the case where ${\bf 1}<\sigma<\tau$, so both  $\tau / \sigma$ and $\sigma$ contain at least one noise. 

\medskip

 We aim at estimating $\mathbb{E} \big[ (\overline{\sf g}^\epsilon_z)^{-1}(\tau/\sigma)\,(\overline{\sf \Pi}^\epsilon \sigma)(y) \big]$ using Theorem \ref{ThmWorkhorse} and the propositions that followed. We note as a preliminary remark that we actually have $(\overline{\sf g}^\epsilon_z)^{-1}(\tau/\sigma) = ({\sf g}^\epsilon_z)^{-1}(\tau/\sigma)$ for all $\sigma\leq\tau$ and all $\tau\in\mcB^-$, from the explicit description of $\mcB^-$ in Section \ref{SubsectionRS} and the fact that the constants are in the kernel of the operator $K$. 

\ssk

We deal by hand with trees $\tau$ with two noises using the explicit expression for this expectation and the fact that the constants are in the kernel of the operator $K$. This allows to estimate some quantities of the form 
\[
\int \varphi^\lambda K*\rho^\epsilon*\rho^\epsilon  = \int \big(\varphi^\lambda - \varphi^\lambda(0)\big) K*\rho^\epsilon*\rho^\epsilon
\]
uniformly in $0<\epsilon\leq 1$ by $\lambda^{-1-}$, since $K(z)\simeq \vert z\vert_{\frak{s}}^{-1}$.
 
For trees with four noises, we see from the list of negative trees that there are no cases where $\tau/\sigma$ contains three noises, so $\sigma$ has either two or three noises. The possible trees for $\sigma$ have the form
\begin{equation} \label{EqListTreeForms}
\tau_1 \Xi \mathcal{I}(\Xi \tau_2), \quad \mathcal{I}_1(\Xi \tau_1 ) \mathcal{I}_1(\Xi \tau_2), \quad \mathcal{I}(\Xi \tau_1) \mathcal{I}_1(\Xi \tau_2), \quad \Xi \mathcal{I}(\tau_1 \mathcal{I}_{1}(\Xi \tau_2 )), \quad \mathcal{I}_1(\Xi) \mathcal{I}_1(\tau_1 \mathcal{I}_{1}(\Xi \tau_2))
\end{equation}
where the trees $ \tau_1 $ or $ \tau_2 $ can be the empty tree and the tree $ \tau_1 \tau_2$ contains at most  one noise. Their negative subtrees are one of the trees    
\begin{equation} \label{EqElementaryDivergingTrees}
\begin{split}
(-1)^- &: 
\begin{tikzpicture}[scale=0.3,baseline=0cm]	
	\node at (0,0) [dot] (dot) {};		
	\node at (-0.8,0.8) [noise] (noise1) {};
	\node at (0.8,0.8) [noise] (noise2) {};
	\draw[DK] (dot) to (noise1);
	\draw[DK] (dot) to (noise2);	
\end{tikzpicture}
\hspace{0.3cm}, \hspace{0.3cm}
 \begin{tikzpicture}[scale=0.3,baseline=0cm]
	\node at (0,0)  [noise] (noise1) {};
	\node at (0,1.1)  [noise] (noise2) {};
	\draw[K] (noise1) to (noise2);
\end{tikzpicture}   \\
0^- &:  \begin{tikzpicture}[scale=0.3,baseline=0cm]	
	\node at (0,0) [dot] (dot) {};		
	\node at (-0.8,0.8) [noise] (noise1) {};
	\node at (0.8,0.8) [noise] (noise2) {};
	\draw[K] (dot) to (noise1);
	\draw[DK] (dot) to (noise2);	
\end{tikzpicture}
\hspace{0.3cm}, \hspace{0.3cm}
\begin{tikzpicture}[scale=0.3,baseline=0cm]
	\node at (0,0)  [noise] (noise1) {};
	\node at (0,1)  [dot] (dot) {};
	\node at (-0.8,1.6)  [noise] (noise2) {};	
	\draw[K] (noise1) to (dot);
	\draw[DK] (dot) to (noise2);
	\end{tikzpicture}
	\hspace{0.3cm}, \hspace{0.3cm}
	\begin{tikzpicture}[scale=0.3,baseline=0cm]
	\node at (0,0)  [dot] (dot1) {};
	\node at (-0.8,0.8)  [dot] (dot2) {};
	\node at (0,1.6)  [noise] (noise2) {};
	\node at (0.8,0.8)  [noise] (noise3) {};	
	\draw[DK] (dot1) to (dot2);
	\draw[DK] (dot2) to (noise2);
	\draw[DK] (dot1) to (noise3);
	\end{tikzpicture}
	\hspace{0.3cm}.
\end{split} \end{equation}
The first two trees have degree $(-1)^-$; the last three three have degree $ 0^{-} $. Note that the diverging subtrees in \eqref{EqElementaryDivergingTrees} have at most one incoming edge when seen as a subtree of the tree $\sigma$. 

\ssk

\textbf{Integrability bounds --} We explain the mechanics in the situation where 
$$
\tau = \begin{tikzpicture}[scale=0.3,baseline=0cm]	
	\node at (0,0)  [dot] (0) {};	
	\node at (0.8,0.8)  [noise] (noise1) {};
	\node at (-0.8,0.8)  [noise] (noise2) {};
	\node at (0,1.6)  [noise] (noise3) {};
	\node at (-0.8,2.4)  [noise] (noise4) {};		
	\draw[DK] (0) to (noise1);
	\draw[DK] (0) to (noise2);
	\draw[K] (noise2) to (noise3);
	\draw[K] (noise3) to (noise4);
	\end{tikzpicture}, 
\quad 
\sigma = \begin{tikzpicture}[scale=0.3,baseline=0cm]	
	\node at (0,0)  [dot] (0) {};	
	\node at (0.8,0.8)  [noise] (noise1) {};
	\node at (-0.8,0.8)  [noise] (noise2) {};
	\node at (0,1.6)  [noise] (noise3) {};
	\draw[DK] (0) to (noise1);
	\draw[DK] (0) to (noise2);
	\draw[K] (noise2) to (noise3);
	\end{tikzpicture},
\quad
\tau/\sigma =  \begin{tikzpicture}[scale=0.3,baseline=0cm]
	\node at (0,0)  [dot] (noise1) {};
	\node at (0,1.1)  [noise] (noise2) {};
	\draw[K] (noise1) to (noise2);
\end{tikzpicture}
$$ 
This case corresponds to the second pattern in \eqref{EqListTreeForms} with $\tau_1=\emptyset$ and $\tau_2=\begin{tikzpicture}[scale=0.3,baseline=0cm]
	\node at (0,0)  [dot] (noise1) {};
	\node at (0,1.1)  [noise] (noise2) {};
	\draw[K] (noise1) to (noise2);
\end{tikzpicture}$. The other cases are dealt with in the very same way and we let  the reader repeat the argument for them. 

The expectation of $ (\overline{\sf g}_z^\epsilon)^{-1}(\tau/\sigma)\,  ({\sf \Pi}^\epsilon \sigma)(\varphi_z^{\lambda}) $ is given by
\begin{equation} \label{telescopic_sumç3}
	\mathbb{E} \left[ (\overline{\sf g}^\epsilon_z)^{-1}(\tau/\sigma)\,({\sf \Pi}^\epsilon \sigma)(\varphi_z^{\lambda}) \right] = - \begin{tikzpicture}[scale=0.3,baseline=0cm]
		\node at (0,-1)  [root, label=below:\tiny{$ z$} ] (rootb) {};
		\node at (0,0)  [dot] (root) {};
		\node at (-0.5,2)  [dot] (root1) {};
		\node at (-1,1)  [dot] (dot1) {};
		\node at (1,1)  [dot] (dot2) {};	
		\node at (1,3)  [dot] (dot21) {};
		\node at (3,3)  [dot] (dot22) {};
		\draw[testfcn] (root) to (rootb);
		\draw[K] (dot2) to node[labl, pos=1.24] {\tiny $$} (dot21);
		\draw[K] (dot22) to[bend left=60] node[labl, pos=-0.2] {} (rootb);
		\draw[DK] (root) to node[labl, pos=-0.5] {\tiny $ y $} (dot1);
		\draw[DK] (root) to node {}(dot2);
		\draw[K,purple] (dot1) to (dot2);
		\draw[K,purple] (dot21) to (dot22);
		\draw (root1) to node[labl,pos=0] {\tiny $$} (root1) ;
	\end{tikzpicture} - 2  \begin{tikzpicture}[scale=0.3,baseline=0cm]
	\node at (0,-1)  [root, label=below:\tiny{$ z$} ] (rootb) {};
	\node at (0,0)  [dot] (root) {};
	\node at (-0.5,2)  [dot] (root1) {};
	\node at (-1,1)  [dot] (dot1) {};
	\node at (1,1)  [dot] (dot2) {};	
	\node at (-1,3)  [dot] (dot21) {};
	\node at (3,1)  [dot] (dot22) {};
	\draw[testfcn] (root) to (rootb);
	\draw[K] (dot1) to node[labl, pos=1.24] {\tiny $$} (dot21);
	\draw[K,purple] (dot1) to[bend left=60] node[labl, pos=1.24] {\tiny $$} (dot21);
	\draw[K] (dot22) to[bend left=60] node[labl, pos=-0.2] {} (rootb);
	\draw[DK] (root) to node[labl, pos=-0.5] {\tiny $ y $} (dot1);
	\draw[DK] (root) to node {}(dot2);
	\draw[K,purple] (dot2) to (dot22);
	\draw (root1) to node[labl,pos=0] {\tiny $$} (root1) ;
	\end{tikzpicture}
\end{equation}
In the diagrams above, we have omitted the lables of $r_e$ and $v_e$ when $ (r_e,v_e)= (0,z) $.
The second term is zero due to the fact that one can detached the loop via translation-invariance and then one integrates the constant one to zero via the kernel $K$.
One has to focus on the first term. We use the following short hand notation:
\begin{equation} \label{main_graph}
\begin{tikzpicture}[scale=0.3,baseline=0cm]
	\node at (0,0)  [dot] (root) {};
	\node at (-1,1) [dot] (dot1) {};
	\node at (1,1)  [dot] (dot2) {};	
	\node at (1,3)  [] (dot21) {};
	\draw[K] (dot2) to node[labl, pos=1.24] {\tiny $G$} (dot21);
	\draw[DK] (root) to node[labl, pos=-0.5] {\tiny $ y $} (dot1);
	\draw[DK] (root) to node {}(dot2);
	\draw[K,purple] (dot1) to (dot2);
\end{tikzpicture}
\end{equation}
where $G$ is the graph obtained by removing the green vertex and the green edge.
 We have changed our notation {\color{purple} $\bullet$} for the spacetime convolution operator with ${\color{purple} \rho^\epsilon*\rho^\epsilon}$ for a {\color{purple} purple edge}. The subdivergent diagram in \eqref{main_graph} is given by
 \begin{equation*}
\int  \partial K(y-z_1) \partial K(y-z_2) ({ \color{purple}\rho^{\eps} * \rho^{\eps}}) ( z_1 - z_2 ) d z_1 d z_2
\eqdef
	\begin{tikzpicture}[scale=0.3,baseline=0cm]
 		\node at (0,0)  [dot] (root) {};
 		\node at (-1,1)  [dot] (dot1) {};
 		\node at (1,1)  [dot] (dot2) {};	
 		\draw[DK] (root) to node {}(dot1);
 		\draw[DK] (root) to node {}(dot2);
 		\draw[DK] (root) to node[labl, pos=-0.5] {\tiny $y$} (dot1);
 		\draw[K,purple] (dot1) to (dot2);
 	\end{tikzpicture}.
 \end{equation*} 

We compensate the subdivergence in \eqref{main_graph} by a {\color{orange} Taylor expansion} on some incoming edge. We use a telescoping  sum to make appear this renormalisation
\begin{equation} \label{telescopic_sum}
\begin{tikzpicture}[scale=0.3,baseline=0cm]
	\node at (0,0)  [dot] (root) {};
	\node at (-1,1)  [dot] (dot1) {};
	\node at (1,1)  [dot] (dot2) {};	
	\node at (1,3)  [] (dot21) {};
	\draw[K] (dot2) to node[labl, pos=1.24] {\tiny $G$} (dot21);
	\draw[DK] (root) to node[labl, pos=-0.5] {\tiny $ y $} (dot1);
	\draw[DK] (root) to node {}(dot2);
	\draw[K,purple] (dot1) to (dot2);
\end{tikzpicture} 
= 
\begin{tikzpicture}[scale=0.3,baseline=0cm]
	\node at (0,0)  [dot] (root) {};
	\node at (-0.16 , 2)  [dot] (root1) {};
	\node at (-1,1)  [dot] (dot1) {};
	\node at (1,1)  [dot] (dot2) {};	
	\node at (1,3)  [] (dot21) {};
	\draw[K] (dot2) to node[labl, pos=1.24] {\tiny $G$} (dot21);
	\draw[DK] (root) to node[labl, pos=-0.5] {\tiny $ y $} (dot1);
	\draw[DK] (root) to node {}(dot2);
	\draw[K,purple] (dot1) to (dot2);
	\draw (root1) to node[labl,pos=0] {\tiny $[y,2]$} (root1) ;
\end{tikzpicture}  
+ 
\begin{tikzpicture}[scale=0.3,baseline=0cm]
	\node at (0,0)  [dot] (root) {};
	\node at (-1,1)  [dot] (dot1) {};
	\node at (1,1)  [dot] (dot2) {};	
	\draw[DK] (root) to node {}(dot1);
	\draw[DK] (root) to node {}(dot2);
	\draw[DK] (root) to node[labl, pos=-0.5] {\tiny $y$} (dot1);
	\draw[K,purple] (dot1) to (dot2);
\end{tikzpicture} \, \,\, 
\begin{tikzpicture}[scale=0.3,baseline=0cm]
	\node at (0,0)  [dot] (root) {};
	\node at (0,2)  [dot] (dot1) {};
	\draw[K] (root) to node[labl, pos=1.24] {\tiny $G$} (dot1);
	\draw to node[labl, pos=-0.5] {\tiny $y$} (root) ;
\end{tikzpicture}	 
+ 	
\begin{tikzpicture}[scale=0.3,baseline=0cm]
\node at (0,0)  [dot] (root) {};
\node at (-1,1)  [dot] (dot1) {};
\node at (1,1)  [dot] (dot2) {};	
\draw[DK] (root) to node {}(dot1);
\draw[DK] (root) to node {}(dot2);
\draw[DK] (root) to node[labl, pos=-0.5] {\tiny $y$} (dot1);
\draw[K,purple] (dot1) to (dot2);
\draw[K,orange] (dot2) to [bend left=+50] (root);
\end{tikzpicture} \, \,\, 
\begin{tikzpicture}[scale=0.3,baseline=0cm]
	\node at (0,0)  [dot] (root) {};
	\node at (0,2)  [dot] (dot1) {};
	\draw[DK] (root) to node[labl, pos=1.24] {\tiny $G$} (dot1);
	\draw to node[labl, pos=-0.5] {\tiny $y$} (root) ;
\end{tikzpicture} 
\end{equation}
The first term is well defined as the decoration $[y,2]$ means a kernel of the form
$$
K(z_1-z_3) - K(y-z_3) - (x_1-x) \partial K(y-z_3)
$$
that behaves like $ (x_1-x)^2 $ when $ z_1 $ is close to $y=(t,x)$.  (The integration variable $z_3$ corresponds to $G$ and the integration variable $z_1$ corresponds to the upper right of the triangle diagram.) One also writes
\begin{equation*}
	\begin{tikzpicture}[scale=0.3,baseline=0cm]
		\node at (0,0)  [dot] (root) {};
		\node at (-1,1)  [dot] (dot1) {};
		\node at (1,1)  [dot] (dot2) {};	
		\draw[DK] (root) to node {}(dot1);
		\draw[DK] (root) to node {}(dot2);
		\draw[DK] (root) to node[labl, pos=-0.5] {\tiny $y$} (dot1);
		\draw[K,purple] (dot1) to (dot2);
		\draw[K,orange] (dot2) to [bend left=+50] (root);
	\end{tikzpicture} \defeq \int {\color{orange} (x_1-x)}  \partial K(y-z_1) \partial K(y-z_2) ({ \color{purple}\rho^{\eps} * \rho^{\eps}}) ( z_1 - z_2 ) d z_1 d z_2.
\end{equation*}
For the first term of \eqref{telescopic_sum} one can check the integrability condition~\eqref{condition_integrability}. Now the two problematic terms are the last two terms of \eqref{telescopic_sum}. They are removed by the renormalization. Indeed one has 
\begin{equation*}
	R(y)(\sigma) = \sigma + \ell(
	\begin{tikzpicture}[scale=0.3,baseline=0cm]
		\node at (0,0)  [dot] (root) {};
		\node at (-0.7,0.7)  [noise] (dot1) {};
		\node at (0.7,0.7)  [noise] (dot2) {};	
		\draw[DK] (root) to node {}(dot1);
		\draw[DK] (root) to node {}(dot2);
		\draw[DK] (root) to node {} (dot1);
	\end{tikzpicture} \,
	)(y) \, \begin{tikzpicture}[scale=0.3,baseline=0cm]
	\node at (0,0)  [dot] (noise1) {};
	\node at (0,1.1)  [noise] (noise2) {};
	\draw[K] (noise1) to (noise2);
\end{tikzpicture} 
+ 
\ell(
\begin{tikzpicture}[scale=0.3,baseline=0cm]	
	\node at (0,0)  [dot] (0) {};	
	\node at (0.8,0.8)  [dot] (noise1) {};
	\node at (-0.8,0.8)  [noise] (noise2) {};	
	\node at (0.8,0.8)  [noiseblue] (1) {};
	\draw[DK] (0) to (noise1);
	\draw[DK] (0) to (noise2);
	\end{tikzpicture}
)(y) \, 
\begin{tikzpicture}[scale=0.3,baseline=0cm]
	\node at (0,0)  [dot] (noise1) {};
	\node at (0,1.1)  [noise] (noise2) {};
	\draw[DK] (noise1) to (noise2);
\end{tikzpicture} 
\end{equation*}
by using the definition \eqref{def_R} of $R_{\ell}$. Then we conclude by noticing that\footnote{We kept track of the point $ y $ in the above computations to stress that the formulation is robust in the non-translation invariant case. In this case we obtain renormalization functions instead of renormalization constants. We could have removed the dependency in $ y $.}
\begin{equation*}
	\ell(
	\begin{tikzpicture}[scale=0.3,baseline=0cm]
		\node at (0,0)  [dot] (root) {};
		\node at (-0.7,0.7)  [noise] (dot1) {};
		\node at (0.7,0.7)  [noise] (dot2) {};	
		\draw[DK] (root) to node {}(dot1);
		\draw[DK] (root) to node {}(dot2);
		\draw[DK] (root) to node {} (dot1);
	\end{tikzpicture} \,
	)(y) = - \begin{tikzpicture}[scale=0.3,baseline=0cm]
		\node at (0,0)  [dot] (root) {};
		\node at (-1,1)  [dot] (dot1) {};
		\node at (1,1)  [dot] (dot2) {};	
		\draw[DK] (root) to node {}(dot1);
		\draw[DK] (root) to node {}(dot2);
		\draw[DK] (root) to node[labl, pos=-0.5] {\tiny $y$} (dot1);
		\draw[K,purple] (dot1) to (dot2);
	\end{tikzpicture}, 
\quad
\ell(
\begin{tikzpicture}[scale=0.3,baseline=0cm]	
	\node at (0,0)  [dot] (0) {};	
	\node at (0.8,0.8)  [dot] (noise1) {};
	\node at (-0.8,0.8)  [noise] (noise2) {};	
	\node at (0.8,0.8)  [noiseblue] (1) {};
	\draw[DK] (0) to (noise1);
	\draw[DK] (0) to (noise2);
	\end{tikzpicture}
)(y)
= - 
\begin{tikzpicture}[scale=0.3,baseline=0cm]
	\node at (0,0)  [dot] (root) {};
	\node at (-1,1)  [dot] (dot1) {};
	\node at (1,1)  [dot] (dot2) {};	
	\draw[DK] (root) to node {}(dot1);
	\draw[DK] (root) to node {}(dot2);
	\draw[DK] (root) to node[labl, pos=-0.5] {\tiny $y$} (dot1);
	\draw[K,purple] (dot1) to (dot2);
	\draw[K,orange] (dot2) to [bend left=+50] (root);
\end{tikzpicture}
\end{equation*}
One can observe that the telescoping sums are a local procedure that detaches incoming edges from a diverging subgraph associated to a subtree of negative degree. These edges are therefore attached to the root of the diverging subtree and then we obtain a diverging constant that can be detached if one has translation-invariance. Moreover one has also the same changes in the decorations provided by the order of the telescoping sum. One gets in the end a one to one correspondence with the procedure of extraction-contraction of a subtree given by $ R_{\ell} $. This has been previously observed in \cite[Sec. 5.3]{BN1}. In the end on the example considered, one has
\begin{equation*} 
	 \label{telescopic_sumç3}
		\mathbb{E} \left[ (\overline{\sf g}^\epsilon_z)^{-1}(\tau/\sigma)\,(\overline{{\sf \Pi}}^\epsilon \sigma)(\varphi_z^{\lambda}) \right] =  \begin{tikzpicture}[scale=0.3,baseline=0cm]
			\node at (0,-1)  [root, label=below:\tiny{$ z$} ] (rootb) {};
			\node at (0.5,2)  [dot] (root1) {};
			\node at (0,0)  [dot] (root) {};
			\node at (-0.5,2)  [dot] (root1) {};
			\node at (-1,1)  [dot] (dot1) {};
			\node at (1,1)  [dot] (dot2) {};	
			\node at (1,3)  [dot] (dot21) {};
			\node at (3,3)  [dot] (dot22) {};
			\draw[testfcn] (root) to (rootb);
			\draw[K] (dot2) to node[labl, pos=1.24] {\tiny $$} (dot21);
			\draw[K] (dot22) to[bend left=60] node[labl, pos=-0.2] {} (rootb);
			\draw[DK] (root) to node[labl, pos=-0.5] {\tiny $ y $} (dot1);
			\draw[DK] (root) to node {}(dot2);
			\draw[K,purple] (dot1) to (dot2);
			\draw[K,purple] (dot21) to (dot22);
				\draw (root1) to node[labl,pos=0] {\tiny $[y,2]$} (root1) ;
			\draw (root1) to node[labl,pos=0] {\tiny $$} (root1) ;
		\end{tikzpicture}.
	\end{equation*}

\ssk

\textbf{Recentering bounds --} One can check by direct inspection that the recentering bound \eqref{condition_recentering} is satisfied for 
\[
\mathbb{E} \left[ (\overline{\sf g}^\epsilon_z)^{-1}(\tau/\sigma)\,(\overline{\sf \Pi}^\epsilon \sigma)(\varphi_z^{\lambda}) \right].
\]
One can check by hand that the bound \eqref{condition_recentering} is satisfied. 
One proceeds then by direct inspection on trees with two and three noises given by the terms $\sigma$. Then, one does the same when $\tau / \sigma$ has three noises which is quite similar to $\sigma$ except that now one can have a recentering on the edges. We consider the case
\begin{equation*}
	\tau = \begin{tikzpicture}[scale=0.3,baseline=0cm]	
		\node at (0,0) [noise] (dot) {};		
		\node at (-0.8,0.8) [noise] (noise1) {};
		\node at (0,1.6) [noise] (noise2) {};
		\node at (-0.8,2.4) [noise] (noise3) {};
		\draw[K] (dot) to (noise1);
		\draw[K] (noise1) to (noise2);
		\draw[K] (noise2) to (noise3);
	\end{tikzpicture}, 
	\quad 
	\sigma = \begin{tikzpicture}[scale=1,baseline=-0.1cm] \node at (0,0)  [noise] (1) {};\end{tikzpicture},
	\quad
	\tau/\sigma =  \begin{tikzpicture}[scale=0.3,baseline=0cm]	
		\node at (0,0) [dot] (dot) {};		
		\node at (-0.8,0.8) [noise] (noise1) {};
		\node at (0,1.6) [noise] (noise2) {};
		\node at (-0.8,2.4) [noise] (noise3) {};
		\draw[K] (dot) to (noise1);
		\draw[K] (noise1) to (noise2);
		\draw[K] (noise2) to (noise3);
	\end{tikzpicture}
\end{equation*} 
Then, one has 
\begin{equation*}
	\mathbb{E} \left[ ({\sf g}^\epsilon_z)^{-1}(\tau/\sigma)\,(\overline{{\sf \Pi}}^\epsilon \sigma)(\varphi_z^{\lambda}) \right] = - \begin{tikzpicture}[scale=0.3,baseline=0cm]
		\node at (0,-1)  [root, label=below:\tiny{$ z$} ] (root) {};
		\node at (0,1)  [dot] (root1) {};
		\node at (0,3)  [dot] (root2) {};
		\node at (0,5)  [dot] (root3) {};
			\node at (2,1)  [dot] (root4) {};
				\node at (-1.5,4)  [] (root3b) {\tiny $[z,1]$};
				\node at (-1.5,2)  [] (root2b) {\tiny $[z,1]$};
		\draw[K] (root) to (root1);
		\draw[K] (root1) to (root2);
		\draw[K] (root2) to (root3);
		\draw[K,purple] (root1) to [bend right=60] (root2);
		\draw[K,purple] (root4) to [bend right=60] (root3);
		\draw[testfcn] (root) to  (root4);
	\end{tikzpicture} - \begin{tikzpicture}[scale=0.3,baseline=0cm]
	\node at (0,-1)  [root, label=below:\tiny{$ z$} ] (root) {};
	\node at (0,1)  [dot] (root1) {};
	\node at (0,3)  [dot] (root2) {};
	\node at (0,5)  [dot] (root3) {};
	\node at (2,1)  [dot] (root4) {};
	\node at (-1.5,4)  [] (root3b) {\tiny $[z,1]$};
	\node at (-1.5,2)  [] (root2b) {\tiny $[z,1]$};
	\draw[K] (root) to (root1);
	\draw[K] (root1) to (root2);
	\draw[K] (root2) to (root3);
	\draw[K,purple] (root2) to [bend right=60] (root3);
	\draw[K,purple] (root1) to  (root4);
	\draw[testfcn] (root) to  (root4);
	\end{tikzpicture} - \begin{tikzpicture}[scale=0.3,baseline=0cm]
	\node at (0,-1)  [root, label=below:\tiny{$ z$} ] (root) {};
	\node at (0,1)  [dot] (root1) {};
	\node at (0,3)  [dot] (root2) {};
	\node at (0,5)  [dot] (root3) {};
	\node at (2,1)  [dot] (root4) {};
	\node at (1,4)  [] (root3b) {\tiny $[z,1]$};
	\node at (1,1.5)  [] (root2b) {\tiny $[z,1]$};
	\draw[K] (root) to (root1);
	\draw[K] (root1) to (root2);
	\draw[K] (root2) to (root3);
	\draw[K,purple] (root1) to [bend left=60] (root3);
	\draw[K,purple] (root4) to [bend right=60] (root2);
	\draw[testfcn] (root) to  (root4);
	\end{tikzpicture}.
\end{equation*}
 One has a recentering $[z,1]$ on one edge connecting the node $z_0$. This is due to the definition of the model for $ \tau = \begin{tikzpicture}[scale=0.3,baseline=0cm]	
		\node at (0,0) [noise] (dot) {};		
		\node at (-0.8,0.8) [noise] (noise1) {};
		\node at (0,1.6) [noise] (noise2) {};
		\node at (-0.8,2.4) [noise] (noise3) {};
		\draw[K] (dot) to (noise1);
		\draw[K] (noise1) to (noise2);
		\draw[K] (noise2) to (noise3);
	\end{tikzpicture}$. Indeed, the subtree $	\begin{tikzpicture}[scale=0.3,baseline=0cm]
		\node at (0,0)  [dot] (1) {};
		\node at (0,1.1)  [noise] (2) {};
		\node at (0,2.2)  [noise] (3) {};	
		\draw[K] (1) to (2);
		\draw[K] (2) to (3);	
	\end{tikzpicture}$ is of degree $ 1 - 2 \kappa $ and the subtree $ \begin{tikzpicture}[scale=0.3,baseline=0cm]
		\node at (0,0)  [dot] (dot) {};
		\node at (0,1)  [noise] (noise) {};
		\draw[K] (dot) to (noise);
	\end{tikzpicture} $ is of degree $\frac{1}{2} - \kappa$, with $\kappa \in (0,1)$.
We will consider the first diagram and use the following short hand notation:
\begin{equation} \label{tree_loop} \begin{tikzpicture}[scale=0.3,baseline=0cm]
		\node at (0,1)  [dot, label=below:\tiny{$ z_0$} ] (root1) {};
		\node at (0,3)  [dot,label=left:\tiny{$ z_1$}] (root2) {};
		\node at (0,5)  [dot,label=above:\tiny{$ G$},label=left:\tiny{$ z_2$} ] (root3) {};
		\node at (-1.5,4)  [] (root3b) {\tiny $[z,1]$};
		\node at (-1.5,2)  [] (root2b) {\tiny $[z,1]$};
		\draw[K] (root1) to (root2);
		\draw[K] (root2) to (root3);
		\draw[K,purple] (root1) to [bend right=60] (root2);
	\end{tikzpicture}
\end{equation}
where $G$ contains all the other nodes and edges except the edge connecting $z_0$ to $z$. We focus on the part of the diagram where one has a subdivergence given by the loop with the purple edge. One notices that there is a recentering on the black edge entering the loop and connected to $G$.
One can adapt the telescoping sums method to this new case. For example, one has
 \begin{equation*}
 	\begin{aligned}
 		K(z_1-z_2) - K(z-z_2) &= K(z_1-z_2) - K(z_0-z_2) - (x_1 - x_0) \partial K(z_0 - z_2)   \\ 
 		&\quad+ K(z_0-z_2) - K(z-z_2) + (x_1 - x_0) \partial K(z_0 - z_2)
 	\end{aligned} 
 \end{equation*}
where $ z_0 = (t_0,x_0) $.  The extra term $ - K(z-z_2) $ is coming from the definition \eqref{recentering_Xi} of the model. Now \eqref{tree_loop} becomes
 \begin{equation} \label{telescopic_sum_bis}
 	\begin{tikzpicture}[scale=0.3,baseline=0cm]
 		\node at (0,1)  [dot, label=below:\tiny{$ z_0$} ] (root1) {};
 		\node at (0,3)  [dot,label=left:\tiny{$ z_1$}] (root2) {};
 		\node at (0,5)  [dot,label=above:\tiny{$ G$},label=left:\tiny{$ z_2$} ] (root3) {};
 		\node at (-1.5,4)  [] (root3b) {\tiny $[z,1]$};
 		\node at (-1.5,2)  [] (root2b) {\tiny $[z,1]$};
 		\draw[K] (root1) to (root2);
 		\draw[K] (root2) to (root3);
 		\draw[K,purple] (root1) to [bend right=60] (root2);
 	\end{tikzpicture} =
 	 \begin{tikzpicture}[scale=0.3,baseline=0cm]
 	 	\node at (0,1)  [dot, label=below:\tiny{$ z_0$} ] (root1) {};
 	 	\node at (0,3)  [dot,label=left:\tiny{$ z_1$}] (root2) {};
 	 	\node at (0,5)  [dot,label=above:\tiny{$ G$},label=left:\tiny{$ z_2$} ] (root3) {};
 	 	\node at (-1.5,4)  [] (root3b) {\tiny $[z_0,2]$};
 	 	\node at (-1.5,2)  [] (root2b) {\tiny $[z,1]$};
 	 	\draw[K] (root1) to (root2);
 	 	\draw[K] (root2) to (root3);
 	 	\draw[K,purple] (root1) to [bend right=60] (root2);
 	 \end{tikzpicture} +  \begin{tikzpicture}[scale=0.3,baseline=0cm]
 	 \node at (0,1)  [dot, label=below:\tiny{$ z_0$} ] (root1) {};
 	 \node at (0,3)  [dot,label=left:\tiny{$ z_1$}] (root2) {};
 	 \node at (-1.5,2)  [] (root2b) {\tiny $[z,1]$};
 	 \draw[K] (root1) to (root2);
 	 \draw[K,purple] (root1) to [bend right=60] (root2);
 	 \end{tikzpicture}  \, \,\, 
 	\begin{tikzpicture}[scale=0.3,baseline=0cm]
 		\node at (0,0)  [dot] (root) {};
 			\node at (1.5,1)  [dot] (root1) {};
 		\node at (0,2)  [dot] (dot1) {};
 		\draw[K] (root) to node[labl, pos=1.24] {\tiny $G$} (dot1);
 		\draw to node[labl, pos=-0.5] {\tiny $z_0$} (root) ;
 		\draw (root1) to node[labl,pos=0] {\tiny $[z,1]$} (root1) ;
 	\end{tikzpicture}	 + 	\begin{tikzpicture}[scale=0.3,baseline=0cm]
 	\node at (0,1)  [dot, label=below:\tiny{$ z_0$} ] (root1) {};
 	\node at (0,3)  [dot,label=left:\tiny{$ z_1$}] (root2) {};
 	\node at (-1.5,2)  [] (root2b) {\tiny $[z,1]$};
 	\draw[K] (root1) to (root2);
 	\draw[K,purple] (root1) to [bend right=60] (root2);
 	\draw[K,orange] (root1) to [bend left=60] (root2);
 	\end{tikzpicture}\, \,\, 
 	\begin{tikzpicture}[scale=0.3,baseline=0cm]
 		\node at (0,0)  [dot] (root) {};
 		\node at (0,2)  [dot] (dot1) {};
 		\draw[DK] (root) to node[labl, pos=1.24] {\tiny $G$} (dot1);
 		\draw to node[labl, pos=-0.5] {\tiny $z_0$} (root) ;
 	\end{tikzpicture} 
 \end{equation}
 From Proposition \ref{PropPreservedRecenteringBounds} one sees that the recentering bounds \eqref{condition_recentering} are preserved for all the terms of the right hand side of \eqref{telescopic_sum_bis}.  In fact, one observes that the telescoping sums method for performing the renormalisation is compatible with the recentering provided by the model.  An issue seems to be the presence of some recentering inside what we expect to be renormalisation constants. One has the following decomposition:
	\begin{equation*}
		\begin{aligned}
		\begin{tikzpicture}[scale=0.3,baseline=0cm]
			\node at (0,1)  [dot, label=below:\tiny{$ z_0$} ] (root1) {};
			\node at (0,3)  [dot,label=left:\tiny{$ z_1$}] (root2) {};
			\node at (-1.5,2)  [] (root2b) {\tiny $[z,1]$};
			\draw[K] (root1) to (root2);
			\draw[K,purple] (root1) to [bend right=60] (root2);
		\end{tikzpicture}   & = \int   K(z_0-z_1) \left(   K(z_0-z_2) -  K(z-z_2)\right) ({ \color{purple}\rho^{\eps} * \rho^{\eps}}) ( z_1 - z_0 ) d z_1 .
		\\ & = \int  K(z_0-z_1)   K(z_0-z_2) ({ \color{purple}\rho^{\eps} * \rho^{\eps}}) ( z_1 - z_0 ) d z_1 \\ & -  \int   K(z_0-z_1)   K(z-z_2) ({ \color{purple}\rho^{\eps} * \rho^{\eps}}) ( z_1 - z_0 ) d z_1.
		\\ & =
		 		\begin{tikzpicture}[scale=0.3,baseline=0cm]
		 	\node at (0,1)  [dot, label=below:\tiny{$ z_0$} ] (root1) {};
		 	\node at (0,3)  [dot,label=left:\tiny{$ z_1$}] (root2) {};
		 	\draw[K] (root1) to (root2);
		 	\draw[K,purple] (root1) to [bend right=60] (root2);
		 	\end{tikzpicture} -	\begin{tikzpicture}[scale=0.3,baseline=0cm]
		 	\node at (0,1)  [dot, label=below:\tiny{$ z_0$} ] (root1) {};
		 	\node at (-2,1)  [dot, label=below:\tiny{$ z$} ] (root1b) {};
		 	\node at (0,3)  [dot,label=left:\tiny{$ z_1$}] (root2) {};
		 	\draw[K] (root1b) to (root2);
		 	\draw[K,purple] (root1) to [bend right=60] (root2);
		 	\end{tikzpicture}.
		\end{aligned}
	\end{equation*}
The condition \eqref{condition_recentering} is easily satisfied on the extra term with the node $z$ -- see Proposition 5.10 in \cite{BN1}. This term is actually not divergent.   

\ssk

\textbf{(b)} We analyse the case $n=1$ or $n=2$ where there is one or two Malliavin derivatives. The reasoning of point {\bf (a)} works with some minor modifications as the noise derivative operators $D_{\Xi_j}$ from Section \ref{SubsectionMalliavin} commute with the model from Lemma \ref{lemma_commutation}.  When $n=1$ and $\vert\tau\vert_\zeta=3$, the divergences are of same nature as in the case $ n=0 $ and $ \tau \neq \sigma $. The only difference is that we work with a mirror graph as a noise $ \Xi $ has been replaced by its Malliavin derivative. 

When two Malliavin derivatives are involved and $\vert\tau\vert_\zeta=4$, the trees  concerned have the form \eqref{EqListTreeForms} but now  with $\tau_1 \tau_2$ a product of at most two planted trees. Note that these diverging diagrams have at most two incoming edges.   \vspace{0.1cm}

We now explain on the example of the tree $ \mathcal{I}_1(\Xi \tau_1 ) \mathcal{I}_1(\Xi \tau_2) $ how such a renormalization works via a mere telescoping sum. One deals with the other terms of the list \eqref{EqListTreeForms} in exactly the same way. 

One gets the following diagram for $\bbE\big[ ( d^2 {\sf \Pi}^\epsilon \mathcal{I}_1(\Xi \tau_1 ) \mathcal{I}_1(\Xi \tau_2))(h_1,h_2)(y) \big]$

\begin{equation} \label{Hairer_Quastel_reno}
\begin{tikzpicture}[scale=0.3,baseline=0cm]
	\node at (0,0)  [dot] (root) {};
	\node at (-1,1)  [dot] (dot1) {};
	\node at (1,1)  [dot] (dot2) {};	
	\node at (-1,3)  [] (dot11) {};
	\node at (1,3)  [] (dot21) {};
	\draw[K] (dot1) to node[labl, pos=1.24] {\tiny $h_1$} (dot11);
	\draw[K] (dot2) to node[labl, pos=1.24] {\tiny $h_2$} (dot21);
	\draw[DK] (root) to node[labl, pos=-0.5] {\tiny $y$} (dot1);
	\draw[DK] (root) to node {}(dot2);
	\draw[K,purple] (dot1) to (dot2);
\end{tikzpicture}
\end{equation}
where here, without loss of generality, we have supposed that $ \tau_1 = \mathcal{I}(\Xi_1) $ and $ \tau_2 = \mathcal{I}(\Xi_2) $.  (The situation where one $h_i$ replaces one of the noises $\xi$ in $\mathcal{I}_1(\Xi \tau_1 ) \mathcal{I}_1(\Xi \tau_2))$ are dealt by elementary means.) We recognize the same sub-divergent diagram  as in {\bf (a)} but this time we have two incoming edges.  It is this type of graphs that is not covered by the Hairer-Quastel theorem from \cite[Theorem A.3]{HairerQuastel}. We proceed with several telescoping sums as below
 \begin{equation} \label{decomp_1}
 	\begin{tikzpicture}[scale=0.3,baseline=0cm]
 		\node at (0,0)  [dot] (root) {};
 		\node at (-1,1)  [dot] (dot1) {};
 		\node at (1,1)  [dot] (dot2) {};	
 		\node at (-1,3)  [] (dot11) {};
 		\node at (1,3)  [] (dot21) {};
 		\draw[K] (dot1) to node[labl, pos=1.24] {\tiny $h_1$} (dot11);
 		\draw[K] (dot2) to node[labl, pos=1.24] {\tiny $h_2$} (dot21);
 		\draw[DK] (root) to node[labl, pos=-0.5] {\tiny $y$} (dot1);
 		\draw[DK] (root) to node {}(dot2);
 		\draw[K,purple] (dot1) to (dot2);
 	\end{tikzpicture} = 	\begin{tikzpicture}[scale=0.3,baseline=0cm]
 	\node at (0,0)  [dot] (root) {};
 		\node at (-2.5,2)  [dot] (root1) {};
 	\node at (-1,1)  [dot] (dot1) {};
 	\node at (1,1)  [dot] (dot2) {};	
 	\node at (-1,3)  [] (dot11) {};
 	\node at (1,3)  [] (dot21) {};
 	\draw[K] (dot1) to node[labl, pos=1.24] {\tiny $h_1$} (dot11);
 	\draw[K] (dot2) to node[labl, pos=1.24] {\tiny $h_2$} (dot21);
 	\draw[DK] (root) to node[labl, pos=-0.5] {\tiny $y$} (dot1);
 	\draw[DK] (root) to node {}(dot2);
 	\draw[K,purple] (dot1) to (dot2);
 	\draw (root1) to node[labl,pos=0] {\hspace{0.1cm}\tiny $(z_2,2)$} (root1) ;
 \end{tikzpicture} +
\begin{tikzpicture}[scale=0.3,baseline=0cm]
	\node at (0,0)  [dot] (root) {};
	\node at (-1,1)  [dot] (dot1) {};
	\node at (1,1)  [dot] (dot2) {};	
	\node at (0,3)  [] (dot11) {};
	\node at (2,3)  [] (dot21) {};
	\draw[K] (dot2) to node[labl, pos=1.24] {\tiny $h_1$} (dot11);
	\draw[K] (dot2) to node[labl, pos=1.24] {\tiny $h_2$} (dot21);
	\draw[DK] (root) to node[labl, pos=-0.5] {\tiny $y$} (dot1);
	\draw[DK] (root) to node {}(dot2);
	\draw[K,purple] (dot1) to (dot2);
\end{tikzpicture} +
 \begin{tikzpicture}[scale=0.3,baseline=0cm]
	\node at (0,0)  [dot] (root) {};
	\node at (-1,1)  [dot] (dot1) {};
	\node at (1,1)  [dot] (dot2) {};	
	\node at (0,3)  [] (dot11) {};
	\node at (2,3)  [] (dot21) {};
	\draw[DK] (dot2) to node[labl, pos=1.24] {\tiny $h_1$} (dot11);
	\draw[K] (dot2) to node[labl, pos=1.24] {\tiny $h_2$} (dot21);
	\draw[DK] (root) to node[labl, pos=-0.5] {\tiny $y$} (dot1);
	\draw[DK] (root) to node {}(dot2);
	\draw[K,purple] (dot1) to (dot2);
	\draw[K,orange] (dot1) to [bend left=+50] (dot2);
\end{tikzpicture}
 \end{equation}
The first term is well-defined as the edge decoration $ (z_2,2) $ means a kernel of the form
$$
K(z_1-z_3) - K(z_2-z_3) - (x_1-x_2) \partial K(z_2-z_3)
$$
that behaves like $ (x_1-x_2)^2 $ when $ z_1 $ is close to $ z_2 $. This is sufficient for making sense of the first term in \eqref{decomp_1} when the $\epsilon$-regularization is removed. More precisely, one can check the bounds \eqref{condition_integrability}. The {\color{orange} orange line} is encoding a multiplicative term of the form $ (x_2-x_1) $
 \begin{equation*}
	\begin{tikzpicture}[scale=0.3,baseline=0cm]
		\node at (0,0)  [dot] (root) {};
		\node at (-1,1)  [dot] (dot1) {};
		\node at (1,1)  [dot] (dot2) {};	
		\draw[DK] (root) to node {}(dot1);
		\draw[DK] (root) to node {}(dot2);
		\draw[DK] (root) to node[labl, pos=-0.5] {\tiny $y$} (dot1);
		\draw[K,purple] (dot1) to (dot2);
		\draw[K,orange] (dot1) to [bend left=+50] (dot2);
	\end{tikzpicture} = \int (x_2-x_1)  \partial K(y-z_1) \partial K(y-z_2) ({\color{purple} \rho^{\eps} * \rho^{\eps}}) ( z_1 - z_2 ) d z_2 d z_1.
\end{equation*}
Note that this term is equal to zero when the mollifier $ \rho $ is symmetric. The last two terms in \eqref{decomp_1} are still ill-defined and one has to move the two edges to the root. One gets
\begin{equation*}
	\begin{tikzpicture}[scale=0.3,baseline=0cm]
		\node at (0,0)  [dot] (root) {};
		\node at (-1,1)  [dot] (dot1) {};
		\node at (1,1)  [dot] (dot2) {};	
		\node at (0,3)  [] (dot11) {};
		\node at (2,3)  [] (dot21) {};
		\draw[K] (dot2) to node[labl, pos=1.24] {\tiny $h_1$} (dot11);
		\draw[K] (dot2) to node[labl, pos=1.24] {\tiny $h_2$} (dot21);
		\draw[DK] (root) to node[labl, pos=-0.5] {\tiny $y$} (dot1);
		\draw[DK] (root) to node {}(dot2);
		\draw[K,purple] (dot1) to (dot2);
	\end{tikzpicture}  =  \begin{tikzpicture}[scale=0.3,baseline=0cm]
	\node at (0,0)  [dot] (root) {};
	\node at (-0.5,2)  (root1) {};
	\node at (-1,1)  [dot] (dot1) {};
	\node at (1,1)  [dot] (dot2) {};	
	\node at (0,3)  [] (dot11) {};
	\node at (2,3)  [] (dot21) {};
	\draw[K] (dot2) to node[labl, pos=1.24] {\tiny $h_1$} (dot11);
	\draw[K] (dot2) to node[labl, pos=1.24] {\tiny $h_2$} (dot21);
	\draw[DK] (root) to node[labl, pos=-0.5] {\tiny $y$} (dot1);
	\draw[DK] (root) to node {}(dot2);
	\draw[K,purple] (dot1) to (dot2);
	\draw (root1) to node[labl,pos=-1] {\hspace{-0.5cm}\tiny $(y,2)$} (root1) ;
\end{tikzpicture} +  \begin{tikzpicture}[scale=0.3,baseline=0cm]
	\node at (0,0)  [dot] (root) {};
	\node at (-1,1)  [dot] (dot1) {};
	\node at (1,1)  [dot] (dot2) {};	
	\node at (-2,0)  [] (dot11) {};
	\node at (2,3)  [] (dot21) {};
	\draw[K] (root) to node[labl, pos=1.24] {\tiny $h_1$} (dot11);
	\draw[K] (dot2) to node[labl, pos=1.24] {\tiny $h_2$} (dot21);
	\draw[DK] (root) to node[labl, pos=-0.5] {\tiny $y$} (dot1);
	\draw[DK] (root) to node {}(dot2);
	\draw[K,purple] (dot1) to (dot2);
\end{tikzpicture} 
+ \begin{tikzpicture}[scale=0.3,baseline=0cm]
	\node at (0,0)  [dot] (root) {};
	\node at (-1,1)  [dot] (dot1) {};
	\node at (1,1)  [dot] (dot2) {};	
	\node at (-2,0)  [] (dot11) {};
	\node at (2,3)  [] (dot21) {};
	\draw[DK] (root) to node[labl, pos=1.24] {\tiny $h_1$} (dot11);
	\draw[K] (dot2) to node[labl, pos=1.24] {\tiny $h_2$} (dot21);
	\draw[DK] (root) to node[labl, pos=-0.5] {\tiny $y$} (dot1);
	\draw[DK] (root) to node {}(dot2);
	\draw[K,purple] (dot1) to (dot2);
	\draw[K,orange] (dot2) to [bend left=+50] (root);
\end{tikzpicture} 
	\end{equation*}
Then we apply the same reduction to the last term and obtain
	\begin{equation*}
	\begin{tikzpicture}[scale=0.3,baseline=0cm]
	\node at (0,0)  [dot] (root) {};
	\node at (-1,1)  [dot] (dot1) {};
	\node at (1,1)  [dot] (dot2) {};	
	\node at (-2,0)  [] (dot11) {};
	\node at (2,3)  [] (dot21) {};
	\draw[DK] (root) to node[labl, pos=1.24] {\tiny $h_1$} (dot11);
	\draw[K] (dot2) to node[labl, pos=1.24] {\tiny $h_2$} (dot21);
	\draw[DK] (root) to node[labl, pos=-0.5] {\tiny $y$} (dot1);
	\draw[DK] (root) to node {}(dot2);
	\draw[K,purple] (dot1) to (dot2);
	\draw[DK,orange] (dot2) to [bend left=+50] (root);
\end{tikzpicture} = \begin{tikzpicture}[scale=0.3,baseline=0cm]
	\node at (0,0)  [dot] (root) {};
	\node at (-1,1)  [dot] (dot1) {};
	\node at (1,1)  [dot] (dot2) {};	
	\node at (-2,0)  [] (dot11) {};
	\node at (2,3)  [] (dot21) {};
	\draw[DK] (root) to node[labl, pos=1.24] {\tiny $h_1$} (dot11);
	\draw[K] (dot2) to node[labl, pos=1.24] {\tiny $h_2$} (dot21);
	\draw[DK] (root) to node[labl, pos=-0.5] {\tiny $y$} (dot1);
	\draw[DK] (root) to node {}(dot2);
	\draw[K,purple] (dot1) to (dot2);
	\draw[K,orange] (dot2) to [bend left=+50] (root);
	\draw (root1) to node[labl,pos=0] {\tiny $(y,1)$} (root1) ;
\end{tikzpicture} + \begin{tikzpicture}[scale=0.3,baseline=0cm]
\node at (0,0)  [dot] (root) {};
\node at (-1,1)  [dot] (dot1) {};
\node at (1,1)  [dot] (dot2) {};	
\draw[DK] (root) to node {}(dot1);
\draw[DK] (root) to node {}(dot2);
\draw[DK] (root) to node[labl, pos=-0.5] {\tiny $y$} (dot1);
\draw[K,purple] (dot1) to (dot2);
\draw[K,orange] (dot2) to [bend left=+50] (root);
\end{tikzpicture} 
\begin{tikzpicture}[scale=0.3,baseline=0cm]
	\node at (0,0)  [dot] (root) {};
	\node at (-1,2)  [dot] (dot1) {};
	\node at (1,2)  [dot] (dot2) {};	
	\draw[DK] (root) to node[labl, pos=1.24] {\tiny $h_1$} (dot1);
	\draw[K] (root) to node[labl, pos=1.24] {\tiny $h_2$} (dot2);
	\draw to node[labl, pos=-0.5] {\tiny $y$} (root) ;
\end{tikzpicture}
	\end{equation*}
Repeating the operation, one has in the end
\begin{equation*}
	\begin{tikzpicture}[scale=0.3,baseline=0cm]
		\node at (0,0)  [dot] (root) {};
		\node at (-1,1)  [dot] (dot1) {};
		\node at (1,1)  [dot] (dot2) {};	
		\node at (-1,3)  [] (dot11) {};
		\node at (1,3)  [] (dot21) {};
		\draw[K] (dot1) to node[labl, pos=1.24] {\tiny $h_1$} (dot11);
		\draw[K] (dot2) to node[labl, pos=1.24] {\tiny $h_2$} (dot21);
		\draw[DK] (root) to node[labl, pos=-0.5] {\tiny $y$} (dot1);
		\draw[DK] (root) to node {}(dot2);
		\draw[K,purple] (dot1) to (dot2);
	\end{tikzpicture} = (\cdots) \, + \begin{tikzpicture}[scale=0.3,baseline=0cm]
	\node at (0,0)  [dot] (root) {};
	\node at (-1,1)  [dot] (dot1) {};
	\node at (1,1)  [dot] (dot2) {};	
	\draw[DK] (root) to node {}(dot1);
	\draw[DK] (root) to node {}(dot2);
	\draw[DK] (root) to node[labl, pos=-0.5] {\tiny $y$} (dot1);
	\draw[K,purple] (dot1) to (dot2);
\end{tikzpicture} 	\begin{tikzpicture}[scale=0.3,baseline=0cm]
	\node at (0,0)  [dot] (root) {};
	\node at (-1,2)  [dot] (dot1) {};
	\node at (1,2)  [dot] (dot2) {};	
	\draw[K] (root) to node[labl, pos=1.24] {\tiny $h_1$} (dot1);
	\draw[K] (root) to node[labl, pos=1.24] {\tiny $h_2$} (dot2);
		\draw to node[labl, pos=-0.5] {\tiny $y$} (root) ;
\end{tikzpicture} + 	\begin{tikzpicture}[scale=0.3,baseline=0cm]
\node at (0,0)  [dot] (root) {};
\node at (-1,1)  [dot] (dot1) {};
\node at (1,1)  [dot] (dot2) {};	
\draw[DK] (root) to node {}(dot1);
\draw[DK] (root) to node {}(dot2);
\draw[DK] (root) to node[labl, pos=-0.5] {\tiny $y$} (dot1);
\draw[K,purple] (dot1) to (dot2);
\draw[K,orange] (dot2) to [bend left=+50] (root);
\end{tikzpicture} 
\begin{tikzpicture}[scale=0.3,baseline=0cm]
	\node at (0,0)  [dot] (root) {};
	\node at (-1,2)  [dot] (dot1) {};
	\node at (1,2)  [dot] (dot2) {};	
	\draw[DK] (root) to node[labl, pos=1.24] {\tiny $h_1$} (dot1);
	\draw[K] (root) to node[labl, pos=1.24] {\tiny $h_2$} (dot2);
	\draw to node[labl, pos=-0.5] {\tiny $y$} (root) ;
\end{tikzpicture} + \begin{tikzpicture}[scale=0.3,baseline=0cm]
\node at (0,0)  [dot] (root) {};
\node at (-1,1)  [dot] (dot1) {};
\node at (1,1)  [dot] (dot2) {};	
\draw[DK] (root) to node {}(dot1);
\draw[DK] (root) to node {}(dot2);
\draw[DK] (root) to node[labl, pos=-0.5] {\tiny $y$} (dot1);
\draw[K,purple] (dot1) to (dot2);
\draw[K,orange] (dot2) to [bend left=+50] (root);
\end{tikzpicture} 
\begin{tikzpicture}[scale=0.3,baseline=0cm]
\node at (0,0)  [dot] (root) {};
\node at (-1,2)  [dot] (dot1) {};
\node at (1,2)  [dot] (dot2) {};	
\draw[K] (root) to node[labl, pos=1.24] {\tiny $h_1$} (dot1);
\draw[DK] (root) to node[labl, pos=1.24] {\tiny $h_2$} (dot2);
\draw to node[labl, pos=-0.5] {\tiny $y$} (root) ;
\end{tikzpicture}
\end{equation*}
The terms $ (\cdots) $ are well-defined with the correct Taylor expansions added. The other terms correspond exactly to the counter-terms that are removed by the BPHZ renormalization. Indeed if we apply a `local product' $ R_{\ell}(y) $ depending possibly on the spacetime point one has
\begin{equation*}
R_{\ell}(y)\big(\mathcal{I}_1(\Xi \tau_1 ) \mathcal{I}_1(\Xi \tau_2)\big) =  \mathcal{I}_1(\Xi \tau_1 ) \mathcal{I}_1(\Xi \tau_2) 
+ 
\ell(
	\begin{tikzpicture}[scale=0.3,baseline=0cm]
		\node at (0,0)  [dot] (root) {};
		\node at (-0.7,0.7)  [noise] (dot1) {};
		\node at (0.7,0.7)  [noise] (dot2) {};	
		\draw[DK] (root) to node {}(dot1);
		\draw[DK] (root) to node {}(dot2);
		\draw[DK] (root) to node {} (dot1);
	\end{tikzpicture} \,
	)(y)
\tau_1 \tau_2 
+ 
\ell(
\begin{tikzpicture}[scale=0.3,baseline=0cm]	
	\node at (0,0)  [dot] (0) {};	
	\node at (0.8,0.8)  [dot] (noise1) {};
	\node at (-0.8,0.8)  [noise] (noise2) {};	
	\node at (0.8,0.8)  [noiseblue] (1) {};
	\draw[DK] (0) to (noise1);
	\draw[DK] (0) to (noise2);
	\end{tikzpicture}
)(y)
\mathcal{D} \left( \tau_1 \tau_2 \right)
\end{equation*}
where we recall that $\mathcal{D}$ is the derivation on the symbols given in \eqref{derivation_symbol}.
Extraction elsewhere in the trees is not possible as $ \tau_1 $ and $ \tau_2 $ may contain only the noises $ \Xi_1 $ and $ \Xi_2 $. One uses Propositions \ref{PropPreservedRecenteringBounds}, \ref{prop:invcond3} to check that the recentering bounds are not modified by the telescoping sum.
\end{Dem}

\medskip

The reader may wonder why we did not directly move the edges toward the base point $z$ in case (d). Actually the purple line corresponding to the mollifier can be understood as being as singular as a Dirac mass at a first approximation and therefore it is a sub-divergence of order $ 0 $. By adding the red line it becomes well-defined and if there is no incoming edge on the node associated with the variable $ z $ one can perform the convolution and say that $ K * \rho^{\eps} * \rho^{\eps} $ has the same behaviour as $ K $.

Let us briefly illustrate the proof above on the first two diagrams given in \eqref{graph_examples}:
\begin{equation} \label{graph_examples_bis}
	\begin{tikzpicture}[scale=0.35,baseline=.4cm]
		\node at (0,0)  [root] (root) {};
		\node at (0,1)  [dot] (dot1) {};
		\node at (0,2)  [var] (var) {};
		\node at (0,3)  [dot] (dot2) {};
		\node at (-0.5,4)  [noise] (n3) {};
		\node at (0.5,4)  [noise] (n4) {};
		
		\draw[testfcn] (dot1) to  (root);
		\draw[DK] (dot1) to [bend left=+70] (var);
		\draw[DK] (dot1) to [bend left=-70] (var);
		\draw[K] (var) to [bend left=0] (dot2); 
		\draw[DK] (dot2) to node[labl, pos=1.24] {\tiny $h_1$} (n3);
		\draw[DK] (dot2) to node[labl, pos=1.24] {\tiny $h_2$} (n4);   
	\end{tikzpicture} 
	\hspace{0.3cm}, \hspace{0.3cm}
	\begin{tikzpicture}[scale=0.3,baseline=.4cm]
		\node at (0,0)  [root] (root) {};
		\node at (0,1)  [dot] (dot1) {};
		\node at (-1,2)  [noise] (n1) {};
		\node at (1,2.3)  [var] (var) {};
		\node at (-1,3.5)  [dot] (dot2) {};
		\node at (-1.5,4.5)  [noise] (n3) {};
		
		\draw[testfcn] (dot1) to (root);
		\draw[DK] (dot1) to (n1);
		\draw[DK] (dot1) to (var);
		\draw[K] (n1) to node[labl, pos=0.05] {\tiny $h_2$} (dot2); 
		\draw[DK] (dot2) to node[labl, pos=1.24] {\tiny $h_1$} (n3);
		\draw[DK] (dot2) to (var);   
	\end{tikzpicture}.
\end{equation}
The second diagram has no divergence thus the integrability condition~\eqref{condition_integrability} is satisfied. The first diagram has the same type of divergence at the root as in the illustrative example of the proof. One can perform a similar renormalisation as described above for getting the convergence of this term.

\medskip

Given the results of Corollary \ref{PropMirrorDiagrams} and Proposition \ref{PropBoundDerivatives} the following fact follows as a consequence of Stroock's formula \eqref{EqStroockFormula}

$$
\big\Vert \big\langle \overline{\sf \Pi}^\epsilon_z\tau , \varphi_z^\lambda\big\rangle \big\Vert_{L^2(\Omega)} \leq \sum_{k=0}^{\vert\tau\vert_\zeta-1} \big\Vert \bbE\big[\big\langle d^k \overline{\sf \Pi}^\epsilon_z\tau , \varphi_z^\lambda\big\rangle\big] \big\Vert + \big\Vert \big\langle d^{\vert\tau\vert_\zeta} \overline{\sf \Pi}^\epsilon_z\tau , \varphi_z^\lambda\big\rangle \big\Vert.
$$

\medskip

\begin{cor}
One has for every $\tau\in \mcB^-$ the $\epsilon$-uniform estimates
\begin{equation} \label{EqEpsilonUniformBound}
\big\Vert \big\langle \overline{\sf \Pi}^\epsilon_z\tau , \varphi_z^\lambda\big\rangle \big\Vert_{L^2(\Omega)} \lesssim \lambda^{\vert\tau\vert}   \qquad (\forall\,0<\lambda\leq 1).
\end{equation}
\end{cor}

\medskip

We know from Theorem 10.7 of \cite{Hai14} that the convergence of the BHZ model to a limit admissible model follows from some estimates of the form
\begin{equation} \label{EqFinalEstimate}
\big\Vert \big\langle \overline{\sf \Pi}_{z}^\epsilon \tau - \overline{\sf \Pi}_{z}^{\epsilon'}\hspace{-0.07cm}\tau \,,\,\varphi_{z}^\lambda \big\rangle \big\Vert_{L^2(\Omega)} \lesssim o_{\epsilon\vee\epsilon'}(1)\,\lambda^{\vert\tau\vert},
\end{equation}
for $\tau\in\mcB^-$. To prove these bounds from the previous bounds we note that the scaling bound \eqref{EqScalingBoundExpectation} from Proposition \ref{PropBoundDerivatives} holds if the kernel $K$ behaved not as $\vert z\vert_{\frak{s}}^{-1}$, as $z$ goes to $0$, but as $\vert z\vert_{\frak{s}}^{-1-\eta}$, for $\eta>0$ small enough. Recall from Proposition 10.17 of \cite{Hai14} the bound
\begin{equation} \label{EqEpsilonGain}
|\!|\!| \rho^\epsilon*\rho^\epsilon - \rho^{\epsilon'}*\rho^{\epsilon'} |\!|\!|_{3+\eta,1} \lesssim \max(\epsilon,\epsilon')^\eta.
\end{equation}
As in Section 5.2 of Hairer \& Pardoux work \cite{HP} one can decompose $\overline{\sf \Pi}_{z}^\epsilon \tau - \overline{\sf \Pi}_{z}^{\epsilon'}$ as a telescoping sum where we replace one after the other the $\zeta^\epsilon$ by the $\zeta^{\epsilon'}$. This amounts to replacing some convolution operators with $\rho^\epsilon$ by the corresponding convolution operator with $\rho^{\epsilon'}$, so the corresponding term in the telescoping sum involves the difference of two terms that differ at only one place. We gain the $o_{\epsilon\vee\epsilon'}(1)$ term in \eqref{EqFinalEstimate} from using the estimate \eqref{EqEpsilonGain} for this edge while repeating the proof of \eqref{EqEpsilonUniformBound}. We leave the details to the reader.

\medskip

One can adapt the present work to prove the convergence of the BHZ model associated with a system of generalized (KPZ) equations, such as the system describing the random motion of a rubber on a manifold \cite{BGHZ}. The diagrams to be considered are essentially the same as we have more noises $\xi_i$ that are some independent copies of the same space-time white noise. The derivative in space will be the same as we increase the number of equations but our solution is still defined on the one dimenstional torus. It is just that the counter-terms for renormalizing the system are more complex to write as they have to take into account directional derivatives.

\bigskip

\noindent \textcolor{gray}{$\bullet$} {\sf I. Bailleul} -- Univ Brest, CNRS, LMBA - UMR 6205, F-29238 Brest, France.   \\
\noindent {\it E-mail}: ismael.bailleul@univ-brest.fr   

\medskip

\noindent \textcolor{gray}{$\bullet$} {\sf Y. Bruned} --   Universite de Lorraine, CNRS, IECL, F-54000 Nancy, France.   \\
{\it E-mail}: yvain.bruned@univ-lorraine.fr


\bigskip



\begin{thebibliography}{99}

\bibitem{BailleulBruned}
I. Bailleul and Y. Bruned,
\newblock {\em Locality for singular stochastic PDEs}.
\newblock arXiv:2109.00399, To appear in Ann. Probab. (2026).

\bibitem{RSGuide}
I. Bailleul and M. Hoshino,
\newblock {\em A tourist guide to regularity structures and singular stochastic PDEs}.
EMS Surv. Math. Sci. (2025).

\bibitem{BH23}
I. Bailleul and M. Hoshino,
\newblock {\em Random models on regularity-integrability structures}.
\newblock arXiv:2310.10202.

\bibitem{Berglund}
N. Berglund, 
\newblock {\em An Introduction to Singular Stochastic PDEs. Allen-Cahn Equations, Metastability, and Regularity Structures}.
\newblock EMS Series of Lectures in Mathematics, (2022).

\bibitem{BOT24}
L. Broux, F. Otto and M. Tempelmayr,
\newblock {\em Lecture notes on Malliavin calculus in regularity structures}.
Stoch. PDEs: Anal. Comp., (2025).

\bibitem{bruned:tel-01306427}
Y.~Bruned.
\newblock \emph{{Singular KPZ Type Equations}}.
\newblock 205 pages, {PhD thesis}, {Universit{\'e} Pierre et Marie Curie - Paris VI}, (2015).
\newblock \url{https://tel.archives-ouvertes.fr/tel-01306427}.

\bibitem{BrunedRecursive}
Y. Bruned,
\newblock {\em Recursive formulae for regularity structures}.
\newblock Stoch. PDEs: Anal. Comp., {\bf 6}(4):525--564, (2018).

\bibitem{BCCH}
Y. Bruned and A. Chandra and I. Chevyrev and M. Hairer,
\newblock {\em Renormalising SPDEs in regularity structures}.
\newblock J. Europ. Math. Soc., {\bf 23}(3):869--947, (2021).

\bibitem{BGHZ}
Y. Bruned and F. Gabriel and M. Hairer and L. Zambotti,
\newblock {\em Geometric stochastic heat equations}.
\newblock J. Am. Math. Soc., {\bf 35}(1):1--80, (2022).

\bibitem{BHZ}
Y. Bruned and M. Hairer, and L. Zambotti,
\newblock {\em Algebraic renormalization of regularity structures},
\newblock Invent. Math., {\bf 215}(3):1039--1156, (2019).

\bibitem{BN1}
Y. Bruned and U. Nadeem,
\newblock {\em Convergence of space-discretised gKPZ via Regularity Structures}. Ann. Appl. Probab. \textbf{34}(2), 2488--2538, (2024).

\bibitem{BrunedNadeem}
Y. Bruned and U. Nadeem,
\newblock {\em Diagram-free approach for convergence of tree-based models in regularity structures}.
J. Math. Soc. Japan, 1--31, (2024).

\bibitem{CH}
A. Chandra and M. Hairer,
\newblock {\em An analytic BPHZ theorem for regularity structures}.
\newblock arXiv:1612.08138, (2016).

\bibitem{ChandraWeber}
A. Chandra and H. Weber,
\newblock {\em Stochastic PDEs, regularity structures and interacting particle systems}.
\newblock Ann. Fac. Sci. Toulouse, {\bf 26}(4):847--909, (2017).

\bibitem{CorwinShen}
I. Corwin and H. Shen,
\newblock {\em Some recent progress in singular stochastic PDEs}.
\newblock Bull. Amer. Math. Soc., {\bf 57}(3):409--454, (2020).

\bibitem{FrizHairer}
P. Friz and M. Hairer,
\newblock {\em A Course on Rough Paths, With an introduction to regularity structures}.
\newblock Universitext, Springer, (2020).

\bibitem{GIP}
M. Gubinelli and P. Imkeller and N. Perkoswki,
\newblock {\em Paracontrolled distributions and singular PDEs}.
\newblock Forum Math. Pi, {\bf 3}(e6):1--75, (2015).

\bibitem{Hai14}
M. Hairer,
\newblock {\em A theory of regularity structures}.
\newblock Invent. Math., {\bf 198}:269--504, (2014).

\bibitem{HairerTakagi}
M. Hairer,
\newblock {\em Renormalization of parabolic stochastic PDEs}.
\newblock The 20th Takagi Lectures, arXiv:1803.03044, (2018).

\bibitem{HP}
M. Hairer and E. Pardoux.
\newblock {\em A Wong-Zakai theorem for stochastic PDEs}.
\newblock J. Math. Soc. Japan, {\bf 67}(4):1551--1604, (2015).

\bibitem{HairerQuastel}
M. Hairer and J. Quastel.
\newblock {\em A class of growth models rescaling to KPZ}.
\newblock Forum Math. Pi, {\bf 6}(e3):1--112, (2018).

\bibitem{HS23}
M. Hairer and R. Steele.
\newblock {\em The BPHZ Theorem for Regularity Structures via the Spectral Gap Inequality}.
Arch. Rational Mech. Anal. \textbf{248}(9), 1--81, (2024).


\bibitem{LOTT}
P. Linares and F. Otto and M. Tempelmayr and P. Tsatsoulis,
\newblock{\em A diagram-free approach to the stochastic estimates in regularity structures}. 
Invent. math., {\bf 237}:1469--1565, (2024).


\bibitem{MourratWeberXu}
J.C. Mourrat and H. Weber and W. Xu,
\newblock {\em Construction of $\Phi^4_3$ diagrams for pedestrians}.
\newblock From particle systems to partial differential equations, 1--46, Springer Proc. Math. Stat., 209, Springer, Cham, (2017).

\bibitem{NourdinPeccati}
I. Nourdin and G. Peccati,
\newblock {\em Normal approximations with Malliavin calculus, from Stein's method to universality}.
\newblock Cambridge tracts in Mathematics, {\bf 192}, (2012).

\bibitem{OSSW2}
F. Otto and J. Sauer and S. Smith and H. Weber,
\newblock {\em A priori bounds for quasi-linear SPDEs in the full sub-critical regime}.
 J. Eur. Math. Soc. (JEMS) {\bf 27}(1):71--118, (2025). 

\bibitem{OST}
F. Otto and K. Seong and M. Tempelmayr,
\newblock {\em Lecture notes on tree-free regularity structures}. Mat. Contemp.,
{\bf 58}:150--196, (2023). 


\bibitem{Wein} 
S. Weinberg. 
\newblock {\em High-energy behavior in quantum field-theory}. 
\newblock Phys. Rev., {\bf 118}(2):838--849, (1960).

\end{thebibliography}
\end{document}